\documentclass[leqno]{article}
\usepackage{latexsym}
\begin{document}
\addtolength{\baselineskip}{0.1cm}
\title{{\huge Bicomplex algebra and function theory}}
\author{{\Large  Stefan R\"{o}nn}\\[4mm]
        Helsinki University of Technology\\
        Department of Computer Science and Engineering\\
        PBox 5400, FIN-02015 HUT, Finland\\[2mm]
        Email: Stefan.Ronn@hut.fi\\
        }
\date{\today}
\maketitle
\thispagestyle{empty}
\begin{abstract}
{\normalsize This treatise investigates holomorphic functions defined on
the space of bicomplex numbers introduced by Segre.
The theory of these functions is associated with
Fueter's theory of regular, quaternionic functions.
The algebras of quaternions and bicomplex
numbers are developed by making use of so-called complex pairs.
Special attention is paid to singular bicomplex numbers that
lack an inverse. The elementary bicomplex functions
are defined and their properties studied. 
The derivative of a bicomplex function is defined as the limit of a
fraction with nonsingular denominator. 
The existence of the derivative amounts to the validity of
the complexified Cauchy-Riemann equations,
which characterize the holomorphic bicomplex functions.
It is proved that such a function has derivatives of all orders.
The bicomplex integral is defined as a line integral.
The condition for path independence and the bicomplex
generalizations of Cauchy's theorem and integral formula are
given. Finally, the relationship between the
bicomplex functions and different forms of the Laplace equation is
considered. In particular, the four-dimensional Laplace equation
is factorized using quaternionic differential operators. The outcome is new classes of
bicomplex functions including Fueter's regular functions.
It is shown that each class contains differentiable functions.}
\end{abstract}
\newpage
\pagenumbering{roman}
\tableofcontents
\newpage
\pagenumbering{arabic}

\newcommand{\tvasteg}{\hspace*{0.3cm}}
\newcommand{\tresteg}{\hspace*{0.6cm}}
\newcommand{\femsteg}{\hspace*{1.0cm}}
\newcommand{\smallcup}{\scriptscriptstyle \cup}
\newcommand{\smallbox}{\scriptscriptstyle \Box}
\newcommand{\smalldiamond}{\scriptscriptstyle \Diamond}
\newcommand{\smallbot}{\scriptscriptstyle \bot}

\section{Introduction}
A common technique in mathematics is to define a new object as an aggregate
consisting of previously defined objects, usually of the same
type. A well-known example is to regard a space vector as
consisting of three real numbers. The number concept is also
extended with this technique. Thus, a rational number is viewed as
a pair of integers and a complex number as a pair of real numbers.
It is only natural to ask whether one can carry
this idea a step further by defining numbers that are represented by
pairs of complex numbers, or {\em complex pairs} as we shall call them.
It turns out that both the {\em quaternions} and the {\em bicomplex numbers}
can be conceived this way. The former were discovered by William Rowan
Hamilton in 1843, the latter by Corrado Segre in 1892~\cite{hamilton2, segre}.
Both types of numbers can be seen as generalizations of the complex numbers
and they lead to corresponding generalizations of the complex functions.

The development of the theory of {\em quaternionic functions} was initiated
by Rudolf Fueter in 1935~\cite{fueter}. He proposed
that in this theory the objects of interest
ought to be the so-called {\em regular functions}, which
are defined by means of a close analogue of the Cauchy-Riemann
equations. More specifically, with $q=x+yi+zj+uk$
a quaternion written in standard notation and  
$\psi\! : q \rightarrow \psi(q)$ a quaternionic function,
$\psi$ is regular if it satisfies the partial differential equation
\begin{equation}
\label{CRF-equations}
\frac{\partial\psi}{\partial x} + i\times\frac{\partial\psi}{\partial y} + j\times\frac{\partial\psi}{\partial z} + k\times\frac{\partial\psi}{\partial u} \;=\; 0
\end{equation}
where $\times$ denotes quaternionic multiplication. In the complex
plane the corresponding formula is
\[
\frac{\partial f}{\partial x} + i\cdot \frac{\partial f}{\partial y} = 0
\]
which is an equivalent way of writing the Cauchy-Riemann
equations and characterizes a complex holomorphic (analytic) function
$f\!: a \rightarrow f(a)$, where $a=x+iy$.

In the other direction there is the theory of {\em bicomplex functions}.
It has fairly recently been investigated by G.Baley Price~\cite{price}, whose
book also contains a brief history of the development of the main
ideas behind the theory (including an account of Segre's work). 
The bicomplex functions of interest are the
{\em holomorphic} ones, which are characterized by the fact that
they are differentiable. They are almost isomorphic to the complex
holomorphic functions, not surprisingly, because the
operations of the bicomplex algebra are almost isomorphic to those
of the complex algebra. The bicomplex algebra has an anomaly, however,
which has a bearing on the associated function theory: there exist so-called
singular numbers other than 0 that lack an inverse.

The purpose of this treatise is to give an overview of and further
develop the theory of holomorphic bicomplex functions. We shall among others
demonstrate that these functions can be divided into eight, slightly different
classes. Some of the functions satisfying (\ref{CRF-equations}) belong to them, hence Fueter's theory is
related to our theory.

Our purpose requires that we study the algebras of the bicomplex
numbers and the quaternions first. Although both number types
have several representations, we shall mostly represent them
by complex pairs, such a pair being of
the form $(a,b)$, where $a$ and $b$ are complex numbers.
Complex pairs will in fact be used throughout the treatise,
because in general they greatly facilitate calculations.

In Chapter 3 we introduce the bicomplex functions that we
regard as elementary: the exponential and logarithm function, 
the hyperbolic and trigonometric functions, 
polynomials and the quotient function.
The properties of the exponential and logarithm function are
specially important, since they throw light on the role played
by the singular numbers in bicomplex function theory.
We shall show that the exponential function gets only
nonsingular values and that, as a consequence, the logarithm
function is defined only for nonsingular arguments.
The periods that we obtain for the transcendental functions
in this chapter are also of interest.

In Chapter 4 we consider the differentiation of bicomplex functions.
The derivative is defined in the usual way as the limit of a
fraction, but the limit must be taken for nonsingular values of
the denominator only. Rendering the derivative as a complex
pair yields two alternative formulas, the equivalence of which amounts
to the validity of the complexified Cauchy-Riemann equations.
Fulfillment of these equations is the main characteristic of
the holomorphic bicomplex functions. We prove that such a
function possesses derivatives of all orders.
In addition, we give ${\rm R^{4}}$-representations
of the bicomplex derivative.

In Chapter 5 we consider the integration of bicomplex
functions. The basic bicomplex integral is defined as
a line integral, which is path
independent if and only if the integrand is the derivative of a differentiable 
bicomplex function. We show that a holomorphic bicomplex function
can be expanded into a Taylor series. In order to prove the bicomplex versions of
Cauchy's theorem and integral formula we generalize Green's theorem
and derive the so-called twining number, the bicomplex analogue of
the complex winding number.

In the last chapter we clarify the relationship between 
bicomplex functions and the Laplace equation. For this purpose we 
first study the derivatives of complex holomorphic and conjugate 
holomorphic functions. We factorize the four-dimensional 
Laplace equation using quaternionic differential operators and 
concepts obtained from the similar but simpler factorization of 
the two-dimensional version of the equation. This leads
to the aforementioned classes of bicomplex functions,
each class containing also differentiable functions.
In this context we shall have the occasion to define
the derivative related to Fueter's regular functions, a topic
that has received attention in the literature. Our proposal is
to define the derivative as the limit of a bicomplex fraction where
both the numerator and the denominator are regular in Fueter's sense.

Theorems, lemmas, definitions and formulas are numbered separately
within each chapter. The number of a theorem, lemma or definition
is always preceded by the chapter number. Numbered formulas
are referred to by the assigned number when needed in the same chapter,
otherwise the chapter number is added.

\newpage
\setcounter{equation}{0}

\section{Quaternionic and bicomplex algebra}
In this chapter we shall make
extensive use of complex pairs to develop the algebras of quaternions
and bicomplex numbers. Although both algebras generalize
the complex algebra, we can of course not expect them to
retain all the properties of this algebra:
in the quaternionic algebra the commutativity of the multiplication operation
is lost, and in the bicomplex algebra we must accept the existence of
the singular numbers. The bicomplex numbers are more 
important for our work, but the
quaternions are not unimportant. They are needed in the
factorization of the four-dimensional Laplace equation and,
above all, provide us with a norm that directly can be taken over to
the bicomplex algebra. Thus, to begin
with we shall investigate the quaternions.

\subsection{Representation of quaternions}
A quaternion is in general written as a four-dimensional vector of
the form
\[
x+yi+zj+uk
\]
where $x,y,z,u$ are real numbers and $i,j,k$ unit
vectors (generalized imaginary units). Given two quaternions
\begin{eqnarray*}
 q &=& x+yi+zj+uk \\ [1mm]
 r &=& X+Yi+Zj+Uk
\end{eqnarray*}
of this kind we may add them by applying component-wise addition, formally:
\begin{equation}
\label{quataddr4}
 q+r \;=\; (x+X)+(y+Y)i+(z+Z)j+(u+U)k
\end{equation}
To denote quaternionic multiplication we shall
consistently use the symbol $\times$.
The well-known multiplication table for the unit vectors is
\addtocounter{equation}{1}
\begin{tabbing}
 XXXXXXXXX \= XXXXXXXXXXXX\=XXXXXXXXXXXX\=XXXXXXXXXXXXXXXXXXXXXXXXXXXXXXXXXXXXXXXXXX \kill
 (\theequation a)\> $i \times i=-1$ \> $j \times j=-1$ \> $k \times k=-1$ \\
 (\theequation b)\> $i \times j=k$ \> $j \times k=i$ \> $k \times i=j$ \\
 (\theequation c)\> $k \times j=-i$ \> $j \times i=-k$ \> $i \times k=-j$
\end{tabbing}
which yields the following expression for the non-commutative product $q\times r$:
\begin{eqnarray}
\label{quatmultr4}
 q \times r & = & (xX-yY-zZ-uU)\: + \:(xY+yX+zU-uZ)i\:+ \\ 
 & & (xZ-yU+zX+uY)j\: + \:(xU+yZ-zY+uX)k \nonumber
\end{eqnarray}
We now wish to represent the quaternions $q$ and $r$ as pairs of
complex numbers:
 \[
 \begin{array}{rcl}
 q &\!\! =\!\! & (a,b) \mbox{\femsteg , } a=x+yi \mbox{ and } b=z+ui \\[2mm]
 r &\!\! =\!\! & (c,d) \mbox{\femsteg , } c=X+Yi \mbox{ and } d=Z+Ui
 \end{array}
 \]
How should the multiplication $\times$ be expressed in
terms of the complex numbers $a,b,c,d$? We find
\begin{equation}
\label{quatmultc2}
 q \times r = (a \cdot c - b \cdot d^{*}\,,\, b \cdot c^{*} + a \cdot d)
\end{equation}
where the operators $\cdot$ and * stand for
complex multiplication and conjugation.
Writing the right-hand side of this formula in terms of $x,y,z,u$
and $X,Y,Z,U$ yields
\begin{eqnarray}
\label{quatmultr4cp}
 q \times r & = & (\, xX-yY-zZ-uU+(xY+yX+zU-uZ)i\:, \\
& & \;\;xZ-yU+zX+uY + (xU+yZ-zY+uX)i\,) \nonumber
\end{eqnarray}
which is another way of rendering (\ref{quatmultr4}).

It is sometimes convenient to represent the
quaternion $x+yi+zj+uk$ as a quadruple
\[ (x,y,z,u) \]
With $q=(x,y,z,u)$ and $r=(X,Y,Z,U)$, the addition and
multiplication rules then look like
\begin{eqnarray*}
 q+r &=& (x+X,y+Y,z+Z,u+U)\\ [2mm]
 q \times r &=& (\,xX-yY-zZ-uU\,,\, xY+yX+zU-uZ\,, \\ 
 & & \;\;xZ-yU+zX+uY\,,\,xU+yZ-zY+uX\,) \nonumber
\end{eqnarray*}
In the context of quaternionic algebra
the four-dimensional real space ${\rm R^{4}}$ is usually called the
quaternionic or Hamiltonian space ${\rm H}$. It can also be
regarded as the Cartesian product of two complex planes, or ${\rm C^{2}}$.
The different views of the space ${\rm R^{4}}$
enable us to classify the representations we have introduced
for a quaternion $q=(a,b)=(x+yi,z+ui)$:
\begin{equation}
\label{quatrepr}
\begin{array}{ll}
\mbox{${\rm H}$-representation:} &q \\ [1mm]
\mbox{${\rm C^{2}}$-representation, complex pair:} &(a,b) \\ [1mm]
\mbox{${\rm R^{4}}$-representation, complex pair:} &(x+yi,z+ui) \\ [1mm]
\mbox{${\rm R^{4}}$-representation, vector-form:} &x+yi+zj+uk  \\ [1mm]
\mbox{${\rm R^{4}}$-representation, quadruple-form:} &(x,y,z,u)
\end{array}
\end{equation}
As a rule the ${\rm C^{2}}$-representation
leads to the shortest calculations.

\subsection{Basic quaternionic algebra with complex pairs}
Let $(a,b)$ and $(c,d)$ be two complex pairs
representing quaternions. The pairs are equal if and only if
their components are equal, formally:
\begin{equation}
\label{quatequal}
 (a,b)=(c,d)\;\: \equiv \;\: (a=c) \wedge (b=d)
\end{equation}
{\bf Remark.} The logical operators are designated by
$\equiv$ (equivalence), $\Rightarrow$ (implication), $\vee$
(disjunction), $\wedge$ (conjunction), and $\neg$ (negation).
Their priority order is from weaker to stronger: $\equiv$,
$\Rightarrow$, $\vee$ and $\wedge$, $\neg$. $\vee$ and $\wedge$
have the same priority. $\Box$

\medskip

\noindent For quaternionic addition and multiplication we found in
the previous section
\begin{eqnarray}
\label{quataddc2}
(a,b)+(c,d)&=&(a+c,b+d) \\ [1mm]
\label{quatmultc2second}
(a,b)\times (c,d) &=& (ac - bd^{*}, bc^{*} + ad)
\end{eqnarray}
With these ${\rm C^{2}}$-formulas it is not difficult to verify the following basic
properties of $+$ and $\times$ (using one of the
${\rm R^{4}}$-representations for the quaternions the
proofs of the first three properties are rather laborious):
\begin{eqnarray*}
[(a,b)\times (c,d)]\times (e,f) &=& (a,b)\times [(c,d)\times (e,f)] \;\mbox{\femsteg \tresteg $\lbrace$$\times$ is associative$\rbrace$} \\ [1mm]
(a,b)\times [(c,d)+(e,f)] &=& (a,b)\times (c,d)\,+\, (a,b)\times (e,f) \mbox{\tresteg $\lbrace$left-distribution$\rbrace$} \\ [1mm]
{[}(c,d)+(e,f)] \times (a,b) &=& (c,d)\times (a,b)\,+\, (e,f)\times (a,b) \mbox{\tresteg $\lbrace$right-distribution$\rbrace$}
\end{eqnarray*}
\begin{eqnarray*}
(a,b)\times (1,0) &=& (a,b) \mbox{\femsteg \tresteg $\lbrace$identity element of $\times$$\rbrace$} \\ [1mm]
(a,b)+ (0,0) &=& (a,b) \mbox{\femsteg \tresteg $\lbrace$identity element of $+$$\rbrace$} \\ [1mm]
(a,b)\times (0,0) &=& (0,0) \mbox{\femsteg \tresteg $\lbrace$zero element of $\times$$\rbrace$}
\end{eqnarray*}
The complex pairs
$(1,0)$ and $(0,0)$ are examples of {\em real} or {\em scalar}
quaternions, whose general form is $(\lambda,0)$, where $\lambda \in {\rm R}$.
Multiplication of a quaternion by a real quaternion is in a sense a
distributive operation, since for any $(a,b) \in {\rm C^{2}}$
\begin{equation}
\label{realxquat}
(\lambda,0)\times (a,b) = (\lambda a,\lambda b)
\end{equation}
Moreover, it is also a commutative operation or
\[
(\lambda,0)\times (a,b) = (a,b)\times (\lambda,0)
\]
The quaternion $(\lambda,0)$ is designated by just
$\lambda$, because it has all the properties of the real number $\lambda$.
The $\times$-operator can be omitted from a quaternionic product if one
of the factors is a real quaternion. For example, the previous
formula is written $\lambda (a,b) = (a,b)\lambda$.

\subsection{Quaternionic conjugation, square norm, and division}
In this section we shall introduce some concepts based on
the quaternionic conjugate. 
The terminology and notation are partly taken from~\cite{dieudonne},~\cite{gurlebeck} and \cite{porteous}.
We first recall the properties
of the complex conjugate. With $a$ and $b$
two complex numbers the conjugation operator * satisfies:
\begin{eqnarray}
\label{cstaroverplus}
(a+b)^{*}&=&a^{*}+b^{*} \mbox{\femsteg $\lbrace$* distributes over +$\rbrace$} \\
\label{cstarovermult}
(a\cdot b)^{*}&=&a^{*}\cdot b^{*} \mbox{\femsteg $\;\;\lbrace$* distributes over $\cdot \rbrace$} \\
\label{cstaridempot}
(a^{*})^{*}&=&a \mbox{\femsteg \tresteg $\;\;\;\lbrace$* is an involution$\rbrace$}
\end{eqnarray}
We shall use the symbol * for quaternionic conjugation, too.
Its definition in ${\rm C^{2}}$ is
\begin{equation}
\label{qstardef}
(a,b)^{*}=(a^{*},-b)
\end{equation}
The equivalent ${\rm R^{4}}$-formulations are:
\begin{eqnarray*}
(x+yi,z+ui)^{*}&=&(x-yi,-z-ui)\\
(x+yi+zj+uk)^{*}&=&x-yi-zj-uk \\
(x,y,z,u)^{*}&=&(x,-y,-z,-u)
\end{eqnarray*}
Quaternionic conjugation is an involution and distributes
over addition, but not over multiplication. Hence, of the formulas
(\ref{cstaroverplus}), (\ref{cstarovermult}), and (\ref{cstaridempot}) only the first and last may be taken over to
quaternionic space:
\begin{eqnarray}
\label{qstaroverplus}
(q+r)^{*}&\!=\!&q^{*}+r^{*} \mbox{\femsteg $q,r \in {\rm H}$}\\ [1mm]
\label{qstaridempot}
(q^{*})^{*}&\!=\!&q
\end{eqnarray}
The {\em square norm} $N(q)$ of the quaternion $q$
is defined by
\begin{equation}
\label{squarenormh}
N(q)=q\times q^{*} \mbox{\femsteg $\lbrace$${\rm H}$-representation$\rbrace$}
\end{equation}
With $q=(a,b)$ we obtain the equivalent complex pair formula,
which reveals that the square norm is a real number:
\begin{equation}
\label{squarenormc2}
N((a,b))= aa^{*}+bb^{*} \mbox{\femsteg $\lbrace$${\rm C^{2}}$-representation$\rbrace$}
\end{equation}
In ${\rm R^{4}}$ the same formula looks like
\begin{equation}
\label{squarenormr4}
N((x+yi,z+ui))= x^{2}+y^{2}+z^{2}+u^{2} \mbox{\femsteg $\lbrace$${\rm R^{4}}$-representation$\rbrace$}
\end{equation}
If we let $N(a)$ stand for the square norm of a complex number $a$, too,
we have $N(a)=aa^{*}$, and hence, by (\ref{squarenormc2})
\begin{equation}
\label{Ndecomp}
N((a,b))=N(a)+N(b) 
\end{equation}
For two quaternions $q$ and $r$ we also have
\[
N(q\times r) = N(q)N(r) 
\]
The function $N$ enables us to introduce
further quaternionic concepts.
We denote the {\em inverse} of a nonzero
quaternion $q$ by $q\!\uparrow\! -1$
(the standard notation $q^{-1}$ being reserved for another
concept), and define
\begin{equation}
\label{qinverse}
q\!\uparrow\! -1 = \frac{1}{N(q)}q^{*} \mbox{\femsteg , $q\neq 0$}
\end{equation}
With the inverse at our disposal
we can introduce a division operation for two quaternions $r$ and $q$, but due to
the non-commutativity of $\times$ we must distinquish
between left-division, $(q\!\uparrow\! -1)\times r$,
and right-division, $r \times (q\!\uparrow\! -1)$. 

Finally, we define the notion of the
{\em absolute value} or {\em modulus}
$\| q \|$ of $q$ by
\begin{equation}
\label{normh}
\| q \| = \sqrt{N(q)} \mbox{\femsteg , $q\in {\rm H}$}
\end{equation}
The $\| \:\: \|$-function plays an important role in the
bicomplex algebra, too, as we shall soon see.

\subsection{Bicomplex algebra}
The space of {\em bicomplex numbers}, which we denote
by ${\rm B}$, has many properties in common with the
Hamiltonian space. A bicomplex number $q$ is foremost viewed
as a complex pair $(a,b)$, where $a,b\in {\rm C}$, but all the other
quaternion representations (\ref{quatrepr}) apply to it, too:
\begin{equation}
\label{bicrepr}
\begin{array}{ll}
\mbox{${\rm B}$-representation:} &q \\ [1mm]
\mbox{${\rm C^{2}}$-representation, complex pair:} &(a,b) \\ [1mm]
\mbox{${\rm R^{4}}$-representation, complex pair:} &(x+yi,z+ui) \\ [1mm]
\mbox{${\rm R^{4}}$-representation, vector-form:} &x+yi+zj+uk  \\ [1mm]
\mbox{${\rm R^{4}}$-representation, quadruple-form:} &(x,y,z,u)
\end{array}
\end{equation}
The equality condition of two bicomplex numbers is given by
(\ref{quatequal}) and they are added according to (\ref{quataddc2}).
The multiplication of bicomplex numbers, however, differs
from the quaternionic operation. 
If we write a complex number $x+iy$
as a pair of {\em real} numbers $(x,y)$, ordinary
complex multiplication of $a=(x,y)$ by $b=(z,u)$
takes the form
\begin{displaymath}
a \cdot b=(xz-yu,yz+xu)
\end{displaymath}
Replacement of the factors $a$ and $b$ in this scheme
by bicomplex numbers yields bicomplex multiplication, designated by $\odot$.
Thus, for $q=(a,b)$ and $r=(c,d)$ we have the definition
\begin{equation}
\label{bicmultdef}
q \odot r = (a \cdot c-b \cdot d\,,\, b \cdot c+a \cdot d)
\end{equation}
Evaluating this formula for
$q=(x+iy,z+iu)$ and $r=(X+iY,Z+iU)$, we obtain the equivalent
${\rm R^{4}}$-formula
\begin{eqnarray*}
 q \odot r &=& (xX-yY-zZ+uU\,+\, i(xY+yX-zU-uZ)\,, \\
 & & \;\,xZ-yU+zX-uY\,+\,i(xU+yZ+zY+uX))
\end{eqnarray*}
The vector-form $q=x+yi+zj+uk$ of a bicomplex number
is occasionally useful. The unit vectors have the complex
pair representations
\begin{equation}
\label{unitvectors}
i=(i,0) \femsteg , \femsteg j=(0,1) \femsteg , \femsteg k=(0,i)
\end{equation}
and are multiplied according to:
\begin{tabbing}
 XXXXXXXXX \= XXXXXXXXXXXX\=XXXXXXXXXXXX\=XXXXXXXXXXXXXXXXXXXXXXXXXXXXXXXXXXXXXXXXXX \kill
 \> $i \odot i=-1$ \> $j \odot j=-1$ \> $k \odot k=1$ \\
 \> $i \odot j=k$ \> $j \odot k=-i$ \> $k \odot i=-j$ \\
 \> $j \odot i=k$ \> $k \odot j=-i$ \> $i \odot k=-j$
\end{tabbing}
A {\em bicomplex scalar} $(\lambda,0)$, where $\lambda \in {\rm R}$,
will be abbreviated by just $\lambda$ in the same way
as a real quaternion. The identity element
$(1,0)$ and the zero element $(0,0)$ of $\odot$ are therefore written 1 and 0.
We also abbreviate the expression $\lambda\odot(a,b)$ by $\lambda (a,b)$
and note that 
\begin{equation}
\label{scalartimesbic}
\lambda (a,b)=(\lambda a, \lambda b)
\end{equation}
The bicomplex number $(c,0)$, where $c\in {\rm C}$, has all the
properties of the complex number $c$, yet we shall not use any
abbreviation in this case.

\subsection{Further bicomplex operations}
The $\odot$-operator is associative, distributes
over + and, most significantly, is {\em commutative}. It therefore
permits us to generalize some of the algebraic
operations of the complex numbers. The exponentiation
of a bicomplex number is performed in the customary fashion:
\begin{equation}
\label{expdef}
\left\{ \begin{array}{ll}
          q^{0}=1 & \\ [1mm]
          q^{n}=q\odot q^{n-1} & \mbox{if $n$ is an integer $\geq 1$}
        \end{array}
\right.
\end{equation}
Bicomplex conjugation is meaningful, too. We designate it by a 
postfix \raisebox{.6ex}{$\smallcup$} and
give it the ${\rm C^{2}}$-definition:
\begin{equation}
\label{biccupdef}
(a,b)^{\smallcup}=(a,-b)
\end{equation}
It is isomorphic to the definition
of complex conjugation * expressed by real pairs or
$(x,y)^{*}=(x,-y)$. As a result, the \raisebox{.6ex}{$\smallcup\,$}-operator inherits
the properties (\ref{cstaroverplus}),(\ref{cstarovermult}), and (\ref{cstaridempot}) of *:
\begin{eqnarray}
\label{biccupoverplus}
(q+r)^{\smallcup}&=&q^{\smallcup}+r^{\smallcup} \> \mbox{\femsteg $q,r \in {\rm B}$} \\ [1mm]
\label{biccupovermult}
(q\odot r)^{\smallcup}&=&q^{\smallcup}\odot r^{\smallcup} \\ [1mm]
\label{biccupidempot}
(q^{\smallcup})^{\smallcup}&=&q
\end{eqnarray}
Next we define the {\em complex square norm} $C\!N(q)$ of
the bicomplex number $q$:
\begin{equation}
\label{csquarenormb}
C\!N(q)=q\odot q^{\smallcup} \mbox{\femsteg $\lbrace$${\rm B}$-representation$\rbrace$}
\end{equation}
a formula of the same type as (\ref{squarenormh}).
With $q=(a,b)$ we get the ${\rm C^{2}}$-representation:
\begin{equation}
\label{csquarenormc2}
C\!N((a,b))= (a^{2}+b^{2},0) 
\end{equation}
It shows that $C\!N(q)$ is not a real number as $N(q)$ but
a complex number (which we nonetheless prefer to write in bicomplex notation).
$C\!N((a,b))$ is zero if and only if $a^{2}+b^{2}= 0$, which
is the case for $b=ia$ or $b=-ia$. 
Viewed as vectors $a$ and $b$ are then of equal length
and perpendicular to each other. We thus have the relation:
\begin{equation}
\label{cnzero}
C\!N((a,b))=0 \:\equiv\: (b=ia) \,\vee\, (b=-ia)
\end{equation}
The condition $C\!N(q)=0$ is crucially important in
bicomplex algebra and has no essential analogue in
complex algebra. We shall use the following terminology~\cite{price}:

\newtheorem{definition}{Definition}[section]
\begin{definition}
{\rm
A bicomplex number $q$ is called {\em singular} if
$C\!N(q)=0$, otherwise it is called {\em nonsingular}. $\Box$
}
\end{definition}
From the point of view of developing a bicomplex function
theory it is comforting to observe that the bicomplex
space contains simply connected domains that are not
singleton sets and contain merely nonsingular points.
Such domains are consequently termed {\em nonsingular domains}.

\medskip

\noindent {\bf Example.} The domain that consists of all bicomplex
numbers $(x+iy, z+iu)$, such that $x>0 \wedge y>0 \wedge z>0 \wedge u>0$ holds,
is nonsingular. $\Box$

\bigskip

\noindent For two bicomplex numbers $q$ and $r$ the complex square norm satisfies
\begin{equation}
\label{CNfact}
C\!N(q\odot r) = C\!N(q)\odot C\!N(r) 
\end{equation}
from which it follows that the bicomplex product preserves
nonsingularity if and only if both factors are nonsingular.

The singular bicomplex numbers $(a,ia)$ and $(a,-ia)$
are so-called {\em zero divisors}, because $(a,ia)\odot (a,-ia)=0$~\cite{riley}.
Hence, although $q=0 \,\vee\, r=0 \;\Rightarrow\; q\odot r=0$ holds for
all $q$ and $r$, the implication in the other direction is not generally valid.

As the inverse of $C\!N((a,b))$ we can take
\begin{equation}
\label{cninv}
\frac{1}{C\!N((a,b))} = \left(\frac{1}{a^{2}+b^{2}},0\right) \mbox{\femsteg , $a^{2}+b^{2}\neq 0$}
\end{equation}
The {\em inverse}
$q^{-1}$ of the bicomplex number $q$ is now defined by
\begin{equation}
\label{bicinvh}
q^{-1} = \frac{1}{C\!N(q)}\odot q^{\smallcup} \mbox{\femsteg , $C\!N(q)\neq 0$}
\end{equation}
The definition is of the same form as (\ref{qinverse}).
Note that only a nonsingular bicomplex number can
possess an inverse.
Inserting $q=(a,b)$ in (\ref{bicinvh}) yields upon application
of (\ref{bicmultdef}), (\ref{biccupdef}) and (\ref{cninv})
the equivalent ${\rm C^{2}}$-representation
\begin{equation}
\label{bicinvc2}
(a,b)^{-1} = \left(\frac{a}{a^{2}+b^{2}}\,,\,\frac{-b}{a^{2}+b^{2}}\right) \mbox{\femsteg , $a^{2}+b^{2}\neq 0$}
\end{equation}
At this point the definition of the {\em division}
of two bicomplex numbers $r$ and $q$ is straightforward:
\begin{equation}
\label{bicdivh}
\frac{r}{q}=r\odot q^{-1} \mbox{\femsteg , $C\!N(q)\neq 0$}
\end{equation}
The notation $\frac{r}{q}$ is meaningful
due to the commutativity of $\odot$. In ${\rm C^{2}}$ the same
formula looks like
\begin{equation}
\label{bicdivc2}
\frac{(c,d)}{(a,b)} = \left(\frac{c\cdot a+d\cdot b}{a^{2}+b^{2}}\:,\:\frac{d\cdot a-c\cdot b}{a^{2}+b^{2}}\right) \mbox{\femsteg , $a^{2}+b^{2}\neq 0$}
\end{equation}
Using the inverse (\ref{bicinvh}) we are also able
to extend the definition of exponentiation (\ref{expdef})
to negative integer exponents by
\begin{equation}
\label{bicnegexp}
q^{-n}=(q^{-1})^{n} \mbox{\femsteg , $n>0$ and $C\!N(q)\neq 0$}
\end{equation}
The definitions of the inverse and of the division operation
are thus based on the complex square norm $C\!N$.
The definition of the absolute value (modulus) of a bicomplex
number, however, cannot be based on $C\!N$; for this purpose we have
to borrow formula (\ref{normh}) from quaternionic algebra:
\begin{equation}
\label{normb}
\| q \| = \sqrt{N(q)} \mbox{\femsteg , $q\in {\rm B}$}
\end{equation}
Due to (\ref{squarenormr4}) the $\| \:\: \|$-function
may be identified with the Euclidean norm, which is known
to have the following properties:
\begin{eqnarray*}
&&\|q\| \geq 0 \mbox{\femsteg \femsteg \tresteg \tvasteg for all $q \in {\rm B}$} \\ [2mm]
&&\|q\| = 0 \;\equiv \; q=0 \\ [2mm]
&&\|\lambda q\| \;=\; \mid \!\lambda \!\mid \!\cdot \: \|q\| \mbox{\femsteg \tresteg for all $\lambda \in {\rm R}$ and $q\in {\rm B}$} \\ [2mm]
&&\|q+r\| \; \leq \; \|q\| + \|r\| \mbox{$\;\,$ \tresteg for all $q,r \in {\rm B}$}
\end{eqnarray*}
Moreover, formula (\ref{Ndecomp}) gives us for all bicomplex $(a,b)$
\[
\|(a,b)\| \; \leq \; \mid \! a \!\mid + \mid \! b \!\mid
\]
The notation $\mid \! c \!\mid$ stands for the absolute value
$\sqrt{N(c)}$ of the complex number $c$.

If we visualize the bicomplex number $q$
as a vector in ${\rm R^{4}}$, the absolute value $\|q\|$
is the length of the vector.

\subsection{Octonions and tricomplex numbers}
It remains to say something about the generalization
of the quaternions and the bicomplex numbers in the
form of {\em octonions} and {\em tricomplex numbers},
respectively.

An octonion can be expressed as a pair of
quaternions $(q,r)$, thus making the octonionic
space ${\rm O}$ 8-dimensional.
Octonionic addition is performed in the normal
way by component-wise addition, i.e.
$(q,r)+(s,t) = (q+s,r+t)$.
Octonionic multiplication $\otimes$ is defined
in terms of quaternionic multiplication by
\addtocounter{equation}{1}
\begin{tabbing}
XXXXXXXXX\=XXXXXXXXXXXXXXXXXXXXXXXXXXXXXXX\=XXXXXXXXXXXXXXXXXXXXXXXXXXXXXXXXXXXXXXXXXXXXXXXXXXXXXX\kill
(\theequation a)\> $(q,r)\otimes (s,t)=(q\times s - t^{*}\times r\,,\, r\times s^{*}+t\times q)$ \> $q,r,s,t \in {\rm H}$
\end{tabbing}

\noindent This formula should be compared with the
formulas we have employed for complex and quaternionic
multiplication:

\addtolength{\baselineskip}{+0.2cm}
\begin{tabbing}
XXXXXXXXX\=XXXXXXXXXXXXXXXXXXXXXXXXXXXXXXX\=XXXXXXXXXXXXXXXXXXXXXXXXXXXXXXXXXXXXXXXXXXXXXXXXXXXXXX\kill
(\theequation b)\> $(x,y)\cdot (z,u)=(xz - uy, yz+ux)$ \> $x,y,z,u \in {\rm R}$\\
(\theequation c)\> $(a,b)\times (c,d)=(a\cdot c -d^{*}\cdot b\,,\, b\cdot c^{*}+d\cdot a)$ \> $a,b,c,d \in {\rm C}$
\end{tabbing}
\addtolength{\baselineskip}{-0.2cm}

\noindent If the conjugate of a real number is taken to be the
number itself, (\theequation b) could be written
$(x,y)\cdot (z,u)=(xz - u^{*}y, yz^{*}+ux)$.
The formulas (\theequation a), (\theequation b) and (\theequation c) clearly
follow the same scheme and therefore complex numbers, quaternions
and octonions have many properties in common.
However, generalizing from quaternions
to octonions causes the loss of the associativity of
the $\otimes$-operation. Instead
$\otimes$ satisfies the weaker property of {\em alternation}:
$(m\otimes m)\otimes n = m\otimes (m\otimes n)$ and
$m\otimes (n\otimes n) = (m\otimes n)\otimes n$ for
$m,n \in {\rm O}$.

A tricomplex number, in turn, is expressed as a pair of
bicomplex numbers. Addition of two tricomplex numbers $(q,r)$ and $(s,t)$
is performed by component-wise addition and 
multiplication by
\[
(q,r)\odot (s,t)=(q\odot s - t\odot r\,,\, r\odot s+t\odot q) \mbox{\femsteg $q,r,s,t \in {\rm B}$}
\]
Here we have taken the liberty of
overloading the symbol $\odot$ by applying it to the multiplication
of both bicomplex and tricomplex numbers. This operation is 
isomorphic to complex multiplication (\theequation b),
as is bicomplex multiplication (\ref{bicmultdef}).

\newpage
\setcounter{equation}{0}

\section{Elementary bicomplex functions}
In this chapter we study functions of the type 
${\rm B} \rightarrow {\rm B}$, or bicomplex functions
of a bicomplex variable.
The benefits of using
complex pairs will now become so obvious that the alternative to base
one's reasoning entirely on one of the ${\rm R^{4}}$-representations
of the bicomplex numbers seems exceedingly unattractive.

We define the elementary bicomplex
functions, or the exponential function,
the hyperbolic and trigonometric functions, 
polynomials, the quotient function and the logarithm function. These 
functions contain their complex counterparts 
as special cases, as must be required of them if the
bicomplex function theory is to be a proper
generalization of the corresponding theory
in the complex domain~\cite{brown,cerny,goodstein,rade}.  

\subsection{The form of bicomplex functions}
Let the bicomplex variable $p$ be given by:
\begin{eqnarray}
\label{pdef}
p &=& (a,b) \\
(a,b) &=& (x+iy,z+iu) \nonumber
\end{eqnarray}
\noindent A bicomplex function $\psi$ of $p$ then has the form:
\begin{eqnarray}
\label{psidef}
\psi(p) &=& (\phi_{1}(a,b),\phi_{2}(a,b)) \\ [2mm]
(\phi_{1}(a,b),\phi_{2}(a,b)) &=& (\psi_{1}(x,y,z,u)+i\psi_{2}(x,y,z,u), \nonumber \\
& & \;\,\psi_{3}(x,y,z,u)+i\psi_{4}(x,y,z,u)) \nonumber
\end{eqnarray}
The functions $\phi_{k}, k=1,2$, are of type ${\rm C}^{2} \rightarrow {\rm C}$,
the functions $\psi_{k}, k=1,2,3,4$, of type ${\rm R}^{4} \rightarrow {\rm R}$.

Formulas (\ref{pdef}) and (\ref{psidef}) provide us with ${\rm B}$-, ${\rm C}^{2}$- and
${\rm R}^{4}$-representations of $p$ and $\psi$.
Because of our frequent reliance on the ${\rm C}^{2}$-representation
it would not be inappropriate to look at $\psi$ as
a function of two complex variables $a$ and $b$ 
rather than as a function of one bicomplex variable $p$.
In any case, what is essential in the sequel is the requirement
that the two-argument, complex functions $\phi_{k}$ should be holomorphic
in both $a$ and $b$, i.e. the partial derivatives
$\frac{\partial \phi_{k}}{\partial a}$ and $\frac{\partial \phi_{k}}{\partial b}$
should exist. Roughly speaking, the $\phi_{k}$'s meet this condition if their values are
obtained by algebraic operations and applications of
holomorphic functions so that
$a$ and $b$ enter in them only as {\em complex aggregates}
without separation of their real and imaginary parts~\cite{myskis}.

\bigskip

\noindent {\bf Example.} For $a=x+iy$ and $b=z+iu$ the functions $a^{2}+b^{2}$ and
$a\cdot \sin b$ are holomorphic in $a$ and $b$, not so $x^{2}+y^{2}+b^{2}$ and $a\cdot \sin(z+u)$. $\Box$

\bigskip

At this point we recall some basic results from complex analysis.
In general, if a complex function
\[
\omega(a,b) = \xi_{1}(x,y,z,u)+i\,\xi_{2}(x,y,z,u)
\]
is holomorphic in both $a=x+iy$ and $b=z+iu$
the partial derivative 
$\frac{\partial \omega}{\partial a}$ is given by either
of the formulas 
\begin{eqnarray}
\label{dthetadainx}
\frac{\partial \omega}{\partial a} &=& \partial_{x}\xi_{1}+i\partial_{x}\xi_{2} \\ [3mm]    
\label{dthetadainy}
\frac{\partial \omega}{\partial a} &=& \partial_{y}\xi_{2}-i\partial_{y}\xi_{1}
\end{eqnarray}
Their equivalence amounts to the truth of 
the Cauchy-Riemann equations
\begin{equation}
\label{CRequations}
\partial_{x}\xi_{1}=\partial_{y}\xi_{2} \tresteg  , \tresteg \partial_{x}\xi_{2}=-\partial_{y}\xi_{1}
\end{equation}
By making the replacements 
$a \rightarrow b$ , $x \rightarrow z$ , $y \rightarrow u$
in (\ref{dthetadainx}), (\ref{dthetadainy}) and 
(\ref{CRequations}) we obtain the corresponding formulas for
$\frac{\partial \omega}{\partial b}$.

\medskip

\noindent {\bf Remark.} Above we abbreviated the operators 
$\frac{\partial}{\partial x}$ and $\frac{\partial}{\partial y}$ by
$\partial_{x}$ and $\partial_{y}$. We employ this
notation when taking the derivative with respect to a
real variable, but not when taking it with respect to a
complex  variable.

We further left out the arguments of the functions.
We allow ourselves such an omission whenever it
is clear from the context which the arguments are. $\Box$

\bigskip

\noindent Returning to the function $\psi$ in (\ref{psidef}),
it should be emphasized that the holomorphism of its
component functions $\phi_{k}$ in $a$ and $b$
does {\em not} ensure that $\psi$ itself is holomorphic in $p$.
For $\psi$ to be holomorphic the $\phi_{k}$-functions must
fulfill a stronger condition, namely the {\em bicomplex
version} of the Cauchy-Riemann equations. The issue
will be discussed in Chapter 4. But the holomorphism of
the component functions does suffice to establish the continuity
of $\psi$. To do so we need the limit concept
$\lim_{p \rightarrow p_{0}} \psi(p) = w_{0}$,
which in the ordinary fashion is defined to be equivalent to
\begin{equation}
\label{fctlimit}
(\forall p,\epsilon: \epsilon >0: \: (\exists \delta: \delta > 0:\: \|p-p_{0} \|<\delta \:\Rightarrow \: \|\psi(p)-w_{0} \|<\epsilon ))
\end{equation}
where $p\! \in \!{\rm B}$ and $\epsilon, \delta \!\in \!{\rm R}$.
The absolute values should be computed using (2.\ref{normb}).
Geometrically interpreted, the limit is taken for all $p$
in a (deleted) neighbourhood of the bicomplex space or,
more precisely, a four-dimensional ball of radius $\delta$
centered at, but not including, the point $p_{0}$.

Falling back on (\ref{fctlimit}), one can define the
continuity of $\psi$ at a point $p_{0}$ by
\begin{equation}
\label{fctcont}
\lim_{\triangle p \rightarrow 0} \psi(p_{0}+\triangle p) = \psi(p_{0})
\end{equation}
The holomorphism of the functions
$\phi_{k}$, $k=1,2$, at the point $(a_{0},b_{0})$
means that they possess
partial derivatives
$\frac{\partial \phi_{k}}{\partial a}$ and $\frac{\partial \phi_{k}}{\partial b}$
at $(a_{0},b_{0})$.
Consequently, the functions are
continuous there, or
\begin{equation}
\label{phicont}
\lim_{(\triangle a, \triangle b) \rightarrow 0} \phi_{k}(a_{0}+\triangle a, b_{0}+\triangle b) = \phi_{k}(a_{0},b_{0}) \mbox{\femsteg , $k=1,2$}
\end{equation}
The limit of a complex function is assumed to be a
well-understood notion. If we rewrite (\ref{fctcont})
in the form of complex pairs using $p_{0}=(a_{0},b_{0})$ and 
$\triangle p = (\triangle a, \triangle b)$ we get an
equality, whose validity
follows from (\ref{phicont}). In other words, $\psi$ is continuous
at $p_{0}$.

In Chapter 4 we shall also need the following formula.
Let $\phi\!: (a,b) \rightarrow \phi(a,b)$ be a complex function
that is holomorphic in a domain $G\subseteq {\rm B}$.
The continuity of 
$\frac{\partial \phi}{\partial a}$ and $\frac{\partial \phi}{\partial b}$ 
means that at any point $(a,b)\in G$ we can write
\begin{equation}
\label{differential}
\phi(a+\triangle a,b+\triangle b)-\phi(a,b) \;=\; \frac{\partial \phi}{\partial a} \cdot \triangle a + \frac{\partial \phi}{\partial b} \cdot \triangle b + \epsilon(\triangle a, \triangle b)
\end{equation}
where $\epsilon(\triangle a, \triangle b)$ is a complex-valued function 
that tends to 0 more rapidly than $\triangle a$ and $\triangle b$ tend to 0.

\subsection{The exponential function}
With the help of the complex exponential, sine and cosine
functions we define the bicomplex exponential function by
\begin{equation}
\label{bicexpdef}
e^{p}=(e^{a}\!\cdot \cos b, e^{a}\!\cdot \sin b) \mbox{\femsteg , $p=(a,b)$}
\end{equation}
The component functions $\phi_{1}(a,b)=e^{a}\!\cdot \cos b$ and 
$\phi_{2}(a,b)=e^{a}\!\cdot \sin b$ are obviously holomorphic in
both $a$ and $b$. Definition (\ref{bicexpdef}) is of the same form
as the corresponding definition of the complex exponential function
rendered as a pair of reals:
\begin{equation}
\label{cexpdef}
e^{a}=(e^{x}\!\cdot \cos y, e^{x}\!\cdot \sin y) \mbox{\femsteg , $a=(x,y)$}
\end{equation}
With $b=(z,u)$ the complex functions $\cos b$ and $\sin b$
have the real pair representations:
\begin{eqnarray*}
\cos b &=& (\cos z \cosh u, -\sin z \sinh u) \\ [1mm]
\sin b &=& (\sin z \cosh u, \cos z \sinh u)
\end{eqnarray*}
Together with (\ref{cexpdef}) they enable us to
rewrite (\ref{bicexpdef}) as
a quadruple in ${\rm R^{4}}$. For $p=(a,b)=(x,y,z,u)$
we get
\begin{eqnarray}
\label{bicexpr4}
e^{p} &=& (\psi_{1}, \psi_{2}, \psi_{3}, \psi_{4}) \\
      &=& (\;e^{x}(\cos y \cos z \cosh u + \sin y \sin z \sinh u), \nonumber \\
      & & \;\;e^{x}(-\cos y \sin z \sinh u + \sin y \cos z \cosh u), \nonumber \\
      & & \;\;e^{x}(\cos y \sin z \cosh u - \sin y \cos z \sinh u), \nonumber \\
      & & \;\;e^{x}(\cos y \cos z \sinh u + \sin y \sin z \cosh u)\;) \nonumber
\end{eqnarray}
Turning our attention to the properties of $e^{p}$ we first observe that
for $p=(a,0)$ we get $e^{(a,0)}=(e^{a},0)$.
The bicomplex function thus generalizes the right complex function.

Definition (\ref{bicexpdef}) also gives us the periodicity of $e^{p}$.
In the complex plane $e^{a}$ has the period $im2\pi$, $m=\pm 1, \pm 2, \ldots ,$
while $\cos b$ and $\sin b$ have the period
$n2\pi$, $n=\pm 1, \pm 2, \ldots$. The period
$w$ of $e^{p}$ is therefore $(im,n)2\pi$, formally:
\begin{eqnarray}
\label{bicexpperiod}
e^{p+w}=e^{p} \mbox{\femsteg for } w&=&(im,n)2\pi \\
                                   m&=&0, \pm 1, \pm 2, \ldots \nonumber \\
                                   n&=&0, \pm 1, \pm 2, \ldots \nonumber
\end{eqnarray}
Setting $p=0$ yields as special case
\begin{equation}
\label{bicexpone}
e^{(im,n)2\pi}=1 \mbox{\femsteg , $m=0, \pm 1, \pm 2, \ldots$ and $n=0, \pm 1, \pm 2, \ldots$}
\end{equation}
The familiar {\em addition formula} $e^{a+b}=e^{a}\cdot e^{b}$
satisfied by the complex exponential function is valid in bicomplex
space, too, provided complex multiplication is replaced by bicomplex
multiplication:
\begin{equation}
\label{bicexpaddform}
e^{p+r}=e^{p}\odot e^{r} \mbox{\femsteg , $p=(a,b)$ and $r=(c,d)$}
\end{equation}
For $r=-p$ this formula becomes $e^{p}\odot e^{-p}=1$, which
shows that $e^{p}$ never gets the value zero. In Section
3.6 we shall obtain the stronger result that 
$e^{p}$ never gets a singular value.

\subsection{The hyperbolic and trigonometric functions}
With the exponential function at our disposal
we are now in a position to define the hyperbolic
functions $\cosh$ and $\sinh$ in ${\rm B}$ in the same way as it is
done in ${\rm C}$:
\begin{eqnarray}
\label{biccoshdef}
\cosh p&=&\frac{1}{2}(e^{p}+e^{-p}) \mbox{\femsteg , $p=(a,b)$} \\ [3mm]
\label{bicsinhdef}
\sinh p&=&\frac{1}{2}(e^{p}-e^{-p})
\end{eqnarray}
The use of (\ref{bicexpdef}) in these formulas
yields after straightforward calculation the complex pair representations
\begin{eqnarray}
\label{biccoshc2}
\cosh p&=&(\cosh a \cdot \cos b, \sinh a \cdot \sin b) \mbox{\femsteg , $p=(a,b)$} \\ [2mm]
\label{bicsinhc2}
\sinh p&=&(\sinh a \cdot \cos b, \cosh a \cdot \sin b)
\end{eqnarray}
Note that for $p=(a,0)$ the functions reduce to their
complex analogues.
Moreover, due to their definitions they have the same
period as $e^{p}$.

Bicomplex division (2.\ref{bicdivh}) gives us the $\tanh$-function
\begin{equation}
\tanh p = \frac{\sinh p}{\cosh p} \mbox{\femsteg , $p=(a,b)$}
\end{equation}
By applying (\ref{biccoshc2}), (\ref{bicsinhc2}) and (2.\ref{bicdivc2})
in this formula we get after some manipulation
\begin{equation}
\tanh p = \left(\frac{\tanh a}{(\cos b)^{2}+(\tanh a \cdot \sin b)^{2}}\:,\: \frac{\tan b}{(\cosh a)^{2}+(\sinh a \cdot \tan b)^{2}} \right)
\end{equation}
To define the bicomplex cosine and sine functions we
take the corresponding complex functions as starting point:
\begin{eqnarray}
\label{ccosdef}
\cos a&=&\frac{1}{2}(e^{ia}+e^{-ia}) \mbox{\femsteg , $a=x+iy$} \\ [3mm]
\label{csindef}
\sin a&=&\frac{1}{2i}(e^{ia}-e^{-ia})
\end{eqnarray}
Employing real pairs to express complex numbers we have
$a=(x,y)$, $i=(0,1)$ and
\[
ia\:=\:(0,1)\cdot(x,y)\:=\:(-y,x)
\]
Transferring this to ${\rm B}$ we replace $a=(x,y)$ by
$p=(a,b)$ and $i=(0,1)$ by $j=(0,1)$ to obtain
\begin{equation}
\label{jodotp}
j\odot p\:=\:(0,1)\odot (a,b)\:=\:(-b,a)
\end{equation}
We choose the expression $j\odot p$ as the bicomplex
correspondence of the complex $ia$ and generalize
(\ref{ccosdef}) and (\ref{csindef}) by
\begin{eqnarray}
\label{biccosdef}
\cos p&=&\frac{1}{2}(e^{j\odot p}+e^{-j\odot p}) \mbox{\femsteg , $p=(a,b)$} \\ [3mm]
\label{bicsindef}
\sin p&=&\frac{1}{2j}\odot(e^{j\odot p}-e^{-j\odot p})
\end{eqnarray}
Their ${\rm C^{2}}$-representations are
\begin{eqnarray}
\label{biccosc2}
\cos p &=& (\cos a \cdot \cosh b, -\sin a \cdot \sinh b) \mbox{\femsteg , $p=(a,b)$} \\ [2mm]
\label{bicsinc2}
\sin p &=& (\sin a \cdot \cosh b, \cos a \cdot \sinh b) 
\end{eqnarray}
Because $\cos a$ and $\sin a$ both have the period
$m2\pi$, $m=\pm 1, \pm 2, \ldots ,$
while $\cosh b$ and $\sinh b$ have the period
$in2\pi$, $n=\pm 1, \pm 2, \ldots$, the period of
$\cos p$ and $\sin p$ is $(m,in)2\pi$.

The ${\rm B}$- and ${\rm C^{2}}$-representations of the 
bicomplex $\tan$-function are:
\begin{equation}
\label{bictandef}
\tan p = \frac{\sin p}{\cos p} \mbox{\femsteg , $p=(a,b)$}
\end{equation}
\begin{equation}
\label{bictanc2}
\tan p = \left(\frac{\tan a}{(\cosh b)^{2}+(\tan a \cdot \sinh b)^{2}}\: , \: \frac{\tanh b}{(\cos a)^{2}+(\sin a \cdot \tanh b)^{2}} \right)
\end{equation}
In the definitions of the preceding functions we observe in
every case that they are the same as the definitions of
their complex counterparts but with complex
notions replaced by bicomplex ones. This isomorphism
manifests itself in the properties of the functions, too. For example,
the familiar addition formulas of $\cosh$, $\sinh$, $\cos$ and $\sin$
can directly be taken over to bicomplex space:
\begin{eqnarray*}
\cosh(p+r) &=& \cosh p \odot \cosh r + \sinh p \odot \sinh r \\ [1mm]
\sinh(p+r) &=& \sinh p \odot \cosh r + \cosh p \odot \sinh r \\ [1mm]
\cos(p+r) &=& \cos p \odot \cos r - \sin p \odot \sin r \\ [1mm]
\sin(p+r) &=& \sin p \odot \cos r + \cos p \odot \sin r
\end{eqnarray*}
We also note that the functions are interconnected through
the formulas:
\begin{eqnarray*}
\cos p&=&\cosh(j\odot p) \\ [3mm]
\sin p&=&\frac{1}{j}\odot \sinh(j\odot p)
\end{eqnarray*}

\subsection{Polynomials}
Bicomplex polynomials are formed as sums of
bicomplex integer powers $p^{n}$ multiplied by
bicomplex constants.

As the simplest polynomial we regard the constant function
\begin{equation}
\label{bicconstfunc}
C(p)=p_{0} \mbox{\femsteg , $p_{0}$ a bicomplex constant}
\end{equation}
The bicomplex identity function is
\begin{equation}
\label{bicidfunc}
I(p)=(a,b) \mbox{\femsteg , $p=(a,b)$}
\end{equation}
and the bicomplex power function
\begin{equation}
\label{bicpowfunc}
F(p)=(a,b)^{n} \mbox{\femsteg , $p=(a,b)$ and $n$ an integer $\geq 0$}
\end{equation}
Here the exponentiation should be performed according to (2.\ref{expdef}).
We have in particular
\begin{eqnarray*}
& &p^{0}=(1,0) \\
& &p^{1}=(a,b) \\
& &p^{2}=(a^{2}-b^{2},2ab)
\end{eqnarray*}
The components of the complex pairs of these equations are examples of
so-called {\em harmonic polynomials}
of $a$ and $b$~\cite{abramovitz}. These polynomials are also
recursively definable: if we denote them 
by $G_{n}(a,b)$ and $H_{n}(a,b)$ we have
\begin{equation}
p^{n}=(G_{n}(a,b) , H_{n}(a,b)) \mbox{\femsteg , $p=(a,b)$ and $n\geq 0$}
\end{equation}
From the identities $p^{0}=(1,0)$ and $p^{n+1}=p\odot p^{n}$ we then
derive the recursive scheme:
\begin{eqnarray*}
&&G_{0}(a,b)=1 \; , \; H_{0}(a,b)=0  \\ [1mm] 
&&G_{n+1}(a,b)\;=\;a\cdot G_{n}(a,b)-b\cdot H_{n}(a,b) \\ [1mm]
&&H_{n+1}(a,b)\;=\;a\cdot H_{n}(a,b)+b\cdot G_{n}(a,b) 
\end{eqnarray*}
The general form of an ${\rm n^{th}}$ degree bicomplex polynomial
with argument $p=(a,b)$ is thus
\begin{equation}
\label{bicpolynomial}
Q(p) = d_{n}\odot p^{n}+d_{n-1}\odot p^{n-1}+\ldots +d_{1}\odot p+d_{0}
\end{equation}
where $d_{k}$ stands for a bicomplex constant.

Similarly, one can define {\em bicomplex power series} of the form
$(\Sigma k\!:k\geq 0\!:d_{k}\odot p^{k})$. Their meaningfulness
depends as usual on their convergence.

\subsection{The quotient function}
Let the bicomplex functions $\psi$ and $\theta$ be given by
\begin{eqnarray*}
\psi(p) &=& (\phi_{1}(a,b),\phi_{2}(a,b)) \mbox{\femsteg , $\; p=(a,b)$}\\ [1mm]
\theta(p) &=& (\omega_{1}(a,b),\omega_{2}(a,b))
\end{eqnarray*}
For the {\em inversion} $\theta^{-1}$ we first have
on application of (2.\ref{bicinvc2})
\begin{equation}
\label{bicinvfunc}
\theta^{-1} = \left(\frac{\omega_{1}}{\omega_{1}^{2}+\omega_{2}^{2}}\,,\,\frac{-\omega_{2}}{\omega_{1}^{2}+\omega_{2}^{2}}\right)
\end{equation}
The function is well-defined in the domain of ${\rm B}$ where $\omega_{1}^{2}+\omega_{2}^{2}\neq 0$, that is
where $C\!N((\omega_{1},\omega_{2}))\neq 0$. A bicomplex
inversion can have rather more zeros in the denominator than
the corresponding complex function, a point to remember when
one attempts to generalize e.g. Cauchy's integral formula.

The inversion of the identity function (\ref{bicidfunc}) is a special case of (\ref{bicinvfunc}):
\begin{equation}
\label{bicidinvfunc}
I^{-1}(p) = \left(\frac{a}{a^{2}+b^{2}}\,,\,\frac{-b}{a^{2}+b^{2}}\right) \mbox{\femsteg , $\; p=(a,b)$}
\end{equation}
The {\em quotient function} $\psi / \theta = \psi \odot \theta^{-1}$
has the form
\begin{equation}
\label{bicquotfunc}
\frac{\psi}{\theta} = \left(\frac{\phi_{1}\cdot \omega_{1}+\phi_{2}\cdot \omega_{2}}{\omega_{1}^{2}+\omega_{2}^{2}}\:,\:\frac{\phi_{2}\cdot \omega_{1}-\phi_{1}\cdot \omega_{2}}{\omega_{1}^{2}+\omega_{2}^{2}}\right)
\end{equation}
An important category of quotient functions consists of
{\em rational functions} of the type
\begin{equation}
R(p)=\frac{Q1(p)}{Q2(p)} \mbox{\femsteg , $\; p=(a,b)$}
\end{equation}
where $Q1(p)$ and $Q2(p)$ are bicomplex polynomials
without common factors. Borrowing terminology from complex function
theory we say that the zeros of $Q2(p)$ are the {\em poles} of $R(p)$.

\subsection{The logarithm function}
We shall define the bicomplex logarithm function 
as the inverse of the bicomplex exponential function.
The definition will be based on the
following result.
\newtheorem{theorem}{Theorem}[section]
\begin{theorem}
The bicomplex number $q=(\gamma + i\delta , \epsilon + i\eta)$ 
can equivalently be rendered in the form
\begin{equation}
\label{bicpolarform}
q=(v,0)\odot (\cos w,\sin w) 
\end{equation}
where $v$ and $w$ are complex numbers, if and only if
q is nonsingular or equal to zero.
\end{theorem}

\noindent {\em Proof}. It is without loss of generality possible to express $v$ and $w$ in terms of the real
numbers $r,y,z,u$ so that
\begin{equation}
\label{vandw}
v=re^{iy} \tresteg  , \tresteg w=z+iu
\end{equation}
Our task is therefore to investigate under
what conditions the equation
$(\gamma + i\delta , \epsilon + i\eta)$ = $(re^{iy},0)\odot (\cos(z+iu) ,\sin(z+iu))$ has a 
solution in the unknowns $r,y,z,u$. The right-hand side is equal to
\[
\left(\frac{1}{2}re^{iy}\cdot (e^{i(z+iu)}+e^{-i(z+iu)})\:,\,\frac{1}{2i}re^{iy}\cdot (e^{i(z+iu)}-e^{-i(z+iu)})\right)
\]
and by rewriting this expression using the addition formula as well as other
basic properties of the complex exponential function,
one can show that the equation is equivalent to the 
following system of equations:
\begin{eqnarray*}
\frac{1}{2}r(e^{-u}\cos(y+z) + e^{u}\cos(y-z)) &=& \gamma \\ [3mm]
\frac{1}{2}r(e^{-u}\sin(y+z) + e^{u}\sin(y-z)) &=& \delta \\ [3mm]
\frac{1}{2}r(e^{-u}\sin(y+z) - e^{u}\sin(y-z)) &=& \epsilon \\ [3mm]
\frac{1}{2}r(-e^{-u}\cos(y+z) + e^{u}\cos(y-z)) &=& \eta
\end{eqnarray*}
If $q=0$, i.e. $\gamma=\delta=\epsilon=\eta=0$, these equations are
solved by $r=0$, in which case the values of $y$, $z$ and $u$ are indeterminate.
If $q\neq 0$ and singular, i.e. of the form $(\gamma + i\delta, -\delta+i\gamma)$
or $(\gamma + i\delta, \delta-i\gamma)$, the equations admit no solution
(since they imply either $\sin(y+z)=0 \, \wedge \, \cos(y+z)=0$ or $\sin(y-z)=0 \, \wedge \, \cos(y-z)=0$).
However, if $q$ is nonsingular they are solved by:
\begin{equation}
\label{vandwsolution}
 \begin{array}{rclrcll}
 r&=&(K\cdot L)^{1/4} & \mbox{\tresteg} K&=&(\gamma + \eta)^{2}+(\delta - \epsilon)^{2}& \\ [4mm]
 y&=&\frac{1}{2}(M+N) & \mbox{\tresteg} L&=&(\gamma - \eta)^{2}+(\delta + \epsilon)^{2}&\\ [4mm]
 z&=&\frac{1}{2}(N-M) & \mbox{\tresteg} M&=&\arctan\left(\frac{\delta-\epsilon}{\gamma+\eta}\right) + n\pi & \mbox{\tvasteg , $n=0 \,\vee\, n=1$} \\ [4mm]
 u&=&\frac{1}{4}\ln\left|\frac{K}{L}\right| & \mbox{\tresteg} N&=&\arctan\left(\frac{\delta+\epsilon}{\gamma-\eta}\right) + n\pi & \mbox{\tvasteg , $n=0 \,\vee\, n=1$}
 \end{array}
\end{equation}
The nonsingularity of $q$ ensures that $K$ and $L$ are both nonzero, thus
making $u$ well-defined. When $M$ is evaluated the
principal value in the interval $(-\frac{\pi}{2},\frac{\pi}{2})$
should be chosen for the $\arctan$-function. The integer
$n$ should be set to 0 or 1 depending on whether the denominator
$\gamma + \eta$ is positive or negative. If the denominator is zero
$M=\frac{\pi}{2}$ or $M=\frac{3\pi}{2}$ depending on whether
the numerator $\delta-\epsilon$ is positive or negative. The same
applies, mutatis mutandis, to the evaluation of $N$.
 
The theorem follows from the existence of the solution (\ref{vandwsolution}).

\noindent $\Box$

\medskip

\noindent In the bicomplex number representation $q=(v,0)\odot (\cos w,\sin w)$ we call
$v$ the {\em scale factor}
and $w$ the {\em complex argument} of $q$. We shall use the notations:
\begin{equation}
\label{sfca}
 \begin{array}{rcl}
 v &\!\! =\!\! & s\!f(q) \mbox{\femsteg \femsteg$\lbrace$scale factor of $q$$\rbrace$} \\[2mm]
 w &\!\! =\!\! & ca(q) \mbox{\femsteg \femsteg$\lbrace$complex argument of $q$$\rbrace$}
 \end{array}
\end{equation}
The function $ca(q)$ is multi-valued, because its real part
is not uniquely determined; for a given $q$ its infinitely many
values differ from each other by multiples of $2\pi$. It is
possible to identify a principal branch of $ca(q)$, denoted by
$Ca(q)$, that is characterized by the fact that its real part is
in a certain interval, $[0,2\pi)$ say. Formally:
\[
w=Ca(q) \;\,\equiv\;\, w=ca(q) \,\wedge\, 0\leq {\rm Re}(w) <2\pi
\]
The complex exponential function gets every nonzero value in
${\rm C}$, hence if $q\neq 0$ in formula (\ref{bicpolarform})
$v$ can be replaced by $e^{a}$ for some $a\in {\rm C}$.
This gives $q=(e^{a}\!\cdot \cos w, e^{a}\!\cdot \sin w)$
or the bicomplex exponential function. With
Theorem 3.2 we then deduce that this function {\em gets every nonsingular
value and no singular value}.

We denote the bicomplex logarithm function by ${\rm blog}$ and
define it as the solution of the equation $e^{p}=q$:
\[
{\rm blog}(q)=p \;\; \equiv \;\; e^{p}=q \mbox{\femsteg for $p,q \in {\rm B}$}
\]
Clearly, ${\rm blog}(q)$ is only defined for nonsingular $q$.
If we set $p=(a,b)$ and lean on (\ref{bicpolarform}) and the definition
of the exponential function, $e^{p}=q$ becomes
\[
(e^{a},0)\odot (\cos b,\sin b)\;=\;(v,0)\odot (\cos w,\sin w)
\]
This equation is satisfied if and only if
\begin{eqnarray*}
a&=&\log(v) \\ [1mm]
b&=&w+n2\pi \mbox{\femsteg , $n$ an integer}
\end{eqnarray*}
where $\log$ stands for the complex logarithm function.
By making use of the scale factor and complex argument
of $q$ --- see (\ref{sfca}) --- we further obtain:
\begin{eqnarray*}
a&=&\log(s\!f(q)) \mbox{\femsteg\tvasteg , $q$ nonsingular} \\ [1mm]
b&=&ca(q)+n2\pi \mbox{\femsteg , $n$ an integer} 
\end{eqnarray*}
Thus, we have derived the following expression
for ${\rm blog}$:
\begin{equation}
\label{blogdef}
{\rm blog}(q)=(\log(s\!f(q))\,,\, ca(q)+n2\pi) \mbox{\tresteg , $q$ nonsingular and $n$ an integer}
\end{equation}
Observe the isomorphism between this function and its
complex counterpart $\log(v)=\ln|v|+ i\arg(v)$.

The ${\rm blog}$-function is multi-valued, because
both the $\log$- and the $ca$-function are multi-valued. If
$p$ is a bicomplex logarithm of the nonsingular bicomplex
number $q$, then so are $p+(im,n)2\pi$ for all integers
$m$ and $n$. The principal branch of ${\rm blog}$, which
we denote by ${\rm Blog}$, is a function, however, and consists by definition of the
principal branches of the component functions:
\begin{equation}
\label{Blogdef}
{\rm Blog}(q)=({\rm Log}(s\!f(q))\,,\, Ca(q)) \mbox{\femsteg , $q$ nonsingular}
\end{equation}
The principal branch ${\rm Log}$ of the complex logarithm
function is defined in the normal manner.

The bicomplex logarithm function allows us to
extend the definition of the power function $q^{r}$
to bicomplex exponents:
\begin{equation}
\label{genpowerfct}
q^{r} = e^{r\odot {\rm Blog}(q)} \mbox{\femsteg for $q,r \in {\rm B}$ and $q$ nonsingular}
\end{equation}

\noindent {\bf Example.} Just for fun we shall evaluate $(0,i)^{(0,i)}$.
The bicomplex number $(0,i)$ is nonsingular, hence with the help of
(\ref{vandw}) and (\ref{vandwsolution}) we may render it in the form (\ref{bicpolarform}):
\[
(e^{i\frac{\pi}{2}},0)\odot (\cos \frac{\pi}{2}\,,\, \sin \frac{\pi}{2})
\]
This yields
\[
{\rm Blog}((0,i))\;=\;({\rm Log}(e^{i\frac{\pi}{2}})\,,\,\frac{\pi}{2})\;=\;(i\frac{\pi}{2}\,,\,\frac{\pi}{2})
\]
and further
\[
(0,i)^{(0,i)} \;=\; e^{(0,i)\odot {\rm Blog}((0,i))} \;=\; e^{(-i\frac{\pi}{2},-\frac{\pi}{2})} \;=\; (0,i)
\]
$\Box$

\bigskip

\noindent We are also able to define the inverse functions of
the trigonometric functions in the same way as it is done in ${\rm C}$.
For example, the inverse cosine is obtained as the solution of
\[
\cos p = \frac{1}{2}(e^{j\odot p}+e^{-j\odot p}) = q \mbox{\femsteg for $p,q \in {\rm B}$}
\]
Solving for $e^{j\odot p}$ and taking the bicomplex logarithm
yields the bicomplex version of a familiar formula:
\[
\arccos q = -j\odot {\rm blog}(q \pm \sqrt{q^{2}-1})
\]
It is of course required that ${\rm blog}$'s argument 
$q \pm \sqrt{q^{2}-1}$ is nonsingular.

\newpage
\setcounter{equation}{0}

\section{Differentiation of bicomplex functions}
We give the derivative of a bicomplex function a definition
that is isomorphic to the corresponding definition in
${\rm C}$ and because, in addition, the bicomplex operations
are isomorphic to the complex ones, we may expect that the
normal rules for computing the derivative of polynomials,
sums, products, quotients, etc. remain valid in ${\rm B}$.
However, in defining the bicomplex derivative we have to
reckon with the fact that the bicomplex numbers do not
form a field in the ordinary sense due to the lack of
inverses of the singular numbers. The possibility to
restrict the discussion to nonsingular numbers only
must be excluded, since they do not form a subfield in ${\rm B}$:
if $a$ and $b$ are two nonsingular bicomplex numbers, then
$a+b$ need not in general be nonsingular, although
$a\odot b$ is so (see formula~(\ref{CNfact})). In spite of this,
the problems related to singularities can be solved.

We shall obtain ${\rm B}$-, ${\rm C}^{2}$- and ${\rm R}^{4}$-representations
for the bicomplex derivative. The ${\rm C}^{2}$-representation is
particularly useful and will give us the bicomplex version of the Cauchy-Riemann equations,
which every differentiable, bicomplex function
has to satisfy. We shall show that such a function possesses derivatives
of all orders.

\subsection{Definition of the derivative}
Following the standard procedure one is tempted
to define the derivative of the bicomplex function $\psi$
at a point $p$ as the limit
\[
\psi\,'(p) = \lim_{\triangle p \rightarrow 0} \frac{\psi(p+\triangle p)-\psi(p)}{\triangle p}
\]
A problem then arises, because the fraction 
$(\psi(p+\triangle p)-\psi(p))/\triangle p$ 
is not defined for each singular $\triangle p$
in the neighbourhood  where the limit is taken.
An obvious remedy is to confine limit taking only to
nonsingular $\triangle p$, which demands that
we apply a weaker formula than (3.\ref{fctlimit})
for the limit concept, namely
\begin{eqnarray}
\label{fctlimitnonsing1}
(\forall \triangle p,\epsilon: \;C\!N(\triangle p)\neq 0 \,\wedge \, \epsilon >0: &&\!\!\!\!\!\!\!(\exists \delta: \delta > 0: \| \triangle p \|<\delta \; \:\Rightarrow \\ [2mm]&&\;\;\;\; \left\| \frac{\psi(p+\triangle p)-\psi(p)}{\triangle p} - \psi\,'(p) \right\|<\epsilon )) \nonumber
\end{eqnarray}
The dummies satisfy $\triangle p \in {\rm B}$ and $\epsilon, \delta \in {\rm R}$.
Recall from Section 2.5 that $\triangle p$ is nonsingular
if $C\!N(\triangle p)\neq 0$. The above limit condition is denoted by
\[
\psi\,'(p) = \lim_{\stackrel{{\scriptstyle \triangle p \rightarrow 0}}{C\!N(\triangle p)\neq 0}}\frac{\psi(p+\triangle p)-\psi(p)}{\triangle p}
\]
Formula (\ref{fctlimitnonsing1}) is an instance of the
following limit condition that pertains to bicomplex quotients in general:
\newpage
\begin{eqnarray}
\label{fctlimitnonsing2}
(\forall p,\epsilon: \;C\!N(\theta(p))\neq 0 \,\wedge \, \epsilon >0: &&\!\!\!\!\!\!\!\! (\exists \delta: \delta > 0: \\ 
&& \| p-p_{0} \|<\delta \;\Rightarrow\;  \left\| \frac{\psi(p)}{\theta(p)} - w_{0} \right\|<\epsilon )) \nonumber
\end{eqnarray}
The corresponding limes notation is
\[
w_{0} = \lim_{\stackrel{{\scriptstyle p \rightarrow p_{0}}}{C\!N(\theta(p))\neq 0}}\frac{\psi(p)}{\theta(p)}
\]
After these preliminaries we can state our first version of the
definition of the bicomplex derivative:
\begin{definition}
{\rm
Let $\psi$ be a bicomplex function whose domain of definition
contains a neighbourhood of the point $p$. The derivative
of $\psi$ at the point $p$ is defined by the equation
\begin{equation}
\label{bicderdef}
\psi\,'(p) = \lim_{\stackrel{{\scriptstyle \triangle p \rightarrow 0}}{C\!N(\triangle p)\neq 0}}\frac{\psi(p+\triangle p)-\psi(p)}{\triangle p}
\end{equation}
provided this limit exists.  $\Box$
}
\end{definition}
A bicomplex function that has a derivative at the point $p$
is said to be differentiable or {\em holomorphic} at $p$.
If the function is holomorphic at all points of a domain
$G\subseteq {\rm B}$ it is said to be holomorphic in $G$.
The derivative $\psi\,'(p)$ is as usual also written $\frac{d\psi}{dp}$.

\medskip

\noindent {\bf Remark.} By a {\em domain} 
we understand in the sequel a connected, open set of points in the bicomplex
space (occasionally in the complex plane). $\Box$

\bigskip

\noindent Definition 4.1  could at least formally be strengthened by the
requirement that for every singular $\triangle p$ such that
$p+\triangle p$ is in the neighbourhood of $p$ in question,
the fraction $\Omega(\triangle q) = (\psi(p+\triangle q)-\psi(p))/\triangle q$
should approach a finite limit $w$ as $\triangle q$ approaches $\triangle p$, formally:
\[
w = \lim_{\stackrel{{\scriptstyle \triangle q \rightarrow \triangle p}}{C\!N(\triangle q)\neq 0}}\frac{\psi(p+\triangle q)-\psi(p)}{\triangle q} \mbox{\femsteg , $C\!N(\triangle p)= 0$}
\]
Note that $\triangle q$ must be nonsingular. Note also
that different values of $w$ can be associated with different values of
$\triangle p$. The existence of
such a limit means that it is possible to remove the singularity
of $\Omega$ at each $\triangle p$ by defining $\Omega(\triangle p)=w$,
which makes the condition $C\!N(\triangle p)\neq 0$ superfluous
in (\ref{bicderdef}). Nevertheless, this strengthening just
postpones the moment when one has to compute a limit in a neighbourhood consisting of
nonsingular points only.

\bigskip

\noindent {\bf Example.} Let $\psi$  be the bicomplex square
function $\psi(p)=p^{2}$, which yields
\[
\frac{\psi(p+\triangle p)-\psi(p)}{\triangle p} = \frac{2p \odot \triangle p +(\triangle p)^{2}}{\triangle p}
\]
The limit value $2p +\triangle p$ of the right-hand side exists for
all singular $\triangle p$. Of these $\triangle p=0$ gives the derivative. $\Box$

\bigskip

\noindent Yet another alternative is to introduce the derivative
of $\psi$ by a definition of Fr\'{e}chet-type:
if there exists a bicomplex number $K$ and a bicomplex function
$\epsilon$ that satisfy the equations 
$\lim_{\triangle p \rightarrow 0} \epsilon(\triangle p) = 0$, 
$\epsilon(0)=0$ and
\begin{equation}
\label{stolzcond}
\psi(p+\triangle p)-\psi(p) = K\odot \triangle p + \epsilon(\triangle p)\odot \triangle p
\end{equation}
then the derivative $\psi\,'(p)$ is equal to $K$.
Thanks to the absence of denominators in (\ref{stolzcond}),
the problems related to singular bicomplex numbers are at first sight
circumvented. However, the identification of a suitable
function $\epsilon$ presupposes in general the existence of
removable singularities of the fraction 
$(\psi(p+\triangle p)-\psi(p))/\triangle p$. The equivalence
of this definition and Definition 4.1 has been demonstrated by Price~\cite{price},
who refers to (\ref{stolzcond}) as the {\em strong Stolz condition}.

Independently of which of the above definitions is applied,
one obtains by employing standard techniques the normal rules for computing the derivatives
of sums, products and quotients of functions:
\begin{eqnarray*}
(\psi + \theta)\,'(p)\!\!&=&\!\!\psi\,'(p) + \theta\,'(p) \\ [1mm]
(\psi \odot \theta)\,'(p)\!\!&=&\!\!\psi\,'(p)\odot \theta(p) + \psi(p)\odot \theta\,'(p) \\[2mm]
\left(\frac{\psi}{\theta}\right)'(p)\; \!\!&=&\!\! \; \frac{\psi\,'(p)\odot \theta(p) - \psi(p)\odot \theta\,'(p)}{(\theta(p))^{2}} \mbox{\femsteg , $C\!N(\theta(p))\neq 0$}
\end{eqnarray*}

\bigskip  

\noindent Moreover, for the bicomplex composite function
$(\psi \circ \theta)(p) = \psi(\theta(p))$ the {\em chain rule} holds:
\[
\label{bicderchain}
(\psi \circ \theta)\,'(p) = \psi\,'(\theta(p))\odot \theta\,'(p) 
\]
The derivative of the power function is given by the normal formula
\begin{equation}
\label{pownder}
\frac{dp^{n}}{dp} = np^{n-1}
\end{equation}
provided $n$ is a positive integer and $p$ {\em any} bicomplex number
or $n$ is a negative integer and $p$ a {\em nonsingular} bicomplex number.
In the next section this derivative will be generalized.

\subsection{Complex pair formulas for the derivative}
We shall seek formulas for the derivative
of a bicomplex function $\psi$ having the complex
pair representation
\begin{equation}
\label{psidef42}
\psi(p)=(\phi_{1}(a,b),\phi_{2}(a,b)) \mbox{\femsteg , $\; p=(a,b)$}
\end{equation}
The calculations are essentially the same as the ones that lead
to formulas (3.\ref{dthetadainx}) and (3.\ref{dthetadainy}) for
the complex derivative, but due to the possibility of singular
denominators it is necessary to go through them in detail.

We investigate what limits definition (\ref{bicderdef}) yields if we 
let $\triangle p$ approach 0 in two ways so that $\triangle p =(\triangle a,0)$
or $\triangle p =(0,\triangle b)$. Then, a key observation is
that, provided $\triangle a \neq 0$ and $\triangle b \neq 0$, $\triangle p$
{\em is never singular}. 
In other words, we do not need to worry about singularities
in the subsequent formulas.

On account of (2.\ref{bicinvc2}) the inverse of $\triangle p$ is
\begin{eqnarray}
\label{triangelpinv1}
\frac{1}{\triangle p}&=&\left(\frac{1}{\triangle a}\;,\;0\right)\; \mbox{\femsteg $\;\,$if $\triangle p=(\triangle a,0)$} \\[3mm]
\label{triangelpinv2}
\frac{1}{\triangle p}&=&\left(0\;,\;-\frac{1}{\triangle b}\right) \mbox{\femsteg if $\triangle p=(0,\triangle b)$} 
\end{eqnarray}
With respect to the case $\triangle p=(\triangle a,0)$ we get
using (\ref{psidef42}) and (\ref{triangelpinv1})
\[
\frac{\psi(p+\triangle p)-\psi(p)}{\triangle p} \;\;= \;\;\left(\frac{\phi_{1}(a+\triangle a,b)-\phi_{1}(a,b)}{\triangle a}\;,\;\frac{\phi_{2}(a+\triangle a,b)-\phi_{2}(a,b)}{\triangle a}\right)
\]
As $\triangle p$ and $\triangle a$ approach 0 simultaneously, the components
of the right-hand side approach the complex derivatives
$\frac{\partial \phi_{1}}{\partial a}$ and $\frac{\partial \phi_{2}}{\partial a}$.
The bicomplex derivative therefore takes the form
\begin{equation}
\label{bicderc21}
\frac{d\psi}{dp} = \left(\frac{\partial \phi_{1}}{\partial a}\;, \frac{\partial \phi_{2}}{\partial a} \right) \mbox{\femsteg , $\; p=(a,b)$}
\end{equation} 
If we deal with the second case $\triangle p=(0,\triangle b)$
in the same way we get with (\ref{psidef42}) and (\ref{triangelpinv2})
\[
\frac{\psi(p+\triangle p)-\psi(p)}{\triangle p} \;\;= \;\;\left(\frac{\phi_{2}(a,b+\triangle b)-\phi_{2}(a,b)}{\triangle b}\;,\;-\frac{\phi_{1}(a,b+\triangle b)-\phi_{1}(a,b)}{\triangle b}\right)
\]
Letting $\triangle p$ and $\triangle b$ approach 0 results in
\begin{equation}
\label{bicderc22}
\frac{d\psi}{dp} = \left(\frac{\partial \phi_{2}}{\partial b}\;, -\frac{\partial \phi_{1}}{\partial b} \right) \mbox{\femsteg , $\; p=(a,b)$}
\end{equation}
This formula must be equivalent to (\ref{bicderc21}), a requirement whose consequences
will be examined shortly.

We are now equipped to find more derivatives of
the elementary functions. Applied to the exponential
function $e^{p}=(e^{a}\!\cdot \cos b, e^{a}\!\cdot \sin b)$
both (\ref{bicderc21}) and (\ref{bicderc22}) give
\begin{equation}
\label{bicexpder}
\frac{de^{p}}{dp}=e^{p}
\end{equation}
In order to find the derivative of the principal branch
of the bicomplex logarithm function (3.\ref{Blogdef})
we assume that the four bicomplex numbers $p$, $\triangle p$, $q$, $\triangle q$
satisfy $q=e^{p} \; \wedge\; q+\triangle q = e^{p+\triangle p}$, thereby
implying that $q$ and $q+\triangle q$ are nonsingular. This yields
\[
\frac{{\rm Blog}(q+\triangle q)-{\rm Blog}(q)}{\triangle q} \; = \; \frac{\triangle p}{e^{p+\triangle p}-e^{p}}
\]
As $\triangle p$ tends to 0, the right-hand side tends to $1/e^{p}$. At the same time
$\triangle q$ tends to 0 and the left-hand side becomes the derivative of ${\rm Blog}$:
\begin{equation}
\label{Blogder}
\frac{d{\rm Blog}(q)}{dq} = \frac{1}{q} \mbox{\femsteg , $q$ nonsingular}
\end{equation}
With the help of (\ref{bicexpder}), (\ref{Blogder}) and the
chain rule we obtain the derivative of the generalized 
power function (3.\ref{genpowerfct})
\begin{equation}
\frac{dq^{r}}{dq} = r\odot q^{r-1} \mbox{\femsteg , $q$ nonsingular}
\end{equation}
For integer $r$ this formula is equivalent to (\ref{pownder}).
One should keep in mind that if $r$ is a positive integer,
the requirement that $q$ is nonsingular can be dropped, a fact
of particular interest in the development of power series, for instance.

We also mention that the complex pair representations
(3.\ref{biccoshc2}), (3.\ref{bicsinhc2}), (3.\ref{biccosc2}), (3.\ref{bicsinc2})
of the hyperbolic and trigonometric cosine and sine functions
yield the expected derivatives:
\begin{eqnarray*}
\label{biccoshder}
\frac{d\cosh p}{dp}&=&\sinh p \\[2mm]
\label{bicsinhder}
\frac{d\sinh p}{dp}&=&\cosh p \\[2mm]
\label{biccosder}
\frac{d\cos p}{dp}&=&-\sin p \\[2mm]
\label{bicsinder}
\frac{d\sin p}{dp}&=&\cos p
\end{eqnarray*}

\subsection{The bicomplex differentiability condition}
The requirement that the two complex pair formulas
(\ref{bicderc21}) and (\ref{bicderc22}), which we found for the bicomplex derivative,
should be equivalent is met if and only if the equations
\begin{eqnarray}
\label{BCR1}
\frac{\partial \phi_{1}}{\partial a} &=& \frac{\partial \phi_{2}}{\partial b} \\ [3mm]
\label{BCR2}
\frac{\partial \phi_{2}}{\partial a} &=& -\frac{\partial \phi_{1}}{\partial b}
\end{eqnarray}
hold in the domain $G\subseteq {\rm B}$ where
$\psi(p)=(\phi_{1}(a,b),\phi_{2}(a,b))$ is differentiable.
These equations represent the 
bicomplex differentiability condition and will
be referred to as the {\em bicomplex
Cauchy-Riemann equations} (abbr. {\em bicomplex CR-equations}).
They are the complexification of the CR-equations (3.\ref{CRequations}).

The validity of (\ref{BCR1}) and (\ref{BCR2}) presupposes that 
the partial derivatives $\frac{\partial \phi_{k}}{\partial a}$
and $\frac{\partial \phi_{k}}{\partial b}$, $k=1,2$, exist, 
i.e. that the functions $\phi_{k}$ are holomorphic in $a$ and $b$.
That this is indeed only a necessary condition for $\psi$
to be holomorphic in $p$ is illustrated by the following example.

\bigskip

\noindent {\bf Example.} Let $\theta$ be given by
\[
\theta(p)=(a^{2},b^{2}) \mbox{\femsteg , $\; p=(a,b)$}
\]
The functions $\phi_{1}(a,b)=a^{2}$ and $\phi_{2}(a,b)=b^{2}$
are holomorphic in $a$ and $b$, but $\theta$
is not holomorphic in $p$, because $\phi_{1}$ and $\phi_{2}$
do not fulfill the
bicomplex CR-equations. These are instead fulfilled by the components of
the square function
\[
\psi(p)=p^{2}=(a^{2}-b^{2},2ab) 
\]
which thus is holomorphic in $p$. $\Box$

\bigskip

\noindent Note that all bicomplex functions introduced in
Chapter 3 fulfill the bicomplex CR-equations.

From the above we conclude:
\newtheorem{lemma}{Lemma}[section]
\begin{lemma}
If the bicomplex function
$\psi(p)=(\phi_{1}(a,b),\phi_{2}(a,b)),\; p=(a,b)$,
is holomorphic in the domain $G\subseteq {\rm B}$,
then the complex functions $\phi_{1}$ and $\phi_{2}$ are holomorphic in $G$
and satisfy the bicomplex CR-equations (\ref{BCR1})--(\ref{BCR2}) there. $\Box$
\end{lemma}

\medskip

\noindent We shall demonstrate that the converse also holds:
\begin{lemma}
If the complex functions $\phi_{1}(a,b)$ and $\phi_{2}(a,b)$ are holomorphic
in the domain $G\subseteq {\rm B}$ and satisfy the bicomplex CR-equations (\ref{BCR1})--(\ref{BCR2}) there,
then the bicomplex function $\psi(p)=(\phi_{1}(a,b),\phi_{2}(a,b)),\; p=(a,b)$,
is holomorphic in $G$.
\end{lemma}

\noindent {\em Proof.} The holomorphism of $\phi_{k}$, $k=1,2$, 
means that they have partial derivatives 
$\frac{\partial \phi_{k}}{\partial a}$ and $\frac{\partial \phi_{k}}{\partial b}$ that
are continuous in $G$. At an arbitrary point $p=(a,b)$ of $G$ we may
therefore by virtue of (3.\ref{differential}) write
\begin{eqnarray*}
\phi_{1}(a+\triangle a,b+\triangle b)-\phi_{1}(a,b)&=&\frac{\partial \phi_{1}}{\partial a}\cdot \triangle a+\frac{\partial \phi_{1}}{\partial b}\cdot \triangle b+\epsilon_{1} \\[3mm]
\phi_{2}(a+\triangle a,b+\triangle b)-\phi_{2}(a,b)&=&\frac{\partial \phi_{2}}{\partial a}\cdot \triangle a+\frac{\partial \phi_{2}}{\partial b}\cdot \triangle b+\epsilon_{2}
\end{eqnarray*}
where $\epsilon_{1}$ and $\epsilon_{2}$ are complex-valued
functions of $\triangle a$ and $\triangle b$ such that they tend to 0 more rapidly
than $(\triangle a,\triangle b)$ in the sense that
$(\epsilon_{1}(\triangle a,\triangle b),0)\,/\,(\triangle a,\triangle b)\rightarrow 0$ and
$(0,\epsilon_{2}(\triangle a,\triangle b))\,/\,(\triangle a,\triangle b)\rightarrow 0$ as
$(\triangle a,\triangle b)\rightarrow 0$.
If we set $\triangle p = (\triangle a,\triangle b)$ and make use of
\[
\psi(p+\triangle p)-\psi(p)\: =\: (\phi_{1}(a +\triangle a,b+\triangle b)-\phi_{1}(a,b)\,,\, \phi_{2}(a+\triangle a,b+\triangle b)-\phi_{2}(a,b))
\]
together with (\ref{BCR1}) and (\ref{BCR2}) we are now able to derive
\[
\psi(p+\triangle p)-\psi(p)\: =\: \left(\frac{\partial \phi_{1}}{\partial a}\,, \frac{\partial \phi_{2}}{\partial a} \right)\odot \triangle p+(\epsilon_{1},\epsilon_{2})
\]
Dividing both sides by $\triangle p$ and splitting the last term into two terms gives
\[
\frac{\psi(p+\triangle p)-\psi(p)}{\triangle p}\: =\: \left(\frac{\partial \phi_{1}}{\partial a}\,, \frac{\partial \phi_{2}}{\partial a} \right)+\frac{(\epsilon_{1},0)}{\triangle p} + \frac{(0,\epsilon_{2})}{\triangle p}
\]
The assumptions about $\epsilon_{1}$ and $\epsilon_{2}$
imply that the last two terms tend to 0 as $\triangle p$ tends to 0,
which shows that $\psi\,'(p)$
exists and equals
\[
\left(\frac{\partial \phi_{1}}{\partial a}\:, \frac{\partial \phi_{2}}{\partial a} \right)
\]
Here, the computation of limits is in each case restricted to nonsingular $\triangle p$
in accordance with (\ref{fctlimitnonsing1}) and (\ref{fctlimitnonsing2}).

\noindent $\Box$

\medskip

\noindent Combining the preceding lemmas results in:
\begin{theorem}
Let $\psi(p)=(\phi_{1}(a,b),\phi_{2}(a,b)),\: p=(a,b)$, be
a bicomplex function. Then $\psi$ is holomorphic in the domain $G\subseteq {\rm B}$ if and only
if the complex functions $\phi_{1}(a,b)$ and $\phi_{2}(a,b)$ are holomorphic
in $G$ as well as fulfill the bicomplex CR-equations (\ref{BCR1})--(\ref{BCR2}) there. $\Box$
\end{theorem}
The bicomplex CR-equations are essential to the
theory under consideration. They can be combined into a single equation, namely
\begin{equation}
\label{BCR3}
\frac{\partial \psi}{\partial a} + j \odot \frac{\partial \psi}{\partial b} = 0 \mbox{\femsteg , $j=(0,1)$}
\end{equation}
where $\frac{\partial\psi}{\partial a}$ and $\frac{\partial\psi}{\partial b}$ denote
\begin{equation}
\label{psiderab}
\frac{\partial\psi}{\partial a} = \left(\frac{\partial \phi_{1}}{\partial a}\,, \frac{\partial \phi_{2}}{\partial a} \right) \mbox{\tresteg , \tresteg}
\frac{\partial\psi}{\partial b} = \left(\frac{\partial \phi_{1}}{\partial b}\,, \frac{\partial \phi_{2}}{\partial b} \right)
\end{equation}
Formula (\ref{BCR3}) is the generalization of
$\frac{\partial f}{\partial x} + i\cdot \frac{\partial f}{\partial y} = 0$
satisfied by a complex holomorphic function.

\subsection{Higher order derivatives}
Next we investigate the existence of higher order derivatives of $\psi$.
Application of (\ref{bicderc21}) $n$ times to $\psi$ gives
\begin{equation}
\label{bicdern1}
\frac{d^{n}\psi}{dp^{n}} = \left(\frac{\partial^{n} \phi_{1}}{\partial a^{n}}\,, \frac{\partial^{n} \phi_{2}}{\partial a^{n}} \right) \mbox{\femsteg , $\; n\geq 0$}
\end{equation}
Alternatively, the ${\rm n}^{\rm th}$ derivative can be computed
by repeated application of (\ref{bicderc22}):
\begin{eqnarray}
\label{bicdern2}
\frac{d^{n}\psi}{dp^{n}} &=& \left((-1)^{\left\lfloor \frac{n}{2} \right\rfloor}\frac{\partial^{n} \phi_{2}}{\partial b^{n}}\;, (-1)^{\left\lfloor \frac{n+1}{2} \right\rfloor}\frac{\partial^{n} \phi_{1}}{\partial b^{n}} \right) \mbox{\tresteg , for odd $n \geq 1$} \\ [3mm]
\label{bicdern3}
\frac{d^{n}\psi}{dp^{n}} &=& \left((-1)^{\left\lfloor \frac{n}{2} \right\rfloor}\frac{\partial^{n} \phi_{1}}{\partial b^{n}}\;, (-1)^{\left\lfloor \frac{n+1}{2} \right\rfloor}\frac{\partial^{n} \phi_{2}}{\partial b^{n}} \right) \mbox{\tresteg , for even $n\geq 0$}
\end{eqnarray}

\smallskip

\noindent For the proof of the subsequent theorem we note that (\ref{bicdern1}) can be rewritten
\begin{equation}
\label{bicdern4}
\frac{d^{n}\psi}{dp^{n}} = \frac{\partial^{n} \psi}{\partial a^{n}} \mbox{\femsteg , $\; n\geq 0$}
\end{equation}
The left-hand side is conveniently abbreviated by
\begin{equation}
\label{bicdern5}
\psi^{(n)} = \frac{d^{n}\psi}{dp^{n}} \mbox{\femsteg , $\; n\geq 0$}
\end{equation}
It is understood that $\psi^{(0)}=\psi$.
\begin{theorem}
A bicomplex function
$\psi(p)=(\phi_{1}(a,b),\phi_{2}(a,b)),\: p=(a,b)$, 
that is holomorphic in a domain $G\subseteq {\rm B}$ possesses
derivatives of all orders in $G$, each one holomorphic in $G$.
\end{theorem}
\noindent {\em Proof.} Let the 
${\rm n}^{\rm th}$ derivative of $\psi$ be given by (\ref{bicdern1})
and designated by $\psi^{(n)}$. According to Lemma 4.2
$\psi^{(n)}$ is holomorphic in $G$ if $\frac{\partial^{n} \phi_{1}}{\partial a^{n}}$
and $\frac{\partial^{n} \phi_{2}}{\partial a^{n}}$ are
holomorphic in $G$ and satisfy the bicomplex CR-equations.

By Lemma 4.1 the holomorphism of $\psi$ in $G$ implies that
$\phi_{1}$ and $\phi_{2}$ are holomorphic in $G$.
From complex analysis we know that
a complex holomorphic function has derivatives of all orders,
which by themselves are holomorphic. 
It follows that the holomorphism of $\phi_{1}$ and $\phi_{2}$
is inherited by $\frac{\partial^{n} \phi_{1}}{\partial a^{n}}$
and $\frac{\partial^{n} \phi_{2}}{\partial a^{n}}$.

To prove that $\psi^{(n)}$ satisfies the bicomplex CR-equations we
apply (\ref{BCR3}), instead of (\ref{BCR1})--(\ref{BCR2}), and establish
\begin{equation}
\label{BCR3n}
\frac{\partial \psi^{(n)}}{\partial a} + j \odot \frac{\partial \psi^{(n)}}{\partial b} = 0
\end{equation}
using mathematical induction. For $n=0$ the formula is
identical with (\ref{BCR3}), which holds on account of $\psi$'s
holomorphism. For the validation of the step we assume
(\ref{BCR3n}) and calculate:
\begin{eqnarray*}
& & \frac{\partial \psi^{(n+1)}}{\partial a} + j \odot \frac{\partial \psi^{(n+1)}}{\partial b} = 0 \\[2mm]
&\equiv& \mbox{\femsteg $\lbrace \psi^{(n+1)}=\frac{\partial \psi^{(n)}}{\partial a}$, see (\ref{bicdern4}) and (\ref{bicdern5})$\rbrace$} \\[2mm]
& & \frac{\partial}{\partial a}\left(\frac{\partial \psi^{(n)}}{\partial a}\right) + j \odot \frac{\partial}{\partial b}\left(\frac{\partial \psi^{(n)}}{\partial a}\right) = 0 \\[2mm]
&\equiv& \mbox{\femsteg $\lbrace$induction hypothesis$\rbrace$} \\[2mm]
& & \frac{\partial}{\partial a}\left(- j \odot \frac{\partial \psi^{(n)}}{\partial b}\right) + j \odot \frac{\partial}{\partial b}\left(\frac{\partial \psi^{(n)}}{\partial a}\right) = 0 \\[2mm]
&\equiv& \mbox{\femsteg $\lbrace$taking derivatives in different order, see below$\rbrace$} \\[2mm]
& & - j \odot \frac{\partial}{\partial b}\left(\frac{\partial \psi^{(n)}}{\partial a}\right) + j \odot \frac{\partial}{\partial b}\left(\frac{\partial \psi^{(n)}}{\partial a}\right) = 0 \\[2mm]
&\equiv& \mbox{\femsteg $\lbrace$identity$\rbrace$} \\[2mm]
& & {\rm true}
\end{eqnarray*}
In the penultimate step we changed the order
of the operators $\frac{\partial}{\partial a}$ and
$\frac{\partial}{\partial b}$ when applied to $\psi^{(n)}$.
The interchange is allowed because $\psi^{(n)}$ consists
of holomorphic component functions and
\[
\frac{\partial}{\partial a}\left(\frac{\partial \omega}{\partial b}\right) = \frac{\partial}{\partial b}\left(\frac{\partial \omega}{\partial a}\right)
\]
holds for an arbitrary complex function $\omega: (a,b) \rightarrow \omega(a,b)$
that is holomorphic in both $a$ and $b$.

\noindent $\Box$

\subsection{${\bf R^{4}}$-representations of the bicomplex derivative}
In order to get another handle to bicomplex derivatives
we rewrite (\ref{bicderc21}) and (\ref{bicderc22}) in ${\rm R^{4}}$.
Assume that the full representation of the
function $\psi$ is
\begin{eqnarray*}
\psi(p) &=& (\phi_{1}(a,b),\phi_{2}(a,b)) \\ [2mm]
(\phi_{1}(a,b),\phi_{2}(a,b)) &=& (\psi_{1}(x,y,z,u)+i\psi_{2}(x,y,z,u)\,, \\
& &\; \psi_{3}(x,y,z,u)+i\psi_{4}(x,y,z,u)) \\ [2mm]
p=(a,b)&,&\mbox{$a=x+iy$ and $b=z+iu$}
\end{eqnarray*}
On account of (3.\ref{dthetadainx}) and (3.\ref{dthetadainy}) the partial
derivatives $\frac{\partial \phi_{1}}{\partial a}$ and
$\frac{\partial \phi_{1}}{\partial b}$ are given by:
\begin{eqnarray*}
\frac{\partial \phi_{1}}{\partial a} &=& \partial_{x}\psi_{1}+i\partial_{x}\psi_{2} \\ [2mm]
\frac{\partial \phi_{1}}{\partial a} &=& \partial_{y}\psi_{2}-i\partial_{y}\psi_{1} \\ [4mm]
\frac{\partial \phi_{1}}{\partial b} &=& \partial_{z}\psi_{1}+i\partial_{z}\psi_{2} \\ [2mm]
\frac{\partial \phi_{1}}{\partial b} &=& \partial_{u}\psi_{2}-i\partial_{u}\psi_{1}
\end{eqnarray*}
Analogous formulas hold for 
$\frac{\partial \phi_{2}}{\partial a}$ and
$\frac{\partial \phi_{2}}{\partial b}$. Substitution of the complex
derivatives into (\ref{bicderc21}) and (\ref{bicderc22}) yields
four ${\rm R^{4}}$-representations of $\frac{d\psi}{dp}$:
\begin{eqnarray}
\label{bicderr41}
\frac{d\psi}{dp}&=&(\partial_{x}\psi_{1}+i\partial_{x}\psi_{2}\:,\:\partial_{x}\psi_{3}+i\partial_{x}\psi_{4}) \\ [2mm]
\label{bicderr42}
\frac{d\psi}{dp}&=&(\partial_{y}\psi_{2}-i\partial_{y}\psi_{1}\:,\:\partial_{y}\psi_{4}-i\partial_{y}\psi_{3}) \\ [2mm]
\label{bicderr43}
\frac{d\psi}{dp}&=&(\partial_{z}\psi_{3}+i\partial_{z}\psi_{4}\:,\:-\partial_{z}\psi_{1}-i\partial_{z}\psi_{2}) \\ [2mm]
\label{bicderr44}
\frac{d\psi}{dp}&=&(\partial_{u}\psi_{4}-i\partial_{u}\psi_{3}\:,\:-\partial_{u}\psi_{2}+i\partial_{u}\psi_{1})
\end{eqnarray}
Formulas (\ref{bicderr41}) and (\ref{bicderr42}) were obtained from (\ref{bicderc21}), (\ref{bicderr43}) and
(\ref{bicderr44}) from (\ref{bicderc22}).

\bigskip

\noindent {\bf Example.} The ${\rm R^{4}}$-representation
of the square function $\psi(p)=p^{2}$ is
\[
\psi(p)\;=\; (x^{2}-y^{2}-z^{2}+u^{2}+i(2xy-2zu)\, , \, 2xz-2yu+i(2yz+2xu))
\]
Each of (\ref{bicderr41})--(\ref{bicderr44}) yields the same derivative
\[
\frac{d\psi}{dp}=(2x+i2y,2z+i2u)
\]
This means that $\frac{dp^{2}}{dp}=2p$, as expected. $\Box$

\bigskip

\noindent The formulas (\ref{bicderr41})--(\ref{bicderr44}) should be equivalent. The equivalence
of (\ref{bicderr41}) and (\ref{bicderr42}), in particular, implies the simultaneous validity of the equations
\begin{eqnarray*}
\partial_{x}\psi_{1}&=&\partial_{y}\psi_{2} \\ [1mm]
\partial_{x}\psi_{2}&=&-\partial_{y}\psi_{1} \\ [4mm]
\partial_{x}\psi_{3}&=&\partial_{y}\psi_{4} \\ [1mm]
\partial_{x}\psi_{4}&=&-\partial_{y}\psi_{3}
\end{eqnarray*}
They express that $\psi_{1}+i\psi_{2}$
and $\psi_{3}+i\psi_{4}$ are holomorphic with respect to $a=x+iy$.
By the same token, the equivalence of (\ref{bicderr43}) and (\ref{bicderr44})
amounts to the simultaneous validity of
\begin{eqnarray*}
\partial_{z}\psi_{3}&=&\partial_{u}\psi_{4} \\ [1mm]
\partial_{z}\psi_{4}&=&-\partial_{u}\psi_{3} \\ [4mm]
-\partial_{z}\psi_{1}&=&-\partial_{u}\psi_{2} \\ [1mm]
-\partial_{z}\psi_{2}&=&\partial_{u}\psi_{1}
\end{eqnarray*}
which express that $\psi_{1}+i\psi_{2}$ and $\psi_{3}+i\psi_{4}$ are holomorphic
with respect to $b=z+iu$. The additional
equivalence of (\ref{bicderr41}) and (\ref{bicderr43}) means that the equations
\begin{eqnarray*}
\partial_{x}\psi_{1}&=&\partial_{z}\psi_{3} \\ [1mm]
\partial_{x}\psi_{2}&=&\partial_{z}\psi_{4} \\ [4mm]
\partial_{x}\psi_{3}&=&-\partial_{z}\psi_{1} \\ [1mm]
\partial_{x}\psi_{4}&=&-\partial_{z}\psi_{2}
\end{eqnarray*}
all hold and validate the bicomplex differentiability condition
(\ref{BCR1})--(\ref{BCR2}). Thus, in a straightforward way
formulas (\ref{bicderr41})--(\ref{bicderr44}) exhibit both the
holomorphism of $\psi$ and its component functions
$\psi_{1}+i\psi_{2}$ and $\psi_{3}+i\psi_{4}$.

\newpage
\setcounter{equation}{0}

\section{Bicomplex integration}
In this chapter we define the concept of a bicomplex
integral as a line integral. We inquire under what
conditions it is path independent using an argument analogous
to the one given by Ahlfors for complex functions. By integration
we show that a holomorphic bicomplex function can be written
as a Taylor series. To prove Cauchy's theorem in ${\rm B}$
we lean on a generalization of Green's theorem.
To express Cauchy's integral formula in ${\rm B}$ we first
derive the bicomplex {\em twining number},
which generalizes the complex winding number. 
We are then able to prove the formula with the help of a Taylor expansion.

\subsection{Definition of the line integral}
We consider the integration of the bicomplex function
\begin{eqnarray}
\label{psidef2}
\psi(p) &=& (\phi_{1}(a,b),\phi_{2}(a,b)) \\ [2mm]
(\phi_{1}(a,b),\phi_{2}(a,b)) &=& (\psi_{1}(x,y,z,u)+i\psi_{2}(x,y,z,u), \nonumber \\
& &\; \psi_{3}(x,y,z,u)+i\psi_{4}(x,y,z,u)) \nonumber \\ [2mm]
p=(a,b)&,&\mbox{$a=x+iy$ and $b=z+iu$} \nonumber
\end{eqnarray}
For the time being $\psi$'s domain of definition is left anonymous.
We assume that $\phi_{1}$ and $\phi_{2}$ are holomorphic in $a$ and $b$, thereby
ensuring that $\psi$ is continuous as demonstrated in Section 3.1.

The basic bicomplex integral is essentially isomorphic to the complex integral.
It is to be understood as a {\em line integral} that is evaluated with
respect to some four-dimensional curve $\Gamma$ in ${\rm B}$.
More specifically, the concept to be defined is
\[
\int_{\Gamma} \psi(p)\odot dp \mbox{\femsteg , $dp=(da,db)$}
\]
Henceforth, we shall choose the curve $\Gamma$ so that it is piecewise continuously
differentiable in ${\rm B}$ and has the parametric equation
\begin{equation}
\label{Gamma1}
\Gamma:\mbox{\tresteg} p=p(t) \mbox{\femsteg , $p(t)=(a(t),b(t))\;\;$ for $r\leq t \leq s$}
\end{equation}
The argument $t$ is real. One can look at $\Gamma$ as a curve made up of
two component curves $\gamma_{1}$ and $\gamma_{2}$
in ${\rm C}$
\begin{equation}
\label{Gamma2}
\Gamma=(\gamma_{1},\gamma_{2})
\end{equation}
whose parametric equations are
\begin{eqnarray*}
&&\gamma_{1}:\mbox{\tresteg} a=a(t) \mbox{\femsteg , $a(t)=x(t)+iy(t)\;$ for $r\leq t \leq s$} \\ [2mm]
&&\gamma_{2}:\mbox{\tresteg} b=b(t) \:\mbox{\femsteg , $b(t)=z(t)+iu(t)\;$ for $r\leq t \leq s$}
\end{eqnarray*}
As definition of the line integral of $\psi(p)$ extended over the curve
$\Gamma$ we then take
\begin{equation}
\label{bicintdef}
\int_{\Gamma} \psi(p)\odot dp = \int_{r}^{s} \psi(p(t))\odot p\,'(t) dt
\end{equation}
Because $\psi$ is continuous, $\psi(p(t))$ at the right-hand side is also continuous.
If $p\,'(t)$ is discontinuous at some points the integration has to be performed
in subintervals of $[r,s]$ in the normal manner.

The left-hand side of (\ref{bicintdef}) can be reformulated in the following two ways.
Firstly, bicomplex multiplication of $\psi=(\phi_{1},\phi_{2})$ and
$dp=(da,db)$ yields
\begin{equation}
\label{bicint1}
\int_{\Gamma} \psi(p)\odot dp \;=\; \left(\int_{\Gamma} \phi_{1}\cdot da - \phi_{2}\cdot db , \int_{\Gamma} \phi_{2}\cdot da + \phi_{1}\cdot db \right)
\end{equation}
Secondly, since the differential $dp$ is expressible as
\[
dp=(da,0)+j\odot (db,0) \mbox{\femsteg , $j=(0,1)$} 
\]
we have
\begin{equation}
\label{bicint2}
\int_{\Gamma} \psi(p)\odot dp \;=\; \int_{\Gamma} \psi(p)\odot (da,0) + j\odot \psi(p) \odot (db,0)
\end{equation}
The bicomplex differentials $(da,0)$ and $(db,0)$ have the same
properties as the complex differentials $da$ and $db$. Formulas (\ref{bicint1}) and (\ref{bicint2}) will
be needed in the next two sections.

We shall also need integrals of the type $\int_{\Gamma} \psi(p)\, \|dp \|$
involving the norm (2.\ref{normb}). With $\Gamma$ specified by (\ref{Gamma1})
we define
\begin{equation}
\label{bicintnormdef}
\int_{\Gamma} \psi(p)\,\| dp\| \;=\; \int_{r}^{s} \psi(p(t))\,\| p\,'(t)\| dt
\end{equation}
The special case $\psi(p)=1$ gives us a suitable definition of the
{\em length} of $\Gamma$:
\begin{equation}
\label{curvelengthdef}
\int_{\Gamma} \| dp\| \;=\; \mbox{length of $\Gamma$}
\end{equation}
By regarding the integral (\ref{bicintdef}) as the limit of a Riemann sum
we obtain by applying the triangle inequality
\begin{equation}
\label{bicintinequality}
\left\| \int_{\Gamma} \psi(p)\odot dp \right\| \;\leq\; \int_{\Gamma} \|\psi(p)\|\cdot \| dp\|
\end{equation}
If the integral at the left-hand side has a finite value 
there exists a real constant $M$ such that 
$\|\psi(p)\| \leq M$ for all $p$ on
$\Gamma$. We therefore get the estimate
\begin{equation}
\label{bicintestimate}
\left\| \int_{\Gamma} \psi(p)\odot dp \right\| \;\leq\; M\cdot L
\end{equation}
where $L$ stands for the length of $\Gamma$.

Formulas (\ref{bicintnormdef})--(\ref{bicintestimate}) are direct
generalizations of the corresponding ones in complex analysis~\cite{lehto}.

\subsection{Path independence and the primitive function}
We ask ourselves under what condition the bicomplex
line integral is independent of the path of integration
in a domain $G\subseteq {\rm B}$. With respect to definition (\ref{bicintdef}) this
means that if $\Gamma$ is free to vary in $G$ the integral
depends only on the end points of $\Gamma$. In particular,
if $\Gamma$ is a closed curve the integral reduces to zero.

Our reasoning will be based on formula (\ref{bicint2}). The fact that its right-hand side is
of the form
\[
\int_{\Gamma} \omega_{1}\odot (da,0) + \omega_{2} \odot (db,0)
\]
gives relevance to the following theorem.
\begin{theorem}
Let $\omega_{1}$ and $\omega_{2}$ be two bicomplex functions
of $p=(a,b)$ such that they are continuous in the domain $G\subseteq {\rm B}$.
Then the line integral $\int_{\Gamma} \omega_{1}\odot (da,0) + \omega_{2} \odot (db,0)$,
defined in $G$, depends only on the end points of the curve $\Gamma$
if and only if there exists a function $\Omega: (a,b) \rightarrow \Omega(a,b)$
such that $\frac{\partial \Omega}{\partial a} = \omega_{1}$ and 
$\frac{\partial \Omega}{\partial b} = \omega_{2}$.
\end{theorem}
\noindent {\it Proof.} Omitted, because it is entirely isomorphic to the proof of the corresponding theorem for complex functions
as given by Ahlfors~\cite{ahlfors}, pp. 106--107.

\noindent $\Box$

\bigskip

\noindent We focus on the
evaluation of the right-hand side of (\ref{bicint2}) in a domain
$G\subseteq {\rm B}$. According to the theorem above 
the integral depends only on the end points of $\Gamma$
if and only if there exists a function $\Psi:p \rightarrow \Psi(p)$, where $p=(a,b)$,
such that
\begin{eqnarray}
\label{Psidera}
& &\frac{\partial \Psi(p)}{\partial a} = \psi(p) \\ [2mm]
\label{Psiderb}
& &\frac{\partial \Psi(p)}{\partial b} = j\odot\psi(p)
\end{eqnarray}
Multiplying (\ref{Psiderb}) by $j$ and adding it to (\ref{Psidera}) yields
\[
\frac{\partial \Psi}{\partial a} + j \odot \frac{\partial \Psi}{\partial b} = 0
\]
But this is the bicomplex differentiability condition (4.\ref{BCR3}), which means
that $\Psi$ is a holomorphic function. Formula (\ref{Psidera}) further tells us that
$\psi$ is the derivative of $\Psi$. (Note that 
$\frac{d\Psi}{dp} = \frac{\partial \Psi}{\partial a}$ on account of (4.\ref{bicdern4}), and that
$\psi$, too, is holomorphic on account of Theorem 4.2.)
Hence, we have established:
\begin{theorem}
Let $\psi$ be a bicomplex function of the form {\rm (\ref{psidef2})} that
is continuous in a domain $G\subseteq {\rm B}$. Then for an arbitrary curve
$\Gamma$ in $G$ the integral $\int_{\Gamma} \psi(p)\odot dp$
depends only on the end points of $\Gamma$ if and only if
$\psi$ is the derivative of a holomorphic function $\Psi$ in $G$. $\Box$
\end{theorem}
The function $\Psi$ is uniquely determined up to an
additive constant and may be chosen as 
primitive function (integral function)
of $\psi$. It enables one to evaluate a line integral
of $\psi$ extended over a curve with end points
$p_{1}$ and $p_{2}$ according to
\[
\int_{p_{1}}^{p_{2}} \psi(p)\odot dp \;=\; [\Psi(p)]\begin{array}{c} 
                                                   {\scriptstyle \!\!\!p_{2}}\\ [-1.5mm]
                                                   {\scriptstyle \!\!\!p_{1}}
                                              \end{array} \;=\; \Psi(p_{2})-\Psi(p_{1})
\]
Thanks to path independence the omission of any reference
to the actual curve is justified here.

\subsection{Taylor series}
We shall show that a holomorphic bicomplex function can be developed into a
Taylor series. Our starting point is the formula for
computing the derivative of the product of two bicomplex functions:
\[
(\theta \odot \varphi)\,'(p) \;=\; \theta\,'(p)\odot \varphi(p) + \theta(p)\odot \varphi\,'(p)
\]
In a domain $G$ where both $\theta$ and $\varphi$ are holomorphic the function
$(\theta \odot \varphi)\,'$ certainly has a primitive function, or $\theta \odot \varphi$.
As a result, if we subject both sides to the integral operator
$\int_{p_{0}}^{p}$ and transfer a term between them we obtain
the bicomplex version of the formula for integration by parts:
\[
\int_{p_{0}}^{p} \theta(p)\odot \varphi\,'(p)\odot dp \;=\; [\theta(p) \odot \varphi(p)]\begin{array}{c} 
                                                                                                 {\scriptstyle \!\!\!\!\!p}\\ [-1.5mm]
                                                                                                 {\scriptstyle \!\!\!p_{0}}
                                                                                            \end{array} - \int_{p_{0}}^{p}\theta\,'(p)\odot \varphi(p)\odot dp
\]
It is possible to compute the integrals along any curve in $G$
starting at $p_{0}$ and ending at $p$.

Let the function $\psi$ be holomorphic in $G$ so that
it possesses derivatives  of all
orders there due to Theorem 4.2. Then, because integration by parts is available,
we can expand $\psi$ into the following Taylor polynomial and remainder term about $p_{0}$:
\begin{eqnarray}
\label{taylorseries}
\psi(p) &=& \left(\sum k: 0\leq k \leq n-1: \; \frac{\psi^{(k)}(p_{0})}{k!} \odot (p-p_{0})^{k} \right) + R_{n}(p,p_{0}) \\ [2mm]
R_{n}(p,p_{0}) &=& \frac{1}{(n-1)!} \int_{p_{0}}^{p} \psi^{(n)}(q)\odot (p-q)^{n-1}\odot dq \nonumber
\end{eqnarray}
In order to estimate the remainder term $R_{n}(p,p_{0})$ we study
its behaviour in a $p_{0}$-centered ball that lies in $G$ and is
specified by $\|q-p_{0}\| \leq \|p-p_{0}\|$. We assume that the integral
is evaluated along a curve whose length 
$\int_{p_{0}}^{p} \| dq\|$ equals $L$. The derivative
$\psi^{(n)}(q)$ satisfies $\|\psi^{(n)}(q)\| \leq M$
for a suitable $M$. Furthermore, using the triangle inequality, in particular,
we deduce: 
\[
\|p-q\| \;\:=\;\: \|(p-p_{0})-(q-p_{0})\| \:\;\leq\;\: \|p-p_{0}\| + \|q-p_{0}\| \:\;\leq\;\: 2\|p-p_{0}\|
\]
Hence, we estimate
\begin{equation}
\label{remestimate}
\|R_{n}(p,p_{0})\| \;\leq \; \frac{M\cdot L\cdot (2\|p-p_{0}\|)^{n-1}}{(n-1)!}
\end{equation}
Because $\|R_{n}(p,p_{0})\|$ approaches 0 as $n \rightarrow \infty$
the infinite Taylor series converges.

\subsection{Cauchy's theorem in bicomplex space}
For the purpose of generalizing Cauchy's theorem for a bicomplex
function we recall one of its formulations in the complex plane.
\begin{theorem}
Let $A$ be a simply connected domain in ${\rm C}$.
If the complex function $f(a)=g(x,y)+ih(x,y)$, $a=x+iy$, is holomorphic
in $A$, then
\[
\int_{\gamma} f(a)\cdot da =0
\]
for any closed curve $\gamma$ in $A$. $\Box$
\end{theorem}

\noindent The domain $A$ being simply connected means that
it is without "holes" or, more precisely, that its complement with
respect to the extended plane is connected. It is well-known that
the theorem is provable with the help of {\em Green's theorem},
provided one assumes that $g(x,y)$ and $h(x,y)$ possess continuous
partial derivatives in $x$ and $y$~\cite{apostol, seaborn}. Although the proof is based
on unnecessarily strong assumptions, we take it as starting point because
conducting it in ${\rm B}$ clarifies important concepts. 
In ${\rm R^{2}}$ Green's theorem is written
\begin{equation}
\label{greenr2}
\int_{\gamma} g\,dx + h\,dy = \int_{S(\gamma)} \left(\frac{\partial h}{\partial x} - \frac{\partial g}{\partial y}\right)dx\,dy
\end{equation}
The left-hand side is a line integral extended over a closed curve $\gamma$,
the right-hand side a surface integral extended over $S(\gamma)$, which is
the domain of the $x$-$y$ plane bounded by $\gamma$.

Transferring the above to ${\rm B}$ first demands that
we get some basic understanding of {\em complex} line
and surface integrals of the form
\begin{eqnarray}
\label{lineintc2}
& &\int_{\Gamma} \phi \cdot da + \omega \cdot db \\ [2mm]
\label{surfintc2}
& &\int_{S(\Gamma)} \phi \cdot da \cdot db
\end{eqnarray}
where
\begin{itemize}
\item $\Gamma$ is a closed curve in ${\rm B}$ with the parametric
      equation (\ref{Gamma1}).
\item $S(\Gamma)$ is an arbitrary surface in ${\rm B}$ bounded by $\Gamma$.
\item $\phi$ and $\omega$ are complex functions of $a$ and $b$,
      i.e. of the type ${\rm B} \rightarrow {\rm C}$ (or ${\rm C}^{2} \rightarrow {\rm C}$, if one so wants).
\end{itemize}
Integrals like (\ref{lineintc2}) appeared already in formula (\ref{bicint1}). 
For the present purpose it suffices to define the integrals (\ref{lineintc2}) and (\ref{surfintc2})
by reducing them to real line and surface integrals, which we consider
well-understood. Let $\phi$ and $\omega$ more specifically be
given by
\begin{eqnarray*}
\phi(a,b) &=& \xi_{1}(x,y,z,u)+i\,\xi_{2}(x,y,z,u) \\ [2mm]
\omega(a,b) &=& \eta_{1}(x,y,z,u)+i\,\eta_{2}(x,y,z,u)
\end{eqnarray*}
where $\xi_{1}$, $\xi_{2}$, $\eta_{1}$, $\eta_{2}$ are continuous functions of
the type ${\rm R}^{4} \rightarrow {\rm R}$. With
$da=dx+idy$ and $db=dz+idu$ we obtain
\begin{eqnarray}
\label{lineintr4}
\int_{\Gamma} \phi \cdot da + \omega \cdot db &=&\int_{\Gamma} (\xi_{1}\, dx - \xi_{2}\, dy + \eta_{1}\, dz - \eta_{2}\, du) \\ [1mm]
                                                    & &\!\!\!\!\!\!+\; i \int_{\Gamma} (\xi_{2}\, dx + \xi_{1}\, dy + \eta_{2}\, dz + \eta_{1}\, du) \nonumber \\ [3mm]
\label{surfintr4}
\int_{S(\Gamma)} \phi \cdot da \cdot db &=& \int_{S(\Gamma)} (\xi_{1}\, dx dz - \xi_{1}\, dy du - \xi_{2}\, dx du - \xi_{2}\, dy dz) \\ [1mm]
                                              & &\!\!\!\!\!\!+\; i \int_{S(\Gamma)} (\xi_{1}\, dx du + \xi_{1}\, dy dz + \xi_{2}\, dx dz - \xi_{2}\, dy du) \nonumber
\end{eqnarray}
We generalize Green's theorem (\ref{greenr2}) by substituting
complex functions and variables for real ones as well as the
four-dimensional curve $\Gamma$ for the two-dimensional $\gamma$.
We also have to require that the functions are holomorphic.
\begin{theorem}
Let $\phi_{1}$ and $\phi_{2}$ be two complex functions
of $a$ and $b$ such that they are holomorphic in a domain
$G\subseteq {\rm B}$. Also let $S(\Gamma)$ be an orientable surface
in $G$ bounded by the closed curve $\Gamma$. Then
\begin{equation}
\label{greenc2}
\int_{\Gamma} \phi_{1}\cdot da + \phi_{2}\cdot db = \int_{S(\Gamma)} \left(\frac{\partial \phi_{2}}{\partial a} - \frac{\partial \phi_{1}}{\partial b}\right)\cdot da\cdot db
\end{equation}
\end{theorem}
\noindent {\it Proof.} Uses the generalized Stokes theorem, see Appendix.

\noindent $\Box$

\bigskip

\noindent With this result the proof of Cauchy's theorem in ${\rm B}$ comes out in few
lines.
\begin{theorem}
Let the bicomplex function
$\psi(p) = (\phi_{1}(a,b),\phi_{2}(a,b))$, $p=(a,b)$,
be holomorphic in a domain $G\subseteq {\rm B}$. Then
\[
\int_{\Gamma} \psi(p)\odot dp = 0
\]
for any closed curve $\Gamma$ that is the boundary of
an orientable surface $S(\Gamma)$ in $G$.
\end{theorem}
{\em Proof.} 
\begin{eqnarray*}
& & \int_{\Gamma} \psi(p)\odot dp \\[2mm]
&=& \mbox{\femsteg $\lbrace$(\ref{bicint1})$\rbrace$} \\[2mm]
& & \left(\int_{\Gamma} \phi_{1}\cdot da - \phi_{2}\cdot db , \int_{\Gamma} \phi_{2}\cdot da + \phi_{1}\cdot db \right) \\[2mm]
&=& \mbox{\femsteg $\lbrace$$S(\Gamma)$ is an orientable surface in $G$; (\ref{greenc2}) twice$\rbrace$} \\[2mm]
& & \left(\int_{S(\Gamma)} \left(-\frac{\partial \phi_{2}}{\partial a} - \frac{\partial \phi_{1}}{\partial b}\right)\cdot da\cdot db\;,\;\int_{S(\Gamma)} \left(\frac{\partial \phi_{1}}{\partial a} - \frac{\partial \phi_{2}}{\partial b}\right)\cdot da\cdot db \right) \\[2mm]
&=& \mbox{\femsteg $\lbrace$$\psi$ fulfills the bicomplex CR-equations (4.\ref{BCR1})--(4.\ref{BCR2}) in $G$$\rbrace$} \\[2mm]
& & 0
\end{eqnarray*}
$\Box$

\bigskip

\noindent {\bf Example.} By way of illustration we shall evaluate
the integral of $e^{p}$ along the closed curve
\[
\Gamma:\mbox{\tresteg} p(t)=(Re^{it},Re^{it}) \mbox{\femsteg , $0\leq t\leq 2\pi$}
\]
where $R$ is a positive real constant. On the $x$-$y$- and $z$-$u$-planes
$\Gamma$ is projected as two origo-centered circles of radius $R$.
The component curves of $\Gamma$ have the parametric equations
\begin{eqnarray*}
&&\gamma_{1}:\mbox{\tresteg} a(t)=Re^{it} \mbox{\femsteg , $0\leq t\leq 2\pi$} \\ [2mm]
&&\gamma_{2}:\mbox{\tresteg} b(t)=Re^{it} 
\end{eqnarray*}
The outcome of the following calculation complies with Cauchy's theorem.
\begin{eqnarray*}
& & \int_{\Gamma} e^{p}\odot dp \\[2mm]
&=& \mbox{\femsteg $\lbrace$bicomplex line integral, $0\leq t\leq 2\pi$$\rbrace$} \\[2mm]
& & \int_{0}^{2\pi} e^{p(t)}\odot p\,'(t)\, dt \\[2mm]
&=& \mbox{\femsteg $\lbrace$definition of $\Gamma$ and $e^{p}$$\rbrace$} \\[2mm]
& & \int_{0}^{2\pi} (e^{Re^{it}}\!\cdot \cos(Re^{it})\:,\: e^{Re^{it}}\!\cdot \sin(Re^{it}))\odot (iRe^{it},iRe^{it})\,dt \\[2mm]
&=&  \mbox{\femsteg $\lbrace$multiplication, distribution properties$\rbrace$} \\[6mm]
& & \left(\int_{0}^{2\pi} (e^{Re^{it}}\!\cdot \cos(Re^{it})\:-\: e^{Re^{it}}\!\cdot \sin(Re^{it}))\cdot iRe^{it}\,dt\;,\right. \\[1mm]
& & \left.\;\; \int_{0}^{2\pi} (e^{Re^{it}}\!\cdot \sin(Re^{it})\:+\: e^{Re^{it}}\!\cdot \cos(Re^{it}))\cdot iRe^{it}\,dt \right) \\[2mm]
&=&  \mbox{\femsteg $\lbrace$definition of the curves $\gamma_{1}$ and $\gamma_{2}$$\rbrace$} \\ [2mm]
& & \left(\int_{0}^{2\pi} (e^{a(t)}\!\cdot \cos a(t)\:-\: e^{a(t)}\!\cdot \sin a(t))\cdot a\,'(t)\,dt\;,\right. \\[1mm]
& & \left.\;\; \int_{0}^{2\pi} (e^{b(t)}\!\cdot \sin b(t)\:+\: e^{b(t)}\!\cdot \cos b(t))\cdot b\,'(t)\,dt \right) \\[2mm]
&=&  \mbox{\femsteg $\lbrace$definition of line integral in ${\rm C}$$\rbrace$} \\ [2mm]
& & \left(\int_{\gamma_{1}}(e^{a}\!\cdot \cos a\:-\: e^{a}\!\cdot \sin a)\cdot \,da\;,\;\int_{\gamma_{2}} (e^{b}\!\cdot \sin b\:+\: e^{b}\!\cdot \cos b)\cdot \,db \right) \\[2mm]
&=&  \mbox{\femsteg $\lbrace$$\gamma_{1}$ and $\gamma_{2}$ are closed curves; Cauchy's theorem in ${\rm C}$$\rbrace$} \\ [2mm]
& & 0
\end{eqnarray*}
$\Box$

\subsection{The twining number}
In the complex plane the so-called {\em winding number} $u(\gamma,a_{0})$
of a closed curve $\gamma$ with respect to a point $a_{0}$ is defined by
\[
u(\gamma,a_{0}) = \frac{1}{2\pi i} \int_{\gamma} \frac{da}{a-a_{0}} \mbox{\femsteg , $a_{0}\not\in \gamma$}
\]
We call the corresponding bicomplex concept the {\em twining number}. The
twining number of a closed curve $\Gamma$ with respect to a point $p_{0}$
in ${\rm B}$ is denoted by $v(\Gamma,p_{0})$ and defined by
\begin{equation}
\label{twining1}
v(\Gamma,p_{0}) = \frac{1}{2\pi j} \odot \int_{\Gamma} \frac{dp}{p-p_{0}} \mbox{\femsteg , $p-p_{0}$ nonsingular for all $p\in \Gamma$}
\end{equation}
The integral has a finite value only if $p-p_{0}$ is nonsingular
at every point $p$ on $\Gamma$. The value cannot always be found by
application of Cauchy's theorem, and may therefore differ from zero.
The application of the theorem demands that $\Gamma$ is the boundary of some orientable
surface $S(\Gamma)$ so that $p-p_{0}$ is nonsingular for all $p\in S(\Gamma)$, 
but such a surface does not exist, if $\Gamma$ can
be regarded as enclosing $p_{0}$. The notion of curve enclosure
is more difficult to grasp intuitively in the bicomplex space than in
the complex plane, and the twining number is in fact the right device for
making it precise:

\begin{definition}
{\rm
The curve $\Gamma$ is said to {\em enclose} the point $p_{0}$
in the bicomplex space if $v(\Gamma,p_{0})\neq 0$. $\Box$
}
\end{definition} 

\smallskip

\noindent {\bf Example.} We shall compute the twining number
$v(\Gamma,p_{0})$ when $\Gamma$ is given by
\[
\Gamma:\mbox{\tresteg} p(t)=p_{0}+(e^{it}\cdot \cos t,e^{it}\cdot \sin t) \mbox{\femsteg , $0\leq t\leq 2\pi$}
\]
Notice as an aid to the geometric interpretation  of this curve
that all its points are at a distance of 1 from $p_{0}$.

We first observe that 
$p(t)-p_{0}=(e^{it}\cdot \cos t,e^{it}\cdot \sin t)=e^{(it,t)}$.
Because the exponential function gets only nonsingular values,
$p(t)-p_{0}$ is nonsingular for all $t$. We also have:
\begin{eqnarray*}
&&\frac{1}{p(t)-p_{0}} \;=\; \left(\frac{\cos t}{e^{it}}\:,\:-\frac{\sin t}{e^{it}}\right) \\ [3mm]
&&p\,'(t) \;=\; (i\cdot e^{it}\cdot \cos t - e^{it}\cdot \sin t\,,\,i\cdot e^{it}\cdot \sin t+e^{it}\cdot \cos t)
\end{eqnarray*}
This gives:
\begin{eqnarray*}
& & \int_{\Gamma} \frac{dp}{p-p_{0}} \\[2mm]
&=& \mbox{\femsteg $\lbrace$line integral, $0\leq t\leq 2\pi$$\rbrace$} \\[2mm]
& & \int_{0}^{2\pi} \frac{1}{p(t)-p_{0}}\odot p\,'(t)\, dt \\[2mm]
&=& \mbox{\femsteg $\lbrace$formulas above and bicomplex algebra$\rbrace$} \\[2mm]
& & \int_{0}^{2\pi} (i,1)\, dt \\[2mm]
&=& \mbox{\femsteg $\lbrace$integration$\rbrace$} \\[2mm]
& & (i2\pi,2\pi)
\end{eqnarray*}
Thus, the twining number (\ref{twining1}) becomes
$v(\Gamma,p_{0})=(1,-i)$. It differs from 0, which
means that $\Gamma$ does not bound any surface $S(\Gamma)$ that
would have the property that $1/(p-p_{0})$ is holomorphic
at all points of $S(\Gamma)$. For instance, take the surface
\[
S(\Gamma):\mbox{\tresteg} p(h,t)=p_{0}+(h\cdot e^{it}\cdot \cos t\,,\,h\cdot e^{it}\cdot \sin t) 
\]
where the real parameters $h$ and $t$ are in the intervals
$[0,1]$ and $[0,2\pi]$, respectively. It is bounded by $\Gamma$, but
contains the point $p(0,t)=p_{0}$, where
$1/(p-p_{0})$ is not holomorphic. 

\noindent $\Box$

\medskip

\noindent In order to find a general expression for the twining number
we evaluate the right-hand side of (\ref{twining1}) assuming that
$\Gamma$'s parametric equation is (\ref{Gamma1}) and that $p-p_{0}$
is nonsingular at all points $p\in \Gamma$.
We introduce the function
\begin{equation}
\label{Fdef}
F(t) = \int_{r}^{t} \frac{p\,'(t)}{p(t)-p_{0}} dt
\end{equation}
$F(s)$ is the value of the integral in (\ref{twining1}).
$F$ is continuous on the interval $[r,s]$ and except at those
points where $p\,'(t)$ is discontinuous it has the derivative
\[
F'(t) = \frac{p\,'(t)}{p(t)-p_{0}}
\]
As a result the derivative of the expression
$e^{-F(t)} \odot (p(t)-p_{0})$ vanishes except
at a finite number of points. The expression itself therefore
reduces to a constant, namely $p(r)-p_{0}$. Accordingly, we obtain
\[
e^{F(t)} = \frac{p(t)-p_{0}}{p(r)-p_{0}}
\]
For $t=s$ we have $e^{F(s)}=1$, since $\Gamma$ is closed and $p(s)=p(r)$.
To fulfill this identity $F(s)$ must due to (3.\ref{bicexpone}) be the period
of the exponential function. Consequently, setting $t=s$ in (\ref{Fdef})
yields first
\[
\int_{r}^{s} \frac{p\,'(t)}{p(t)-p_{0}} dt = (im,n)2\pi \mbox{\femsteg , $m$ and $n$ are integers}
\]
and then, on omission of the parametrization
\[
\int_{\Gamma} \frac{dp}{p-p_{0}} = (im,n)2\pi \mbox{\femsteg , $m$ and $n$ are integers}
\]
Hence, (\ref{twining1}) becomes
\begin{equation}
\label{twining2}
v(\Gamma,p_{0}) = (m,-in) \;\mbox{\femsteg , $m$ and $n$ are integers}
\end{equation}
For $n=0$ the twining number reduces to the complex winding number.

\subsection{Cauchy's integral formula in bicomplex space}
We are now ready to tackle the bicomplex generalization of Cauchy's integral formula.
\begin{theorem}
Let $\psi$ be a function, $p_{0}$ a point, $\Gamma$ a closed
curve and $G$ a domain in the bicomplex space
and assume that they satisfy the following conditions:
\newcounter{itnum}
\begin{list}
{\alph{itnum})}{\usecounter{itnum}
\setlength{\rightmargin}{\leftmargin}}
\item $\psi$ is holomorphic in $G$.
\item $\Gamma$ is contained in $G$ and 
the difference $p-p_{0}$ is nonsingular for all $p\in \Gamma$.
\item There exists an orientable surface $S(\Gamma)$
contained in $G$ and bounded by the curve $\Gamma$ so that $p_{0}\not\in S(\Gamma) \:\wedge\: v(\Gamma,p_{0})=0$ holds
and $p-p_{0}$ is nonsingular for all $p\in S(\Gamma)$, or
$p_{0}\in S(\Gamma) \:\wedge\: v(\Gamma,p_{0})\neq 0$ holds and
$p-p_{0}$ is nonsingular for all $p\in S(\Gamma)-\{p_{0}\}$.
\end{list}
Then, we have
\begin{equation} 
\label{cauchyintform}
\psi(p_{0})\odot v(\Gamma,p_{0}) = \frac{1}{2\pi j} \odot \int_{\Gamma} \frac{\psi(p)}{p-p_{0}}\odot dp 
\end{equation}
\end{theorem}

\noindent {\em Proof}. Our argumentation is essentially the same as the one given
by Nevanlinna and Paatero for the complex case~\cite{nevanlinna}, pp. 133--135. It is primarily based on
conditions a) and c), condition b) being a necessary condition of c).

Condition c) gives us two options with regard to the surface $S(\Gamma)$.
If it does not pass through $p_{0}$, the function $\psi(p)/(p-p_{0})$ is
holomorphic at all $p\in S(\Gamma)$. Formula (\ref{cauchyintform}) then holds,
since $v(\Gamma,p_{0})=0$, and the integral at the right-hand side is equal
to 0 due to Cauchy's theorem.

If $S(\Gamma)$ passes through $p_{0}$ and $v(\Gamma,p_{0})\neq 0$, i.e.
$\Gamma$ encloses $p_{0}$, we consider the function
\[
F(p) = \left\{\begin{array}{ll}
                    {\displaystyle \frac{\psi(p)-\psi(p_{0})}{p-p_{0}}} & \mbox{\femsteg if $p\neq p_{0}$} \\ [6mm]
                    \psi\,'(p_{0}) & \mbox{\femsteg if $p=p_{0}$}
             \end{array} \right.
\]
By virtue of conditions a) and c) $F$ is holomorphic
at all points of the punctured surface $S(\Gamma)-\{p_{0}\}$,
while the situation at $p_{0}$ remains open. Our aim is to
apply Cauchy's theorem to $F$ and for this purpose we 
must show that $F$ is holomorphic at $p_{0}$, too.

The function $\psi$ is holomorphic at $p_{0}$, hence
it has a Taylor expansion of the form (\ref{taylorseries}) about this point:
\begin{eqnarray*}
\psi(p) &=& \psi(p_{0})+\psi\,'(p_{0})\odot(p-p_{0})+\frac{1}{2}\psi\,''(p_{0})\odot(p-p_{0})^{2}+ \\[2mm]
&&\mbox{\femsteg \femsteg \femsteg\femsteg} \frac{1}{6}\psi\,'''(p_{0})\odot(p-p_{0})^{3}+R_{4}(p,p_{0})
\end{eqnarray*}
We get
\begin{eqnarray*}
F(p)\;\;=\;\; \frac{\psi(p)-\psi(p_{0})}{p-p_{0}} &=& \psi\,'(p_{0})+\frac{1}{2}\psi\,''(p_{0})\odot(p-p_{0})+ \\[2mm]
&&\mbox{\femsteg} \frac{1}{6}\psi\,'''(p_{0})\odot(p-p_{0})^{2}+\frac{R_{4}(p,p_{0})}{p-p_{0}}
\end{eqnarray*}
and further, because $F(p_{0})=\psi\,'(p_{0})$:
\[
\frac{F(p)-F(p_{0})}{p-p_{0}} \;=\; \frac{1}{2}\psi\,''(p_{0})+\frac{1}{6}\psi\,'''(p_{0})\odot(p-p_{0})+\frac{R_{4}(p,p_{0})}{(p-p_{0})^{2}}
\]
The inequality (\ref{remestimate}) implies that $R_{4}(p,p_{0})/(p-p_{0})^{2}$
tends to 0 as $p\rightarrow p_{0}$. We conclude that $F\,'(p_{0})=\frac{1}{2}\psi\,''(p_{0})$,
or that $F$ is holomorphic at $p_{0}$.

Since $F$ is holomorphic on the whole surface $S(\Gamma)$ it is legitimate
to use Cauchy's theorem
\[
\int_{\Gamma} F(p)\odot dp = 0
\]
which yields
\[
0 \;\;=\;\; \int_{\Gamma} \frac{\psi(p)-\psi(p_{0})}{p-p_{0}} \odot dp \;\;=\;\; \int_{\Gamma} \frac{\psi(p)}{p-p_{0}} \odot dp \;-\; \psi(p_{0})\odot \int_{\Gamma} \frac{dp}{p-p_{0}}
\]
Formula (\ref{cauchyintform}) follows if we apply (\ref{twining1})
and substitute $v(\Gamma,p_{0})\odot 2\pi j$ for
the last integral. \\
$\Box$

\bigskip

\noindent Assuming that $v(\Gamma,p_{0})\neq 0$ in (\ref{cauchyintform}) and that
conditions a) -- c) of the preceding theorem otherwise
also hold, we may select a particular curve $\Gamma$ so that $v(\Gamma,p_{0})=(0,-i)$
in agreement with (\ref{twining2}) to obtain
\begin{equation} 
\label{cauchyintform2}
\psi(p_{0}) = \frac{1}{2\pi i} \odot \int_{\Gamma} \frac{\psi(p)}{p-p_{0}}\odot dp 
\end{equation}
We conclude that for such a curve the common form of Cauchy's integral formula
is valid in ${\rm B}$, too.

\bigskip

\noindent {\bf Example.} Given the closed curve
\[
\Gamma:\mbox{\tresteg} p(t)=(e^{it},e^{it}) \mbox{\femsteg , $0\leq t\leq 2\pi$}
\]
we shall compute the value
\[
W = \frac{1}{2\pi i} \odot \int_{\Gamma} \frac{e^{p}}{p}\odot dp 
\]
By making the identifications $\psi(p)=e^{p}$, $p_{0}=0$ and $G={\rm B}$ we see
that the conditions a)--c) of Theorem 5.6 are met. Especially concerning c)
we note that $v(\Gamma,0)=(0,-i)$  and that all points $p\in \Gamma$ are
nonsingular (since $e^{it}\neq ie^{it} \:\wedge\: e^{it}\neq -ie^{it}$ for all $t\in [0,2\pi]$).
Moreover, the surface
\[
S(\Gamma):\mbox{\tresteg} p(h,t)=(h\cdot e^{it}\,,\,h\cdot e^{it}) \mbox{\femsteg , $0\leq t\leq 2\pi$ and $0\leq h \leq 1$}
\]
passes through the point $p(0,t)=0$ and all points
$p\in S(\Gamma)-\{0\}$ are nonsingular. Therefore, formula
(\ref{cauchyintform2}) is applicable and yields $W=e^{0}=1$.
We arrive at this result by direct calculation, too:
\begin{eqnarray*}
& & \int_{\Gamma} \frac{e^{p}}{p}\odot dp \\[2mm]
&=& \mbox{\femsteg $\lbrace$bicomplex line integral, $0\leq t\leq 2\pi$$\rbrace$} \\[2mm]
& & \int_{0}^{2\pi} \frac{e^{p(t)}}{p(t)}\odot p\,'(t)\, dt \\[2mm]
&=& \mbox{\femsteg $\lbrace$definitions of $e^{p}$ and $\Gamma$; bicomplex inverse$\rbrace$} \\[2mm]
& & \int_{0}^{2\pi} (e^{e^{it}}\!\cdot \cos(e^{it})\:,\: e^{e^{it}}\!\cdot \sin(e^{it}))\odot \left(\frac{e^{it}}{2e^{i2t}}\,,\,\frac{-e^{it}}{2e^{i2t}}\right)\odot(ie^{it},ie^{it})\,dt \\[2mm]
&=&  \mbox{\femsteg $\lbrace$multiplication, distribution of the integral$\rbrace$} \\[2mm]
& & \left(\int_{0}^{2\pi} \frac{e^{e^{it}}\!\cdot \cos(e^{it})}{e^{it}}\cdot ie^{it}\,dt\;,\;\int_{0}^{2\pi} \frac{e^{e^{it}}\!\cdot \sin(e^{it})}{e^{it}}\cdot ie^{it}\,dt \right) \\[2mm]
&=&  \mbox{\femsteg $\lbrace$$\gamma_{1}\!: a(t)=e^{it}$ and $\gamma_{2}\!: b(t)=e^{it}$ are the components of $\Gamma$$\rbrace$} \\ [2mm]
& & \left(\int_{0}^{2\pi} \frac{e^{a(t)}\!\cdot \cos a(t)}{a(t)}\cdot a\,'(t)\,dt\;,\;\int_{0}^{2\pi} \frac{e^{b(t)}\!\cdot \sin b(t)}{b(t)}\cdot b\,'(t)\,dt \right) \\[2mm]
&=&  \mbox{\femsteg $\lbrace$definition of line integral in ${\rm C}$$\rbrace$} \\ [2mm]
& & \left(\int_{\gamma_{1}}\frac{e^{a}\!\cdot \cos a}{a} \cdot \,da\;,\;\int_{\gamma_{2}} \frac{e^{b}\!\cdot \sin b}{b}\cdot \,db \right) \\[2mm]
&=&  \mbox{\femsteg $\lbrace$Cauchy's integral formula in ${\rm C}$$\rbrace$} \\ [2mm]
& & (2\pi i,0)
\end{eqnarray*}
The result is $W=1$, as it should. \\
$\Box$

\newpage
\setcounter{equation}{0}

\section{Bicomplex harmonic analysis}
Our last topic is the relationship between
the bicomplex functions and the Laplace equation. It is well-known
that the complex holomorphic functions satisfy the two-dimensional
Laplace equation, hence one may expect that the bicomplex holomorphic
functions satisfy its four-dimensional version.
We shall show that this is indeed true for eight classes
of bicomplex holomorphic functions. One of them consists of the
functions that fulfill the bicomplex Cauchy-Riemann equations
(4.\ref{BCR1})--(4.\ref{BCR2}), the others of similar
functions that are partly or entirely formed with the {\em conjugates}
of complex holomorphic functions. These functions quite naturally
emerge from the factorization of the Laplace equation in its various forms. We start
by factorizing the two-dimensional equation, since it will provide
us with important concepts for the more general cases.

\subsection{The Cauchy-Riemann operator}
With $f$ a complex function depending on the
real variables $x$ and $y$, the two-dimensional
Laplace equation is written
\[
\frac{\partial^{2}f}{\partial x^{2}}+\frac{\partial^{2}f}{\partial y^{2}} = 0
\]
or, using abbreviated differential operators
\begin{equation}
\label{Lapleqr2}
\partial_{x}^{2}f+\partial_{y}^{2}f = 0
\end{equation}
The two-dimensional
Laplace operator $\bigtriangleup_{2}= \partial_{x}^{2}+\partial_{y}^{2}$
admits the factorization
\begin{equation}
\label{Laplopr2}
\bigtriangleup_{2}=(\partial_{x}+i\partial_{y})\cdot (\partial_{x}-i\partial_{y})
\end{equation}
\noindent which for the complex number $a=x+iy$ gives us the
occasion to introduce the so-called {\em Cauchy-Riemann operator}
$D_{a}$~\cite{brackx}:
\begin{equation}
\label{CRop}
D_{a}=\partial_{x}+i\partial_{y} \mbox{\femsteg , $a=x+iy$}
\end{equation}
Note that $D_{a}^{*}=\partial_{x}-i\partial_{y}$.

We now let the complex function $f$ in (1) more specifically be given by
\begin{equation}
\label{fdefc}
f(a)= g(x,y)+ih(x,y) \mbox{\femsteg , $a=x+iy$}
\end{equation}
where $g$ and $h$ are functions of the type ${\rm R}^{2} \rightarrow {\rm R}$
such that they have continuous partial derivatives in $x$ and $y$. 
From (2) and (3) we obtain the relation
\begin{equation}
\label{Laplsolr2}
\bigtriangleup_{2}f(a)=0 \:\Leftarrow\: (D_{a}\cdot f(a)=0) \vee (D_{a}^{*}\cdot f(a)=0)
\end{equation}
by means of which the Laplace equation of the consequent
can be solved by solving either disjunct of the antecedent. We first evaluate
$D_{a}\cdot f(a)=0$:
\begin{tabbing}
xxxxxxx\=xxx\=xxxxxxxxxxxxxxxxxxxxxxxxxxxxxxxxxxxxxxxxxxxxxxxxxxxxxxxxxxxxxxxxx\kill
 \> \>$D_{a}\cdot f(a)=0$\\
 \> $\equiv$ \> \tresteg $\lbrace$definition of $D_{a}$ and $f$, omitting arguments$\rbrace$\\
 \> \> $(\partial_{x}+i\partial_{y})\cdot (g+ih)=0$\\
 \> $\equiv$ \> \tresteg $\lbrace$complex multiplication$\rbrace$\\
 \> \> $\partial_{x}g-\partial_{y}h+i(\partial_{x}h+\partial_{y}g)=0$\\
 \> $\equiv$ \> \tresteg $\lbrace$complex algebra$\rbrace$\\
 \> \> $(\partial_{x}g-\partial_{y}h=0) \wedge (\partial_{x}h+\partial_{y}g=0)$
\end{tabbing}

\medskip

\noindent The first disjunct of the antecedent of (\ref{Laplsolr2})
obviously resolves into the Cauchy-Riemann equations
\begin{equation}
\label{CReqsr2}
\partial_{x}g=\partial_{y}h \mbox{\tresteg , \tresteg} \partial_{x}h=-\partial_{y}g
\end{equation}
An analogous treatment
of the second disjunct $D_{a}^{*}\cdot f(a)=0$ shows
that it is equivalent to
\begin{equation}
\label{ConjCReqsr2}
\partial_{x}g=-\partial_{y}h \mbox{\tresteg , \tresteg} \partial_{x}h=\partial_{y}g
\end{equation}
We call these the {\em conjugate Cauchy-Riemann equations}.
The name is justified by the fact that they are solved
by the conjugates of the holomorphic complex functions.
These are henceforth referred to as the {\em conjugate
holomorphic functions}. 

\medskip

\noindent {\bf Example.} The exponential function
$g+ih=e^{x}\cos y + ie^{x}\sin y$ is holomorphic and
satisfies (\ref{CReqsr2}). Its conjugate $(g+ih)^{*}=e^{x}\cos y + i(-e^{x}\sin y)$
satisfies (\ref{ConjCReqsr2}). $\Box$

\medskip

\noindent From the above we deduce the following relationship
between the Cauchy-Riemann operator and the holomorphic/conjugate holomorphic
complex functions:

\begin{eqnarray}
\label{CRopf1}
D_{a}\cdot f(a)=0 &\equiv& \mbox{$f$ is a holomorphic function} \\ [2mm]
\label{CRopf2}
D_{a}^{*}\cdot f(a)=0 &\equiv& \mbox{$f$ is a conjugate holomorphic function}
\end{eqnarray}
The application of the Cauchy-Riemann operator or its conjugate
to a complex function thus reveals whether the function is
differentiable or not.

\subsection{The derivatives of a complex function}
It is a noteworthy fact that given a complex number $d$ with non-zero real and imaginary parts
{\em any} complex number $a$ can be expressed in the form
$a=\alpha d+\beta d^{*}$, where $\alpha$ and $\beta$ are real numbers.
Because $d$ and its conjugate $d^{*}$ span the entire complex
plane in this way, one could argue that in the theory of
complex functions complex-valued variables and functions
ought to be treated on an equal footing with their conjugates.
This should especially apply to the derivatives of the functions.
With the complex function $f$ specified by (\ref{fdefc}) it is possible to distinguish
between four derivatives: 
\[
\frac{df}{da}\;,\; \frac{df^{*}}{da^{*}}\;,\;\frac{df}{da^{*}}\;,\;\frac{df^{*}}{da}
\]
We first recall the 
definition of the first of these derivatives at the point $a_{0}=x_{0}+iy_{0}$
\begin{equation}
\label{derdeffa0}
\frac{df}{da}
             \begin{array}{l}
                  \\
               {\scriptstyle \mid a=a_{0}}
             \end{array}
                              = \lim_{a \rightarrow a_{0}} \frac{f(a)-f(a_{0})}{a-a_{0}}
\end{equation}

\noindent Letting $a$ approach $a_{0}$ along the real and imaginary axis
separately leads to the formulas (the evaluation at the
point $a_{0}$ now omitted)
\begin{eqnarray}
\label{derdeffa1}
\frac{df}{da} &=& \partial_{x}g+i\partial_{x}h \\ [3mm]
\label{derdeffa2}
\frac{df}{da} &=& \partial_{y}h-i\partial_{y}g
\end{eqnarray}
\noindent which are equivalent if the CR-equations (\ref{CReqsr2}) hold.

Next we consider the derivative of the conjugate function
$f^{*}$ with respect to the conjugate variable $a^{*}$. We define
\[
\frac{df^{*}}{da^{*}}
             \begin{array}{l}
                  \\
               {\scriptstyle \mid a^{*}=a_{0}^{*}}
             \end{array}
                              = \lim_{a^{*} \rightarrow a_{0}^{*}} \frac{f^{*}(a)-f^{*}(a_{0})}{a^{*}-a_{0}^{*}}
\]

\noindent But if the limit at the right-hand side of (\ref{derdeffa0}) exists we have
\[
\lim_{a^{*} \rightarrow a_{0}^{*}} \frac{f^{*}(a)-f^{*}(a_{0})}{a^{*}-a_{0}^{*}}=\left(\lim_{a \rightarrow a_{0}} \frac{f(a)-f(a_{0})}{a-a_{0}}\right)^{*}
\]

\noindent from which it follows that
\begin{equation}
\label{derdefconjfconja0}
\frac{df^{*}}{da^{*}}=\left(\frac{df}{da}\right)^{*}
\end{equation}
Conjugation of (\ref{derdeffa1}) and (\ref{derdeffa2}) consequently results in two formulas for
$\frac{df^{*}}{da^{*}} $:
\begin{eqnarray}
\label{derdefconjfconja1}
\frac{df^{*}}{da^{*}} &=& \partial_{x}g-i\partial_{x}h \\ [3mm]
\label{derdefconjfconja2}
\frac{df^{*}}{da^{*}} &=& \partial_{y}h+i\partial_{y}g
\end{eqnarray}
\noindent The equivalence of these equations is also implied by the CR-equations.

Then we turn our attention to the derivative of
$f$ with respect to $a^{*}$. Its definition is
\[
\frac{df}{da}
             \begin{array}{l}
                  \\
               {\scriptstyle \mid a^{*}=a_{0}^{*}}
             \end{array}
                              = \lim_{a^{*} \rightarrow a_{0}^{*}} \frac{f(a)-f(a_{0})}{a^{*}-a_{0}^{*}}
\]
The right-hand side can be rewritten
\[
\lim_{\begin{array}{l}
        {\scriptstyle x-iy \rightarrow} \\
        {\scriptstyle x_{0}-iy_{0} }
      \end{array}} \frac{g(x,y)-g(x_{0},y_{0})+i(h(x,y)-h(x_{0},y_{0}))}{x-iy-(x_{0}-iy_{0})}
\]
Letting $a^{*}$ approach $a_{0}^{*}$ along the
real and imaginary axis amounts to setting $y=y_{0}$
and $x=x_{0}$, respectively, in this formula. Taking the
two limits separately yields
\begin{eqnarray}
\label{derdeffconja1}
\frac{df}{da^{*}} &=& \partial_{x}g+i\partial_{x}h \\ [3mm]
\label{derdeffconja2}
\frac{df}{da^{*}} &=& -\partial_{y}h+i\partial_{y}g
\end{eqnarray}
These equations are equivalent if the conjugate CR-equations
(\ref{ConjCReqsr2}) hold, i.e. if $f$ is conjugate holomorphic.

For the derivative $\frac{df^{*}}{da}$, finally, we deduce with the same
technique that lead to (\ref{derdefconjfconja0})
\begin{equation}
\label{derdefconjfa0}
\frac{df^{*}}{da}=\left(\frac{df}{da^{*}}\right)^{*}
\end{equation}
Thus, by conjugating (\ref{derdeffconja1}) and (\ref{derdeffconja2}) we get
\begin{eqnarray}
\label{derdefconjfa1}
\frac{df^{*}}{da} &=& \partial_{x}g-i\partial_{x}h \\ [3mm]
\label{derdefconjfa2}
\frac{df^{*}}{da} &=& -\partial_{y}h-i\partial_{y}g
\end{eqnarray}
whose equivalence likewise follows from the conjugate
CR-equations.

Summarizing these results, we have:
\begin{itemize}
\item Formulas (\ref{derdeffa1})/(\ref{derdeffa2}) and (\ref{derdefconjfconja1})/(\ref{derdefconjfconja2}) are pairwise equivalent
      if $f$ is a {\em holomorphic} complex function.
\item Formulas (\ref{derdeffconja1})/(\ref{derdeffconja2}) and (\ref{derdefconjfa1})/(\ref{derdefconjfa2}) are pairwise equivalent
      if $f$ is a {\em conjugate holomorphic} complex function.
\end{itemize}

\noindent The preceding argument admittedly contains some redundancy,
because we have obtained four pairs of equivalent formulas, although
we have only two types of functions, holomorphic and conjugate holomorphic ones.
The reason is that we have designated a conjugate holomorphic function
by both $f$ and $f^{*}$. The duplicity could of course easily be removed,
but experience has shown that one needs both notations and therefore all the
formulas, too.

It is clear that in the complex plane the conjugate holomorphic functions must possess
properties that are very similar to those of the holomorphic functions.
The conjugate CR-equations (\ref{ConjCReqsr2}) play the same role for the
former functions as the CR-equations play for the latter;
in practice this means that if a relation $R(x,y)$ has been
shown to hold for an holomorphic function $f$ of the form (\ref{fdefc})
the relation $R(x,-y)$ will hold for $f^{*}$. Due to this kind
of close similarity it can be regarded as an effort devoid of interest
to actually develop a theory of conjugate holomorphic functions.
What we shall mainly need from it in the sequel is that
a conjugate holomorphic function $f$ of the form (\ref{fdefc}) has
a continuous derivative $\frac{df}{da^{*}}$ computable
with (\ref{derdeffconja1}) or (\ref{derdeffconja2}). Alternatively, if the function is designated
by $f^{*}$ we can apply (\ref{derdefconjfconja1}) or (\ref{derdefconjfconja2}) to find the derivative.

\subsection{Factorization of the four-dimensional Laplace equation}
We return to the bicomplex space and consider once more a function $\psi$ of the type
\begin{eqnarray}
\label{psidef3}
\psi(p) &=& (\phi_{1}(a,b),\phi_{2}(a,b)) \\ [2mm]
(\phi_{1}(a,b),\phi_{2}(a,b)) &=& (\psi_{1}(x,y,z,u)+i\psi_{2}(x,y,z,u), \nonumber \\
& &\; \psi_{3}(x,y,z,u)+i\psi_{4}(x,y,z,u)) \nonumber \\ [2mm]
p=(a,b)&,&\mbox{$a=x+iy$ and $b=z+iu$} \nonumber
\end{eqnarray}
We want $\psi$ to fulfill the
four-dimensional Laplace equation
\begin{equation}
\label{lapleqr4def}
\partial_{x}^{2}\psi + \partial_{y}^{2}\psi + \partial_{z}^{2}\psi + \partial_{u}^{2}\psi = 0
\end{equation}
which we shall solve by factorization. Using the 
Laplace-operator 
\[
\bigtriangleup_{4} = \partial_{x}^{2} + \partial_{y}^{2} + \partial_{z}^{2} + \partial_{u}^{2}
\]
we render (\ref{lapleqr4def}) in the form
\begin{equation}
\label{lapleqr4op}
\bigtriangleup_{4} \times \psi =0
\end{equation}
where $\times$ denotes quaternionic multiplication (2.\ref{quatmultc2}).
The $\bigtriangleup_{4}$-operator is to
be regarded as a scalar quaternion, hence the multiplication in~
(\ref{lapleqr4op}) merely applies it to each component $\psi_{k}$
of $\psi$ --- see (2.\ref{realxquat}).

We now make the crucial observation that the Laplace operator admits the factorization
\begin{equation}
\label{laplopr4fact}
\bigtriangleup_{4}\; = \; (\partial_{x}+ i\partial_{y},\partial_{z}+ i\partial_{u})\times (\partial_{x}- i\partial_{y},-\partial_{z}- i\partial_{u})
\end{equation}
\noindent If we introduce the new operator $L_{q}$ by
\begin{equation}
\label{lopdef}
L_{q}=(\partial_{x}+ i\partial_{y},\partial_{z}+ i\partial_{u})
\end{equation}
the factorization becomes
\smallskip
\begin{equation}
\label{laplophfact}
\bigtriangleup_{4} = L_{q}\times L_{q}^{*}
\end{equation}
The Cauchy-Riemann operators $D_{a}=\partial_{x}+i\partial_{y}$
and $D_{b}=\partial_{z}+i\partial_{u}$ give us
\begin{equation}
\label{lopc2fact}
L_{q}=(D_{a},D_{b}) \mbox{\tresteg , \tresteg} L_{q}^{*}=(D_{a}^{*},-D_{b})
\end{equation}
and
\smallskip
\begin{equation}
\label{laplopc2fact}
\bigtriangleup_{4}\: = \:(D_{a},D_{b})\times (D_{a}^{*},-D_{b})
\end{equation}
With formulas (\ref{laplopr4fact}), (\ref{laplophfact}) and (\ref{laplopc2fact}) we have factorized the four-dimensional
Laplace operator in ${\rm R}^{4}$, ${\rm B}$, and ${\rm C}^{2}$, respectively.

Because $L_{q}\times L_{q}^{*}= L_{q}^{*}\times L_{q}$ we deduce from (\ref{lapleqr4op})
and (\ref{laplophfact})
\begin{equation}
\label{lapleqhfact}
\bigtriangleup_{4} \times \psi =0 \;\Leftarrow\; (L_{q}\times \psi =0) \vee (L_{q}^{*}\times \psi=0)
\end{equation}
by means of which we can solve the Laplace equation by solving
either disjunct of the antecedent:
\begin{eqnarray}
\label{crfh} 
L_{q}\times \psi&=&0 \\ [1mm] 
\label{ccrfh}
L_{q}^{*}\times \psi&=&0 
\end{eqnarray}
We focus on the first of these. With (\ref{lopc2fact}) we first rewrite it in ${\rm C}^{2}$:
\begin{eqnarray*}
& &L_{q}\times \psi=0 \\
&\equiv& \mbox{ \tresteg $\lbrace L_{q}=(D_{a},D_{b})$, $\psi=(\phi_{1},\phi_{2})\rbrace$} \\
& &(D_{a},D_{b})\times (\phi_{1},\phi_{2})=0 \\
&\equiv& \mbox{\tresteg $\lbrace$quaternionic multiplication$\rbrace$} \\
& &(D_{a}\cdot \phi_{1}-D_{b}\cdot \phi_{2}^{*}\: ,\: D_{b}\cdot \phi_{1}^{*}+ D_{a}\cdot \phi_{2})=0 \\
&\equiv& \mbox{\tresteg $\lbrace$equality of quaternions$\rbrace$} \\
& &(D_{a}\cdot \phi_{1}-D_{b}\cdot \phi_{2}^{*}=0) \wedge (D_{b}\cdot \phi_{1}^{*}+D_{a}\cdot \phi_{2}=0) \\
\end{eqnarray*}
Thus, (\ref{crfh}) is equivalent to the simultaneous validity of the
equations
\begin{eqnarray}
\label{crfc21}
D_{a}\cdot \phi_{1}&=&D_{b}\cdot \phi_{2}^{*} \\
\label{crfc22} 
D_{b}\cdot \phi_{1}^{*}&=&-D_{a}\cdot \phi_{2} 
\end{eqnarray}
Substituting $D_{a}=\partial_{x}+i\partial_{y}$,
$D_{b}=\partial_{z}+i\partial_{u}$ and $\phi_{1}=\psi_{1}+i\psi_{2}$,
$\phi_{2}=\psi_{3}+i\psi_{4}$ in these we obtain with
complex algebra the corresponding ${\rm R}^{4}$-representation:
\begin{eqnarray}
\label{crfr41}
\partial_{x}\psi_{1} - \partial_{y}\psi_{2} - \partial_{z}\psi_{3} - \partial_{u}\psi_{4} &=& 0 \\ 
\label{crfr42}
\partial_{x}\psi_{2} + \partial_{y}\psi_{1} + \partial_{z}\psi_{4} - \partial_{u}\psi_{3} &=& 0 \\ 
\label{crfr43}
\partial_{x}\psi_{3} - \partial_{y}\psi_{4} + \partial_{z}\psi_{1} + \partial_{u}\psi_{2} &=& 0 \\ 
\label{crfr44}
\partial_{x}\psi_{4} + \partial_{y}\psi_{3} - \partial_{z}\psi_{2} + \partial_{u}\psi_{1} &=& 0 
\end{eqnarray}
This partial differential equation was derived in 1935 by Fueter~\cite{fueter}.
As mentioned in the Introduction his starting point
was the formula
\begin{equation}
\label{crfvect}
\partial_{x}\psi + i\times\partial_{y}\psi + j\times\partial_{z}\psi + k\times\partial_{u}\psi = 0
\end{equation}
which we get by substituting the vector representation of $L_{q}$, or
$L_{q}=\partial_{x} + \partial_{y}i + \partial_{z}j + \partial_{u}k$, in $L_{q}\times \psi = 0$.

Equations (\ref{crfh}), (\ref{crfc21})--(\ref{crfc22}), (\ref{crfr41})--(\ref{crfr44}), and (\ref{crfvect}) are four equivalent ways
of expressing the same condition. We shall refer to them as the
{\em Fueter equations}.

The treatment of equation (\ref{ccrfh}), $L_{q}^{*}\times \psi=0$,
which also enables one to solve the Laplace equation, is analogous with
the first case. Written in the form
$(D_{a}^{*},-D_{b})\times (\phi_{1},\phi_{2})=0$ --- see
(\ref{psidef3}) and (\ref{lopc2fact}) --- it becomes equivalent to the simultaneous
validity of
\begin{eqnarray}
\label{ccrfc21}
D_{a}^{*}\cdot \phi_{1}&=&-D_{b}\cdot \phi_{2}^{*} \\
\label{ccrfc22}
D_{b}\cdot \phi_{1}^{*}&=&D_{a}^{*}\cdot \phi_{2}
\end{eqnarray}
We call these the {\em conjugate Fueter equations}.
Expanded in ${\rm R^{4}}$ they look like:
\begin{eqnarray}
\label{ccrfr41}
\partial_{x}\psi_{1} + \partial_{y}\psi_{2} + \partial_{z}\psi_{3} + \partial_{u}\psi_{4} &=& 0 \\
\label{ccrfr42}
\partial_{x}\psi_{2} - \partial_{y}\psi_{1} - \partial_{z}\psi_{4} + \partial_{u}\psi_{3} &=& 0 \\
\label{ccrfr43}
\partial_{x}\psi_{3} + \partial_{y}\psi_{4} - \partial_{z}\psi_{1} - \partial_{u}\psi_{2} &=& 0 \\
\label{ccrfr44}
\partial_{x}\psi_{4} - \partial_{y}\psi_{3} + \partial_{z}\psi_{2} - \partial_{u}\psi_{1} &=& 0
\end{eqnarray}
This partial differential equation was derived in 1929 by Lanczos~\cite{lanczos}.
The fourth, equivalent formulation of the conjugate Fueter equations
is found by substitution of the
vector-form $L_{q}^{*}=\partial_{x} - \partial_{y}i - \partial_{z}j - \partial_{u}k$
in $L_{q}^{*}\times \psi = 0$, which becomes
\begin{equation}
\label{ccrfvect}
\partial_{x}\psi - i\times\partial_{y}\psi - j\times\partial_{z}\psi - k\times\partial_{u}\psi = 0
\end{equation}

\subsection{Regular and conjugate regular functions}
In the solutions of the Fueter equations we identify
a new type of bicomplex functions, hence we define:
\begin{definition}
{\rm
A bicomplex function 
is said to be {\em regular} in a domain of the bicomplex space
if it satisfies the Fueter equations in this domain. $\Box$
}
\end{definition}
To further characterize the regular functions we prove:
\begin{theorem}
Let the bicomplex function $\psi$ be defined in the domain $G\subseteq {\rm B}$ so that
$\psi(p)= (\phi_{1}(a,b),\phi_{2}(a,b))$ for $p=(a,b)$.
If the functions $\phi_{1}$ and $\phi_{2}$ are holomorphic
in $a$ and conjugate holomorphic in $b$ at all points $(a,b)\in G$,
then $\psi$ is regular in $G$.
\end{theorem}
{\em Proof}. The assumptions about $\phi_{1}$ and $\phi_{2}$
imply on account of (\ref{CRopf1}) and (\ref{CRopf2}) that
\begin{eqnarray*}
D_{a}\cdot \phi_{k}=0 & & \mbox{\tresteg , $\; k=1,2$}\\
D_{b}^{*}\cdot \phi_{k}=0 & & \mbox{\tresteg , $\; k=1,2$}
\end{eqnarray*}
hold in $G$. Conjugating the second formula yields $D_{b}\cdot \phi_{k}^{*}=0$,
which combined with the first validates the Fueter equations
(\ref{crfc21}) and (\ref{crfc22}). Hence, $\psi$ is regular in $G$.

\noindent $\Box$

\bigskip

\noindent An example of a regular function is
\begin{equation}
\label{bicEdef}
E(p) = (e^{a}\!\cdot \cos^{*}b, e^{a}\!\cdot \sin^{*}b) \mbox{\femsteg , $\; p=(a,b)$}
\end{equation}
The complex pair components $e^{a}\!\cdot \cos^{*}b$ and $e^{a}\!\cdot \sin^{*}b$
are functions of the
variables $a$ and $b$ such that they are holomorphic in $a$ and
conjugate holomorphic in $b$, as required.

Two-argument, complex functions that are holomorphic in the one argument and
conjugate holomorphic in the other 
are evidently of special interest and deserve an own term:
\begin{definition}
{\rm
A two-argument, complex function $\phi\!:(a,b) \rightarrow \phi(a,b)$,
defined in a domain $G\subseteq {\rm B}$, is said to be {\em coanalytic in $G$}
if it is holomorphic in $a$ and conjugate holomorphic in $b$ at all points $(a,b)\in G$. $\Box$
}
\end{definition}
Due to Theorem 6.1 the coanalyticity of the components of
a bicomplex function $\psi$ is a {\em sufficient} condition
for $\psi$'s regularity.

It is a fundamental property of any conjugate
holomorphic complex function $f^{*}$ that
$f^{*}(b)=f(b^{*})$. The right-hand side of
(\ref{bicEdef}) can therefore be rewritten 
$(e^{a}\!\cdot \cos b^{*}, e^{a}\!\cdot \sin b^{*})$,
which we also obtain by replacing the argument
$p=(a,b)$ of the exponential function (3.\ref{bicexpdef})
by $q=(a,b^{*})$:
\begin{equation}
\label{regexpdef}
e^{q}=(e^{a}\!\cdot \cos b^{*}, e^{a}\!\cdot \sin b^{*}) \mbox{\femsteg , $q=(a,b^{*})$}
\end{equation}
With the same change of argument we derive the regular counterparts
of the other elementary bicomplex functions of Chapter 3, e.g.
the regular identity function and power function --- see (3.\ref{bicidfunc}) and (3.\ref{bicpowfunc}):
\begin{eqnarray}
\label{regiddef}
I(q) &=& (a,b^{*})\;\; \mbox{\femsteg , $q=(a,b^{*})$} \\ [1mm]
\label{regpowfunc}
F(q)&=&(a,b^{*})^{n} \mbox{\femsteg , $q=(a,b^{*})$ and $n$ an integer $\geq 0$}
\end{eqnarray}
and the regular hyperbolic and trigonometric
functions --- see (3.\ref{biccoshc2})--(3.\ref{bicsinhc2})
and (3.\ref{biccosc2})--(3.\ref{bicsinc2}):
\begin{eqnarray}
\label{regcoshc2}
\cosh q&=&(\cosh a \cdot \cos b^{*} , \sinh a \cdot \sin b^{*}) \mbox{\femsteg , $q=(a,b^{*})$} \\
\label{regsinhc2}
\sinh q&=&(\sinh a \cdot \cos b^{*} , \cosh a \cdot \sin b^{*}) \\ [2mm]
\label{regcosc2}
\cos q &=& (\cos a \cdot \cosh b^{*}, -\sin a \cdot \sinh b^{*}) \\
\label{regsinc2}
\sin q &=& (\sin a \cdot \cosh b^{*} , \cos a \cdot \sinh b^{*})
\end{eqnarray}
Note that all the component functions of the complex pairs at the right-hand sides
are coanalytic in ${\rm B}$.

The holomorphism of
one-argument complex functions is preserved by complex
addition, multiplication and division as well as 
by functional composition. The conjugate holomorphism
of one-argument functions is similarly preserved. Being a
combination of these properties, the coanalyticity 
of two-argument complex functions is
thus preserved by the aforementioned operations. 

\bigskip

\noindent {\bf Example.} The functions $a\cdot \cos b^{*}$ and $b^{*}\cdot \sin a$ are
coanalytic in ${\rm C}^{2}$, hence $a\cdot \cos b^{*} + b^{*}\cdot \sin a$
and $(a\cdot \cos b^{*}) \cdot (b^{*}\cdot \sin a)$ are coanalytic in ${\rm C}^{2}$. $\Box$

\medskip

\noindent We are now able to prove that the regular functions form
an algebra under the bicomplex operations and that, in addition,
functional composition preserves the regularity property.
\begin{theorem}
Let the  bicomplex functions $\psi\!: q \rightarrow \psi(q)$ and 
$\theta\!: q \rightarrow \theta(q)$, where $q=(a,b^{*})$,
be regular in a domain $G\subseteq {\rm B}$, and
assume that they can be represented by complex pairs of coanalytic functions.
Then their bicomplex sum $\psi + \theta$,
product $\psi \odot \theta$ and quotient $\psi / \theta$
are regular in $G$, the last only at those points $q_{0}$ where $C\!N(\theta(q_{0}))\neq 0$.
Furthermore, the composite function 
$(\psi \circ \theta)(q) = \psi(\theta(q))$ is regular in $G$,
provided the domain of $\psi$ is in the range of $\theta$.
\end{theorem}
{\em Proof.} We verify the regularity of the product $\psi \odot \theta$, in particular.
According to the assumption $\psi$ and $\theta$ have
representations of the form
\begin{eqnarray*}
\psi(q) &=& (\phi_{1}(a,b^{*}),\phi_{2}(a,b^{*})) \mbox{\femsteg , $q=(a,b^{*})$} \\ [1mm]
\theta(q) &=& (\omega_{1}(a,b^{*}),\omega_{2}(a,b^{*}))
\end{eqnarray*}
where the functions $\phi_{1},\phi_{2},\omega_{1}$ and
$\omega_{2}$ are coanalytic in $G$. Bicomplex multiplication
gives
\[
\psi \odot \theta = (\phi_{1} \cdot \omega_{1} - \phi_{2} \cdot \omega_{2}\,,\, \phi_{2} \cdot \omega_{1}+ \phi_{1} \cdot \omega_{2})
\]
Because complex addition and multiplication preserve coanalyticity,
the component functions $\phi_{1} \cdot \omega_{1} - \phi_{2} \cdot \omega_{2}$
and $\phi_{2} \cdot \omega_{1}+ \phi_{1} \cdot \omega_{2}$
are coanalytic in $G$, thereby making $\psi \odot \theta$ regular.
The regularity of $\psi + \theta$, $\psi / \theta$ and 
$\psi \circ \theta$ is verified by analogous arguments.

\noindent $\Box$

\medskip

\noindent Next we focus on the conjugate Fueter equations, especially
in their form (\ref{ccrfc21})--(\ref{ccrfc22}). A bicomplex
function satisfying these equations is said to be {\em conjugate regular}.
We prove a sufficient condition for a function to be 
conjugate regular (compare Theorem 6.1).
\begin{theorem}
Let the bicomplex function $\psi$ be defined in the domain $G\subseteq {\rm B}$ so that
$\psi(p)= (\phi_{1}(a,b),\phi_{2}(a,b))$ for $p=(a,b)$.
If the functions $\phi_{1}$ and $\phi_{2}$ are conjugate holomorphic
in both $a$ and $b$ at all points $(a,b)\in G$ then $\psi$ is conjugate regular in $G$.
\end{theorem}

\noindent {\em Proof}. The assumption about $\phi_{1}$ and $\phi_{2}$ 
implies on account of (\ref{CRopf2}) that
\begin{eqnarray*}
D_{a}^{*}\cdot \phi_{k}=0 & & \mbox{\tresteg , $\; k=1,2$}\\
D_{b}^{*}\cdot \phi_{k}=0 & & \mbox{\tresteg , $\; k=1,2$}
\end{eqnarray*}
hold in $G$. The first formula and the conjugation
of the second then together validate (\ref{ccrfc21})--(\ref{ccrfc22}).

\noindent $\Box$

\medskip

\noindent If we take the quaternionic conjugate of the
regular identity function (\ref{regiddef}) we get
$I^{*}(q) = (a^{*},-b^{*})$ or, by a change of argument
\begin{equation}
\label{conjregidfunc}
I(q^{*}) = (a^{*},-b^{*})\mbox{\femsteg , $q^{*}=(a^{*},-b^{*})$}
\end{equation}
We regard this as the conjugate regular identity function.
In the same way, by applying the other elementary bicomplex functions
of Chapter 3 to the argument $q^{*}=(a^{*},-b^{*})$
we derive their conjugate regular counterparts, for instance:
\begin{eqnarray}
\label{conjregcosh}
\cosh q^{*}&=&(\cosh a^{*} \cdot \cos b^{*}, -\sinh a^{*} \cdot \sin b^{*}) \mbox{\femsteg , $q^{*}=(a^{*},-b^{*})$} \\ [2mm]
\label{conjregsinh}
\sinh q^{*}&=&(\sinh a^{*} \cdot \cos b^{*}, -\cosh a^{*} \cdot \sin b^{*}) \\ [2mm]
\label{conjregcos}
\cos q^{*}&=& (\cos a^{*} \cdot \cosh b^{*}, \sin a^{*} \cdot \sinh b^{*}) \\ [2mm]
\label{conjregsin}
\sin q^{*}&=& (\sin a^{*} \cdot \cosh b^{*}, -\cos a^{*} \cdot \sinh b^{*})
\end{eqnarray}
The component functions of the complex pairs at the
right-hand sides are conjugate holomorphic in both $a$ and $b$.

Analogously to the regular functions the conjugate regular ones
form an algebra under bicomplex operations.

\subsection{New classes of differentiable functions}
By considering the derivatives of the regular and conjugate
regular functions and of similar functions we shall obtain 
altogether eight classes of holomorphic
bicomplex functions.

Several authors have studied the derivative of regular
functions in the context of quaternionic function theory (see e.g.~\cite{deavours, sudbery}).
A general observation has been that it is unproductive to define 
the derivative of a regular quaternionic
function $\psi$ at the point
$p$ by applying definition (2.\ref{qinverse}) for the quaternionic inverse and taking the limit
\[
\lim_{\triangle p \rightarrow 0} \: [(\psi(p + \triangle p)-\psi(p))\times ((\triangle p)\!\uparrow\! -1)]
\]
because for $p=x+yi+zj+uk$ and $\triangle p=\triangle x+\triangle yi+\triangle zj+\triangle uk$
the class of differentiable functions will become too
restricted to be of interest. 
(The class will contain only first-degree polynomials.
Nothing is gained by using left-division instead of right-division
in the above expression.) The remedy is to employ bicomplex
division and to insist on taking the limit of a fraction
that is regular in Fueter's sense. The latter is achieved by replacing
$p$ and $\triangle p$ by the regular
\begin{eqnarray*}
& &q=(a,b^{*})=(x+iy, z-iu) \\
& &\triangle q=(\triangle a,\triangle b^{*})=(\triangle x+i\triangle y, \triangle z-i\triangle u)
\end{eqnarray*}
Hence, the derivative of a holomorphic regular function
is introduced as follows.
\begin{definition}
{\rm
Let $\psi$ be a regular bicomplex function whose domain of definition
contains a neighbourhood of the point $q=(a,b^{*})$. The derivative
of $\psi$ at the point $q$ is defined by the equation
\[
\psi\,'(q) = \lim_{\stackrel{{\scriptstyle \triangle q \rightarrow 0}}{C\!N(\triangle q)\neq 0}}\frac{\psi(q+\triangle q)-\psi(q)}{\triangle q}
\]
provided this limit exists, when it is computed for nonsingular 
$\triangle q=(\triangle a,\triangle b^{*})$.  $\Box$
}
\end{definition}

\bigskip

\noindent Apart from the function argument this definition
is the same as Definition 4.1 for the derivative
of an "ordinary" bicomplex function of the variable $p=(a,b)$. 
This means that if one
makes the appropriate change of argument, the expected formulas
for computing
the derivative of a sum, product and quotient of bicomplex
functions as well as the chain rule are retained for the
regular holomorphic functions. Furthermore, if $\psi$
has the complex pair representation 
$\psi(q)=(\phi_{1}(a,b^{*}),\phi_{2}(a,b^{*}))$
its derivative is also given by the following
equivalent formulas, analogous to
(4.\ref{bicderc21}) and (4.\ref{bicderc22}):
\begin{eqnarray}
\label{regderc21}
\frac{d\psi}{dq} &=& \left(\frac{\partial \phi_{1}}{\partial a}\;, \frac{\partial \phi_{2}}{\partial a} \right) \mbox{\femsteg , $\; q=(a,b^{*})$} \\ [3mm]
\label{regderc22}
\frac{d\psi}{dq} &=& \left(\frac{\partial \phi_{2}}{\partial b^{*}}\;, -\frac{\partial \phi_{1}}{\partial b^{*}} \right)
\end{eqnarray}
Application of either formula to the regular exponential function 
$e^{q}=(e^{a}\!\cdot \cos b^{*}, e^{a}\!\cdot \sin b^{*})$, for example, yields
the same result or 
\[
\frac{de^{q}}{dq}=e^{q}
\]
The derivatives of the other elementary regular functions are
likewise found with the normal formulas.

\medskip

\noindent {\bf Remark.} In the application of (\ref{regderc22})
it is helpful to remember (\ref{derdefconjfconja0})
\[
\frac{df^{*}}{da^{*}}=\left(\frac{df}{da}\right)^{*}
\]
satisfied by a conjugate holomorphic complex function $f^{*}$.
It enables one to compute the derivative of e.g. $\cos b^{*}$ according to:
\[
\frac{d \cos b^{*}}{d b^{*}} = \frac{d (\cos b)^{*}}{d b^{*}} = \left(\frac{d \cos b}{d b}\right)^{*} = (-\sin b)^{*} = -\sin b^{*}
\]
$\Box$

\bigskip

\noindent The equivalence of (\ref{regderc21}) and (\ref{regderc22}) presupposes the
simultaneous validity of the equations
\begin{eqnarray}
\label{RCR1}
\frac{\partial \phi_{1}}{\partial a} &=& \frac{\partial \phi_{2}}{\partial b^{*}} \\ [3mm]
\label{RCR2}
\frac{\partial \phi_{2}}{\partial a} &=& -\frac{\partial \phi_{1}}{\partial b^{*}}
\end{eqnarray}
They express the differentiability condition of the
regular functions and are termed the {\em regular
Cauchy-Riemann equations}.

\medskip

\noindent {\bf Example.} For $q=(a,b^{*})$ the function
\[
\psi(q)=(a^{2},(b^{*})^{2}) 
\]
is regular but not holomorphic, because it does not
fulfill the above equations. In contrast, the square function
\[
\psi(q)=q^{2}=(a^{2}-(b^{*})^{2},2ab^{*}) 
\]
is both regular and holomorphic. $\Box$

\bigskip

\noindent It is instructive to rewrite (\ref{regderc21}) and (\ref{regderc22})
in ${\rm R^{4}}$ assuming that 
$\phi_{1} = \psi_{1}+i\psi_{2}$ and $\phi_{2} = \psi_{3}+i\psi_{4}$.
Using formulas (\ref{derdeffa1})--(\ref{derdeffa2})
and (\ref{derdeffconja1})--(\ref{derdeffconja2}) of Section 6.2
we get two formulas for each of the derivatives
$\frac{\partial \phi_{1}}{\partial a}$, $\frac{\partial \phi_{2}}{\partial a}$, $\frac{\partial \phi_{1}}{\partial b^{*}}$ and $\frac{\partial \phi_{2}}{\partial b^{*}}$.
Substituted into (\ref{regderc21}) and (\ref{regderc22}) they furnish 
four ${\rm R^{4}}$-representations of $\frac{d\psi}{dq}$:
\begin{eqnarray}
\label{regderr41}
\frac{d\psi}{dq}&=&(\partial_{x}\psi_{1}+i\partial_{x}\psi_{2}\:,\:\partial_{x}\psi_{3}+i\partial_{x}\psi_{4}) \\ [2mm]
\label{regderr42}
\frac{d\psi}{dq}&=&(\partial_{y}\psi_{2}-i\partial_{y}\psi_{1}\:,\:\partial_{y}\psi_{4}-i\partial_{y}\psi_{3}) \\ [2mm]
\label{regderr43}
\frac{d\psi}{dq}&=&(\partial_{z}\psi_{3}+i\partial_{z}\psi_{4}\:,\:-\partial_{z}\psi_{1}-i\partial_{z}\psi_{2}) \\ [2mm]
\label{regderr44}
\frac{d\psi}{dq}&=&(-\partial_{u}\psi_{4}+i\partial_{u}\psi_{3}\:,\:\partial_{u}\psi_{2}-i\partial_{u}\psi_{1})
\end{eqnarray}
The first two of these were obtained from (\ref{regderc21}), the last two
from (\ref{regderc22}).

\noindent {\bf Example.} For $q=(a,b^{*})=(x+iy,z-iu)$
the ${\rm R^{4}}$-representation
of the regular square function $\psi(q)=q^{2}$ is
\[
\psi(q)\; =\; (x^{2}-y^{2}-z^{2}+u^{2}+i(2xy+2zu)\, , \, 2xz+2yu+i(2yz-2xu))
\]
Each of (\ref{regderr41})--(\ref{regderr44}) yields the same derivative or
\[
\frac{d\psi}{dq}=(2x+i2y,2z-i2u)
\]
This means that $\frac{d\psi}{dq}=2q$, as expected. $\Box$

\bigskip

\noindent The central concepts pertaining to the differentiation
of the conjugate regular functions are analogous.
It was shown in the previous section that the application of bicomplex functions
to the argument $q^{*}=(a^{*},-b^{*})$ yields the corresponding
conjugate regular functions, as illustrated by 
(\ref{conjregidfunc})--(\ref{conjregsin}).
The same argument is used in the definition of the derivative
of such a function.

Let $\psi(q^{*})=(\phi_{1}(a^{*},-b^{*}),\phi_{2}(a^{*},-b^{*}))$, where $q^{*}=(a^{*},-b^{*})$,
be a conjugate regular function. Its derivative at the point $q^{*}=(a^{*},-b^{*})$
is defined by the equation
\[
\psi\,'(q^{*}) = \lim_{\stackrel{{\scriptstyle \triangle q^{*} \rightarrow 0}}{C\!N(\triangle q^{*})\neq 0}}\frac{\psi(q^{*}+\triangle q^{*})-\psi(q^{*})}{\triangle q^{*}}
\]
provided this limit exists, when it is computed for nonsingular 
$\triangle q^{*}=(\triangle a^{*},-\triangle b^{*})$.
Analysis of this derivative
leads to the equivalent complex pair representations
\begin{eqnarray}
\label{conjregderc21}
\frac{d\psi}{dq^{*}} &=& \left(\frac{\partial \phi_{1}}{\partial a^{*}}\,, \frac{\partial \phi_{2}}{\partial a^{*}} \right) \mbox{\femsteg , $\; q^{*}=(a^{*},-b^{*})$} \\ [2mm]
\label{conjregderc22}
\frac{d\psi}{dq^{*}} &=& \left(-\frac{\partial \phi_{2}}{\partial b^{*}}\,, \frac{\partial \phi_{1}}{\partial b^{*}} \right)
\end{eqnarray}
The equivalence of (\ref{conjregderc21}) and (\ref{conjregderc22}) presupposes the validity of the
{\em conjugate regular Cauchy-Riemann equations}
\begin{eqnarray}
\label{CRCR1}
\frac{\partial \phi_{1}}{\partial a^{*}} &=& -\frac{\partial \phi_{2}}{\partial b^{*}} \\ [3mm]
\label{CRCR2}
\frac{\partial \phi_{2}}{\partial a^{*}} &=& \frac{\partial \phi_{1}}{\partial b^{*}}
\end{eqnarray}
Let us summarize our present findings concerning differentiable functions.
In Section 6.2 we saw that in the complex plane a holomorphic
or conjugate holomorphic function has an argument of the type
$a=(x,y)$ or $a=(x,-y)$ when expressed as a real pair. 
In Chapter 4 and this section we have obtained holomorphic
and conjugate holomorphic functions in the bicomplex space
by applying the elementary functions of Chapter 3 to the arguments
$p=(a,b)$, $q=(a,b^{*})$, and $q^{*}=(a^{*},-b^{*})$.
Looking at the form of these arguments it is hard to avoid
the thought that there must be differentiable functions
associated with each of the following arguments:
\begin{eqnarray}
\label{argc2}
&&p=(a,b) \;\;\;\;\;\mbox{\tvasteg , \tvasteg} q=(a,b^{*})  \;\;\;\;\mbox{\tvasteg , \tvasteg} r=(a^{*},b)  \;\;\;\;\;\mbox{\tvasteg , \tvasteg} s=(a^{*},b^{*}) \\ [2mm]
\label{biconjargc2}
&&p^{\smallcup}=(a,-b) \mbox{\tvasteg , \tvasteg} q^{\smallcup}=(a,-b^{*})  \mbox{\tvasteg , \tvasteg} r^{\smallcup}=(a^{*},-b)  \mbox{\tvasteg , \tvasteg} s^{\smallcup}=(a^{*},-b^{*})
\end{eqnarray}
Each argument of the second row is the bicomplex conjugate
of the corresponding argument of the first row. 

Bicomplex functions applied to the
arguments (\ref{argc2})--(\ref{biconjargc2}) are evidently of two main types:
\begin{eqnarray}
\label{bicfunctype1}
\psi(t)&=&(\phi_{1}(a^{\smallbox},b^{\smalldiamond}),\phi_{2}(a^{\smallbox},b^{\smalldiamond})) \mbox{\femsteg \tresteg , $t=(a^{\smallbox},b^{\smalldiamond})$} \\ [2mm]
\label{bicfunctype2}
\psi(t^{\smallcup})&=&(\phi_{1}(a^{\smallbox},-b^{\smalldiamond}),\phi_{2}(a^{\smallbox},-b^{\smalldiamond})) \:\mbox{\femsteg , $t^{\smallcup}=(a^{\smallbox},-b^{\smalldiamond})$}
\end{eqnarray}
The variable $t$ should be replaced by $p$, $q$, $r$ or $s$ and
either of the symbols \raisebox{.6ex}{$\smallbox$} and \raisebox{.6ex}{$\smalldiamond$}
by a blank or a \raisebox{.6ex}{$*$}.
With the derivative 
defined as the limit of a fraction
containing the appropriate argument the
functions (\ref{bicfunctype1}) are holomorphic if they fulfill complexified
CR-equations of the form
\begin{eqnarray}
\label{type1CR1}
\frac{\partial \phi_{1}}{\partial a^{\smallbox}} &=& \frac{\partial \phi_{2}}{\partial b^{\smalldiamond}} \\ [3mm]
\label{type1CR2}
\frac{\partial \phi_{2}}{\partial a^{\smallbox}} &=& -\frac{\partial \phi_{1}}{\partial b^{\smalldiamond}}
\end{eqnarray}
while the functions (\ref{bicfunctype2}) are conjugate holomorphic if they fulfill
\begin{eqnarray}
\label{type2CR1}
\frac{\partial \phi_{1}}{\partial a^{\smallbox}} &=& -\frac{\partial \phi_{2}}{\partial b^{\smalldiamond}} \\ [3mm]
\label{type2CR2}
\frac{\partial \phi_{2}}{\partial a^{\smallbox}} &=& \frac{\partial \phi_{1}}{\partial b^{\smalldiamond}}
\end{eqnarray}
Equations (\ref{type1CR1})--(\ref{type1CR2}) 
correspond to the CR-equations (\ref{CReqsr2}), 
equations (\ref{type2CR1})--(\ref{type2CR2}) to the
conjugate CR-equations (\ref{ConjCReqsr2}).

\bigskip

\noindent {\bf Example.} Application of the bicomplex
cosinus function (3.\ref{biccosc2}) to $r^{\smallcup}=(a^{*},-b)$ yields
$\cos r^{\smallcup} = (\cos a^{*} \cdot \cosh b \,,\, \sin a^{*} \cdot \sinh b)$.
The component functions $\cos a^{*} \cdot \cosh b$ and 
$\sin a^{*} \cdot \sinh b$ satisfy (\ref{type2CR1})--(\ref{type2CR2}) if a 
\raisebox{.6ex}{$*$} is substituted for \raisebox{.6ex}{$\smallbox$} and
a blank for \raisebox{.6ex}{$\smalldiamond$}:

\begin{eqnarray*}
\frac{\partial \phi_{1}}{\partial a^{*}} &=& -\frac{\partial \phi_{2}}{\partial b} \\ [3mm]
\frac{\partial \phi_{2}}{\partial a^{*}} &=& \frac{\partial \phi_{1}}{\partial b}
\end{eqnarray*}
$\Box$

\subsection{Bicomplex solutions of the Laplace equation}
A final remark about the solutions of the Laplace equation is in order.

Let the function $\psi$ be given by (\ref{bicfunctype1}) so that
$\phi_{1}(a^{\smallbox},b^{\smalldiamond}) = \psi_{1}(x,y,z,u)+i\psi_{2}(x,y,z,u)$ and
$\phi_{2}(a^{\smallbox},b^{\smalldiamond}) = \psi_{3}(x,y,z,u)+i\psi_{4}(x,y,z,u)$. 
(The corresponding treatment of a function of type (\ref{bicfunctype2})
is analogous.)
If $\phi_{1}$ and $\phi_{2}$ are holomorphic or
conjugate holomorphic in $a$ and $b$, their
component functions fulfill the two-dimensional
Laplace equation:
\begin{eqnarray*}
\partial_{x}^{2}\psi_{i} + \partial_{y}^{2}\psi_{i} &=& 0 \mbox{\femsteg , $i=1,2,3,4$}\\ [1mm]
\partial_{z}^{2}\psi_{i} + \partial_{u}^{2}\psi_{i} &=& 0
\end{eqnarray*}
It follows that $\psi$ fulfills the 
four-dimensional Laplace equation 
$\partial_{x}^{2}\psi + \partial_{y}^{2}\psi + \partial_{z}^{2}\psi + \partial_{u}^{2}\psi = 0$.

Application of the differential operator $\frac{\partial}{\partial a^{\smallbox}}$ to one
of the Cauchy-Riemann equations (\ref{type1CR1})--(\ref{type1CR2}) and 
$\frac{\partial}{\partial b^{\smalldiamond}}$ to the other yields
after computing the sum or difference of the results either
\[
\frac{\partial^{2}\phi_{1}}{\partial a^{\smallbox 2}}+\frac{\partial^{2}\phi_{1}}{\partial b^{\smalldiamond 2}} = 0
\]
or
\[
\frac{\partial^{2}\phi_{2}}{\partial a^{\smallbox 2}}+\frac{\partial^{2}\phi_{2}}{\partial b^{\smalldiamond 2}} = 0
\]
The complex pair components of all bicomplex
holomorphic functions thus satisfy the two-dimensional,
complexified Laplace equation.

In fact, an alternative approach to the understanding
of differentiable bicomplex functions is to take the
previous two equations as starting point and combine them into
\begin{equation}
\label{Lapleqc2}
\left(\frac{\partial^{2}}{\partial a^{\smallbox 2}}+\frac{\partial^{2}}{\partial b^{\smalldiamond 2}}\,,\,0\right)\odot (\phi_{1},\phi_{2})\;=\;0
\end{equation}
This is a generalization of equation (\ref{Lapleqr2}). The operator
to the left admits the factorization
\[
\left(\frac{\partial^{2}}{\partial a^{\smallbox 2}}+\frac{\partial^{2}}{\partial b^{\smalldiamond 2}}\,,\,0\right) \;=\;\left(\frac{\partial}{\partial a^{\smallbox}}\,,\,\frac{\partial}{\partial b^{\smalldiamond}}\right)\odot \left(\frac{\partial}{\partial a^{\smallbox}}\,,\,-\frac{\partial}{\partial b^{\smalldiamond}}\right)
\]
which means that the solutions of
\begin{eqnarray*}
 \left(\frac{\partial}{\partial a^{\smallbox}}\,,\,\frac{\partial}{\partial b^{\smalldiamond}}\right)\odot (\phi_{1},\phi_{2})\;&=&\;0 \\[2mm]
 \left(\frac{\partial}{\partial a^{\smallbox}}\,,\,-\frac{\partial}{\partial b^{\smalldiamond}}\right)\odot (\phi_{1},\phi_{2})\;&=&\;0
\end{eqnarray*}
are solutions of (\ref{Lapleqc2}) as well. Because these two equations are equivalent
to (\ref{type1CR1})--(\ref{type1CR2}) and (\ref{type2CR1})--(\ref{type2CR2}), respectively,
we see as in our treatment of the complex holomorphic
and conjugate holomorphic functions in Section 6.1, that the corresponding
bicomplex functions can be obtained by factorizing the appropriate
version of the Laplace equation.

\newpage
\section{Conclusions}
In the course of this work we have frequently seen that a concept introduced or a result
obtained has been the generalization of the corresponding concept or result in
complex analysis. No doubt there are still plenty
of such generalizations to be made, notably with regard to power series,
Laurent series and analytic continuation.

A somewhat separate matter is the development of a theory
of {\em tricomplex functions} of the type $\Omega =(\Psi_{1},\Psi_{2})$,
where $\Psi_{1}$ and $\Psi_{2}$ are bicomplex functions of two
bicomplex variables. In forming differentiable tricomplex functions  one
has at one's disposal all the holomorphic
and conjugate holomorphic functions of Section 6.5.

A corner stone of our work has been the application of complex pairs to
represent quaternions and bicomplex numbers. With respect to the former
essentially the same idea is found in the book~\cite{maclane} by MacLane and Birkhoff, but
instead of rendering a quaternion as a complex pair $(a,b)$ they
denote it by $a+bj$, which in our notation becomes 
$(a,0)+(b,0)\times j$. In~\cite{price} Price employs the notation $a+i_{2}b$
for the bicomplex number that we would write $(a,0)+ j \odot (b,0)$.
All these notations seem rather cumbersome compared to the concise pair notation.

In this context it may be noted
that the practice of denoting complex numbers by pairs of real numbers
is due to Hamilton~\cite{hamilton1}. 
His objective was to employ these 
algebraic couples or number-couples, as he called them, to remove the
"metaphysical stumbling blocks" that he felt was connected
with the imaginary unit $\sqrt{-1}$. An account of Hamilton's thoughts
on algebraic couples is found in~\cite{kline} , the reading of which inspired
the writing of this treatise. 

\subsection*{Acknowledgements}
I am indebted to Prof. Raimo Lehti for lucid lectures on series and function
theory at the Helsinki University of Technology~\cite{lehti} and for valuable
comments on a previous version of this treatise.
The financial support of the Ella and Georg Ehrnrooth Foundation is
also gratefully acknowledged.

\newpage
\setcounter{equation}{0}
\setcounter{section}{5}
\setcounter{theorem}{3}
\addcontentsline{toc}{section}{Appendix}
\subsection*{Appendix}

\section*{The generalized Green theorem}
We shall give a proof of the generalized
Green theorem as stated in Chapter 5.
\begin{theorem}
Let $\phi_{1}$ and $\phi_{2}$ be two complex functions
of $a$ and $b$ such that they are holomorphic in a domain
$G\subseteq {\rm B}$. Also let $S(\Gamma)$ be an orientable surface
in $G$ bounded by the closed curve $\Gamma$. Then
\begin{equation}
\label{greendef}
\int_{\Gamma} \phi_{1}\cdot da + \phi_{2}\cdot db = \int_{S(\Gamma)} \left(\frac{\partial \phi_{2}}{\partial a} - \frac{\partial \phi_{1}}{\partial b}\right)\cdot da\cdot db
\end{equation}
\end{theorem}
\noindent {\it Proof.} Let the ${\rm R}^{4}$-representations
of $\phi_{1}$ and $\phi_{2}$ be
\begin{eqnarray}
\label{phidef}
\phi_{1}(a,b) &=& \psi_{1}(x,y,z,u)+i\,\psi_{2}(x,y,z,u) \\ [1mm]
\phi_{2}(a,b) &=& \psi_{3}(x,y,z,u)+i\,\psi_{4}(x,y,z,u) \nonumber \\ [2mm]
a=x+iy &,& b=z+iu \nonumber
\end{eqnarray}
Because $\phi_{1}$ and $\phi_{2}$ are holomorphic in $G$ we get
on application of (3.\ref{dthetadainx}) the following derivatives:
\begin{eqnarray}
\label{deriv1}
\frac{\partial\phi_{2}}{\partial a} &=& \partial_{x}\psi_{3}+i\partial_{x}\psi_{4} \\ [3mm]
\label{deriv2}
\frac{\partial\phi_{1}}{\partial b} &=& \partial_{z}\psi_{1}+i\partial_{z}\psi_{2}
\end{eqnarray}
Inserting (\ref{phidef}), (\ref{deriv1}) and (\ref{deriv2}) into (\ref{greendef}), performing the complex multiplications
and demanding the equality of real and imaginary parts change
our proof obligation into:
\begin{eqnarray}
\label{integral1}
& &\int_{\Gamma} (\psi_{1}\, dx - \psi_{2}\, dy + \psi_{3}\, dz - \psi_{4}\, du) = \\
& &\mbox{\tresteg}\int_{S(\Gamma)} [(\partial_{x}\psi_{3}-\partial_{z}\psi_{1})\, dx dz - (\partial_{x}\psi_{3}-\partial_{z}\psi_{1})\, dy du - \nonumber \\
& &\mbox{\tresteg \femsteg} (\partial_{x}\psi_{4}-\partial_{z}\psi_{2})\, dx du - (\partial_{x}\psi_{4}-\partial_{z}\psi_{2})\, dy dz)] \nonumber \\ [3mm]
\label{integral2}
& &\int_{\Gamma} (\psi_{2}\, dx + \psi_{1}\, dy + \psi_{4}\, dz + \psi_{3}\, du) = \\
& &\mbox{\tresteg}\int_{S(\Gamma)} [(\partial_{x}\psi_{3}-\partial_{z}\psi_{1})\, dx du + (\partial_{x}\psi_{3}-\partial_{z}\psi_{1})\, dy dz + \nonumber \\
& &\mbox{\tresteg \femsteg}(\partial_{x}\psi_{4}-\partial_{z}\psi_{2})\, dx dz - (\partial_{x}\psi_{4}-\partial_{z}\psi_{2})\, dy du)] \nonumber 
\end{eqnarray}
The line and surface integrals of these formulas are of the type
found at the right-hand sides of (5.\ref{lineintr4}) and (5.\ref{surfintr4}).

At this point we invoke the {\em generalized Stokes theorem} in
${\rm R}^{4}$~\cite{flanders, teichmann}. It states that four 
real-valued, continuously differentiable functions $A_{1}, A_{2}, A_{3}, A_{4}$ 
of the variables $x,y,z,u$ satisfy
\begin{eqnarray}
\label{stokes}
& &\int_{\Gamma} (A_{1}\, dx + A_{2}\, dy + A_{3}\, dz + A_{4}\, du) = \\
& &\mbox{\tresteg}\int_{S(\Gamma)} (A_{12}\, dx dy + A_{13}\, dx dz + A_{14}\, dx du + A_{23}\, dy dz + A_{24}\, dy du + A_{34}\, dz du) \nonumber \\ [3mm]
& &\mbox{where $S(\Gamma)$ is an orientable surface in ${\rm R}^{4}$ bounded by $\Gamma$ and} \nonumber \\ [3mm]
& &\begin{array}{ccccc}
 A_{12}=\partial_{x}A_{2}-\partial_{y}A_{1}& , &A_{13}=\partial_{x}A_{3}-\partial_{z}A_{1}& , & A_{14}=\partial_{x}A_{4}-\partial_{u}A_{1} \nonumber \\[2mm]
  A_{23}=\partial_{y}A_{3}-\partial_{z}A_{2}& , &A_{24}=\partial_{y}A_{4}-\partial_{u}A_{2}& , & A_{34}=\partial_{z}A_{4}-\partial_{u}A_{3}
 \end{array}
\end{eqnarray}

\medskip

\noindent Our intention is to show that (\ref{stokes}) reduces to (\ref{integral1}) or (\ref{integral2}) when
the $A_{k}$-functions are chosen properly. Focusing on formula (\ref{integral1}) first,
we compare its left-hand side with the left-hand side of (\ref{stokes}) and identify
\[
A_{1}=\psi_{1} \;\;,\;\; A_{2}=-\psi_{2} \;\;,\;\; A_{3}=\psi_{3} \;\;,\;\; A_{4}=-\psi_{4}
\]
The functions $\phi_{1}$ and $\phi_{2}$ are holomorphic in $a$ and $b$,
hence the $\psi_{k}$-functions are continuously differentiable and
satisfy the CR-equations (3.\ref{CRequations}) in the following way:
\[
\begin{array}{ccc}
\partial_{x}\psi_{1}=\partial_{y}\psi_{2} & , & \partial_{x}\psi_{2}=-\partial_{y}\psi_{1} \\ [3mm]
\partial_{z}\psi_{1}=\partial_{u}\psi_{2} & , & \partial_{z}\psi_{2}=-\partial_{u}\psi_{1} \\ [5mm]
\partial_{x}\psi_{3}=\partial_{y}\psi_{4} & , & \partial_{x}\psi_{4}=-\partial_{y}\psi_{3} \\ [3mm]
\partial_{z}\psi_{3}=\partial_{u}\psi_{4} & , & \partial_{z}\psi_{4}=-\partial_{u}\psi_{3}
\end{array}
\]
We thus get
\[
\begin{array}{ccccl}
A_{12} &=& -\partial_{x}\psi_{2}-\partial_{y}\psi_{1} &=& 0 \\ [2mm]
A_{13} &=& \partial_{x}\psi_{3}-\partial_{z}\psi_{1} & &  \\ [2mm]
A_{14} &=& -\partial_{x}\psi_{4}-\partial_{u}\psi_{1} &=& -\partial_{x}\psi_{4}+\partial_{z}\psi_{2} \\ [2mm]
A_{23} &=& \partial_{y}\psi_{3}+\partial_{z}\psi_{2} &=& -(\partial_{x}\psi_{4}-\partial_{z}\psi_{2}) \\ [2mm]
A_{24} &=& -\partial_{y}\psi_{4}+\partial_{u}\psi_{2} &=& -\partial_{x}\psi_{3}+\partial_{z}\psi_{1} \\ [2mm]
A_{34} &=& -\partial_{z}\psi_{4}-\partial_{u}\psi_{3} &=& 0
\end{array}
\]

\medskip

\noindent These identities make the right-hand side of (\ref{stokes}) equal to the right-hand side of (\ref{integral1}), as required.

Comparison of the left-hand side of (\ref{stokes}) with the left-hand side of (\ref{integral2}), in turn, tells us to choose
\[
A_{1}=\psi_{2} \;\;,\;\; A_{2}=\psi_{1} \;\;,\;\; A_{3}=\psi_{4} \;\;,\;\; A_{4}=\psi_{3}
\]
It is as straightforward as above to show that evaluated for this choice
the expressions $A_{ij}$ make 
the right-hand side of (\ref{stokes}) equal to the right-hand side of (\ref{integral2}), which concludes the proof.

\medskip

\noindent $\Box$

\end{document}